\documentclass[11pt,a4paper]{article}

\usepackage[utf8]{inputenc}
\usepackage[OT2,T1]{fontenc}

\usepackage[a4paper]{geometry}
\usepackage[french,english]{babel}

\usepackage{fancyhdr}
\usepackage[final]{graphicx}

\usepackage[figurewithin=none]{caption}
\usepackage{zref}
\usepackage{dsfont}
\usepackage{enumitem}
\usepackage{ifthen}
\usepackage{verbatim}

\usepackage[nodisplayskipstretch]{setspace}
\usepackage{amsmath,amssymb,makeidx,amsthm}
\usepackage{stmaryrd}
\usepackage{mathrsfs}
\usepackage[all]{xy} 
\usepackage{multirow} 
\usepackage{array}
\usepackage{mathtools}
\usepackage{tikz-cd}

\usepackage{url}
\usepackage[hypertexnames=false]{hyperref}

\newcommand{\Z}{{\mathbb Z}}
\newcommand{\Q}{{\mathbb Q}}
\renewcommand{\C}{{\mathbb C}}

\newcommand{\F}{{\mathbb F}}
\newcommand{\rar}{\rightarrow}

\renewcommand{\P}{\mathscr{P}}
\newcommand{\Ell}{\ensuremath{\mathbf{Ell}}}
\newcommand{\Sch}{\ensuremath{\mathbf{Sch}}}
\newcommand{\GpSch}{\ensuremath{\mathbf{GpSch}}}
\newcommand{\Set}{\ensuremath{\mathbf{Set}}}
\newcommand{\Gp}{\ensuremath{\mathbf{Gp}}}

\newcommand{\M}{\ensuremath{\mathscr{M}}}
\newcommand{\MP}{\ensuremath{\M(\P)}}
\newcommand{\MPp}[1]{\ensuremath{\M(#1)}}
\newcommand{\CMP}{\ensuremath{\overline{\M}(\P)}}
\newcommand{\CMPp}[1]{\ensuremath{\overline{\M}(#1)}}
\newcommand{\CP}{\ensuremath{C_{\P}}}
\newcommand{\CPp}[1]{\ensuremath{C_{#1}}}
\newcommand{\mfk}[1]{\mathfrak{#1}}
\newcommand{\cE}{\mathcal{E}}
\newcommand{\bA}{\mathbb{A}}
\newcommand{\bP}{\mathbb{P}}

\renewcommand{\emptyset}{\varnothing}
\newcommand{\mrm}[1]{\mathrm{#1}}
\newcommand{\OO}{\mathcal{O}}
\newcommand{\Qbar}{\overline{\mathbb{Q}}}
\newcommand{\Fr}{\mrm{Frob}}
\renewcommand{\tilde}[1]{\widetilde{#1}}
\newcommand{\Aut}[1]{\ensuremath{\operatorname{Aut}(#1)}}

\newcommand{\GL}[1]{\ensuremath{\operatorname{GL}_2(#1)}}
\newcommand{\SL}[1]{\ensuremath{\operatorname{SL}_2(#1)}}

\newcommand{\PSL}[1]{\ensuremath{\operatorname{PSL}_2(#1)}}

\DeclareMathOperator{\Sp}{Spec}

\usepackage{calc}

\makeatletter
\let\oldlabel\label
\renewcommand{\label}[1]{%
  \zref@labelbylist{#1}{special}
  \oldlabel{#1}
}

\zref@newlist{special}
\zref@newprop{section}{\arabic{section}.\arabic{subsection}}
\zref@addprop{special}{section}
\newcounter{propo}
\renewcommand{\thepropo}{\thesubsection.\arabic{propo}}
\makeatletter \@addtoreset{propo}{subsection}\makeatother

\makeatletter \@addtoreset{rema}{subsection}\makeatother

\newcounter{mytheo}
\renewcommand{\themytheo}{\arabic{mytheo}}

\newcounter{mycor}
\renewcommand{\themycor}{\arabic{mycor}}
\makeatletter \@addtoreset{mycor}{mytheo}\makeatother

\newcommand{\questin}[2][None]{\ifthenelse{\equal{#1}{None}}{}{\label{#1}} \smallskip \noindent \textbf{Question \;} {\it #2}\smallskip}
\newcommand{\conjintro}[2][None]{\ifthenelse{\equal{#1}{None}}{}{\label{#1}} \smallskip \noindent \textbf{Conjecture \;} {\it #2}\smallskip}
\newcommand{\theoin}[2][None]{\ifthenelse{\equal{#1}{None}}{}{\label{#1}} \smallskip \noindent \textbf{Theorem \;} {\it #2}\smallskip}

\newcommand{\lem}[2][None]{\refstepcounter{propo}\ifthenelse{\equal{#1}{None}}{}{\label{#1}} \smallskip \noindent \textbf{Lemma \thepropo\;} {\it #2}\smallskip}
\newcommand{\nott}[1]{\smallskip \noindent \textbf{Notation \;} {\it #1} \smallskip}
\newcommand{\theo}[2][None]{\refstepcounter{propo}\ifthenelse{\equal{#1}{None}}{}{\label{#1}} \smallskip \noindent \textbf{Theorem \thepropo\;} {\it #2} \smallskip}
\newcommand{\theoi}[2][None]{\refstepcounter{theop}\ifthenelse{\equal{#1}{None}}{}{\label{#1}}\smallskip \noindent \textbf{Theorem \thetheop\;} {\it #2} \smallskip}
\newcommand{\cori}[2][None]{\refstepcounter{theop}\ifthenelse{\equal{#1}{None}}{}{\label{#1}}\smallskip \noindent \textbf{Corollary \thetheop\;} {\it #2} \smallskip}

\newcommand{\prop}[2][None]{\refstepcounter{propo}\ifthenelse{\equal{#1}{None}}{}{\label{#1}} \smallskip \noindent \textbf{Proposition \thepropo\;} {\it #2} \smallskip}

\newcommand{\cor}[2][None]{\refstepcounter{propo}\ifthenelse{\equal{#1}{None}}{}{\label{#1}} \smallskip \noindent \textbf{Corollary \thepropo\;} {\it #2} \smallskip}

\newcommand{\defi}[2][None]{\refstepcounter{propo}\ifthenelse{\equal{#1}{None}}{}{\label{#1}} \smallskip \noindent {\bf Definition \thepropo\;} {#2} \smallskip}

\newcommand{\rem}[2][None]{\refstepcounter{propo}\ifthenelse{\equal{#1}{None}}{}{\label{#1}}\smallskip \noindent {\bf Remark \thepropo\;} {#2} \smallskip}
\newcommand{\rems}[1]{\refstepcounter{propo}\smallskip \noindent {\bf Remarks \thepropo\;} {#1} \smallskip}

\newcommand{\demo}[1]{\noindent {\it Proof.\;--\;} #1\hfill$\Box$ \smallskip}

\newcommand{\qnn}[2][None]{\refstepcounter{propo}\ifthenelse{\equal{#1}{None}}{}{\label{#1}} \smallskip \noindent \textbf{Question \thepropo\;} {\it #2} \smallskip}

\newcommand{\mytheo}[2][None]{\refstepcounter{mytheo}\ifthenelse{\equal{#1}{None}}{}{\label{#1}} \smallskip \noindent \textbf{Theorem \themytheo\;} {\it #2}\smallskip}
\newcommand{\mycor}[2][None]{\refstepcounter{mycor}\ifthenelse{\equal{#1}{None}}{}{\label{#1}} \smallskip \noindent \textbf{Corollary \themycor\;} {\it #2}\smallskip}

\title{Compactified moduli spaces and Hecke correspondences for elliptic curves with a prescribed $N$-torsion scheme}
\author{Elie Studnia}

\begin{document}
\maketitle

\tableofcontents

\section*{Introduction}

In the article \cite{FreyMazur}, Mazur asks the following question: 

\questin{(Mazur, \cite{FreyMazur}) Does there exist non-isogenous elliptic curves $E, F/\Q$ and $N \geq 7$ such that the Galois-modules $E[N](\Qbar)$ and $F[N](\Qbar)$ are isomorphic? }

The abelian groups $E[N](\Qbar)$ and $F[N](\Qbar)$ are both isomorphic to $(\Z/N\Z)^{\oplus 2}$, but they are also endowed with a natural perfect bilinear alternating pairing into the set $\mu_N(\Qbar)$ of $N$-th roots of unity, the \emph{Weil pairing}. A refinement suggested by Mazur is to ask for a \emph{symplectic} isomorphism $E[N](\Qbar) \rar F[N](\Qbar)$, that is, an isomorphism respecting the Weil pairing.  

\medskip

The first answer to Mazur's question comes in the article \cite{KO} by Kraus and Oesterl\'e: they show that the elliptic curves with equations $y^2=x^3+7x^2+28$ and $y^2=x^3+x^2-x+3$ (with respective conductors $7448=2^3 \cdot 7^2\cdot 19$ and $152=2^3 \cdot 19$, hence non-isogenous) have symplectically isomorphic $7$-torsion Galois-modules. \\

This work has been improved in \cite{HalbKraus}: in Proposition 6.3 of \emph{op.cit.}, the authors show that there is an infinite family of $6$-uples of pairwise non-isogenous elliptic curves $(E_i)_{1 \leq i \leq 6}$ such that for every $1 \leq i < j \leq 6$, the Galois-modules $E_i[7](\Qbar)$ and $E_j[7](\Qbar)$ are symplectically isomorphic. 

They proceed as follows: by \cite[Proposition 1]{KO}, given an elliptic curve $E/\Q$, there is a smooth affine curve $Y_E(7)$ over $\Q$ such that, for any field $K$ of characteristic zero, $Y_E(7)(K)$ is in bijection with the isomorphism classes of couples $(F,i)$, where $F/K$ is an elliptic curve, and $i: E[7]_K \rar F[7]$ is a symplectic isomorphism of $K$-group schemes. 

The authors then compute equations for a suitable compactification of $Y_E(7)$ \cite[Th\'eor\`emes 2.1, 2.2]{HalbKraus}, using the fact that it is a smooth plane quartic, isomorphic over $\C$ to the Klein quartic $X(7)$ (with projective equation $x^3y+y^3z+z^3x=0$), containing over any field where $E$ has three nontrivial $2$-torsion points $A,B,C$ the four tautological points 
\vspace{-8pt}\[(E,\mrm{id}), (E/\langle A\rangle,5q_A), (E/\langle B\rangle,5q_B), (E/\langle C\rangle,5q_C),\vspace{-8pt}\]
 where $q_P$ denotes the quotient $E \rar E/\langle P\rangle$.  

With an equation for the curve $X_E(7)$, it can become easier to show that the $7$-torsion Galois-module of a given $E'/\Q$ is symplectically isomorphic to that of $E$: in \cite[Remarque 6.1]{HalbKraus}, the authors give an example of $(E,E')$ where the symplectic isomorphism of $7$-torsion follows from the equation, while checking it by comparing their $L$-series as in \cite[Proposition 4]{KO} would require computing coefficients at all prime indices $p \leq 25 \cdot 10^9$.

Using the equation, the authors then construct, for certain elliptic curves $E$, new points on the curve $X_E(7)$ by geometric methods, which is enough to complete the proof of \cite[Proposition 6.3]{HalbKraus}. \\

Similarly, it seems known\footnote{For instance, Mazur gives without reference a slightly stronger definition in \cite[\S I.7]{Mazur-Open}.} that, for any $n \geq 3$, any $\alpha \in (\Z/n\Z)^{\times}$, and any elliptic curve $E/\Q$, there exists a smooth affine curve $Y_E^{\alpha}(n)$ over $\Q$ such that for any field $K$ of characteristic $0$, $Y_E^{\alpha}(n)(K)$ parametrizes the isomorphism classes of pairs $(F,i)$, where $F/K$ is an elliptic curve, and $i: E[n]_K \rar F[n]$ is an isomorphism of $K$-group schemes of \emph{determinant} $\alpha$, i.e. such that $\langle i(-),\,i(\ast)\rangle_{F[n]} = \langle -,\,\ast\rangle_{E[n]}^{\alpha}$, where $\langle \cdot,\,\cdot\rangle$ denotes the Weil pairing. Because we can always multiply $i$ by a constant, the curve $Y_E^{\alpha}(n)$ only depends on the class of $\alpha$ modulo squares. In particular, when $n$ is prime or twice a prime, only two values of $\alpha$ need to be investigated: $1$ and any quadratic non-residue in $(\Z/n\Z)^{\times}$ (corresponding to so-called \emph{anti-symplectic congruences} of elliptic curves). 

\medskip

There are analogues of \cite{HalbKraus} for many small values of $n$. For instance, equations for $X_E^{\alpha}(n)$ over a field of characteristic zero are known for $3 \leq n \leq 13$ by works of:

\begin{itemize}[noitemsep,label=\tiny$\bullet$,topsep=2pt]
\item Silverberg and Rubin--Silverberg \cite{Rubin-Silverberg,Silverberg-CSS} for $n \leq 5$ and $\alpha=1$, 
\item Papadopoulos \cite{Papadopoulos} and Rubin--Silverberg \cite{RS-6} for $(n,\alpha)=(6,1)$, 
\item Poonen--Schaefer--Stoll \cite{PSS} for $(n,\alpha)=(7,-1)$,
\item Fisher \cite{FisherHess,Fisher-5} for $n \leq 5$ and $\alpha=-1$,
\item Fisher \cite{Fisher-9,Fisher-711,Fisher-13} for $n=9,11,13$ and any $\alpha$,
\item Chen \cite{Chen-8,Chen-PhD} for $n \in \{6,8,10\}$ and $(n,\alpha) \in \{(12,1),(12,7)\}$,
\item Frengley \cite{Frengley-12} for $n=12$.
\end{itemize}

\medskip

For some values of $(n,\alpha)$ given above, one can use such equations to exhibit infinite families of couples $(E,F)$ of non-isogenous elliptic curves over $\Q$ such that there is an isomorphism of $\Q$-group schemes $E[n] \rar F[n]$ with determinant $\alpha$. 

But for larger $n$, it becomes easier to work with certain \emph{modular diagonal quotient surfaces} studied among others by Kani and Schanz \cite{MDQS}, coarsely parametrizing triples $(E,F,i)$, where $E,F$ are elliptic curves and $i: E[n] \rar F[n]$ is an isomorphism of group schemes with determinant $\alpha$. The equations attached to these surfaces tend to be more tractable, and the above works thus prove the following:

\theoin[theoin-1]{For any $3 \leq n \leq 13$ and every $\alpha \in (\Z/n\Z)^{\times}$, there exists a smooth integral $\Q$-scheme $X$ of dimension one such that $X(\Q)$ is infinite, non-isogenous elliptic curves $\mathcal{E},\mathcal{F}$ over $X$, and an isomorphism of finite \'etale $X$-group schemes $i: \mathcal{E}[n] \rar \mathcal{F}[n]$ with determinant $\alpha$ such that the map $(j(\mathcal{E}),j(\mathcal{F})): X \rar \mathbb{A}^2$ is non-constant.

In particular, there is an infinite family of pairs $(E,F)$ of non-isogenous elliptic curves over $\Q$ such that there is an isomorphism $i: E[n] \rar F[n]$ of group schemes with determinant $\alpha$. Furthermore, the collection of the $(j(E),j(F)) \in \Q^2$ is infinite.
}

Using modular diagonal quotient surfaces, Fisher \cite{Fisher-17} was also able to produce the first known pair of non-isogenous elliptic curves over $\Q$ with symplectically isomorphic $17$-torsion Galois-modules. \\

On the other hand, it was shown by Cremona and Freitas \cite{Cremona-Freitas} that there is no couple $(E,F)$ of non-isogenous elliptic curves over $\Q$ with conductors at most $500000$ such that, for some prime $p > 18$, $E[p](\Qbar)$ and $F[p](\Qbar)$ are isomorphic $\mrm{Gal}(\Qbar/\Q)$-modules. Moreover, the main result of \cite{BaTs} implies that infinite families of non-isogenous elliptic curves with isomorphic $p$-torsion Galois-modules (in the sense of the above Theorem) do not exist for large enough $p$. 

These results are in line with the currently open \emph{Frey--Mazur conjecture} (as formulated in e.g. \cite{Frey-Conj}\footnote{This name is also used to refer to certain stronger versions of the claim, such as \cite[Conjecture 1.1]{Fisher-17} which claims that $C_E=18$ always works in the statement following this footnote.}):

\medskip

\conjintro[FM-intro]{Let $E/\Q$ be an elliptic curve. There is a constant $C_E > 0$ such that, for any prime $p > C_E$, for any elliptic curve $F/\Q$ such that the Galois-modules $E[p](\Qbar)$ and $F[p](\Qbar)$ are isomorphic, $E$ and $F$ are isogenous.}

\medskip

This conjecture is, in a sense, inverse to what the authors we previously cited did: broadly speaking, they managed to find, for small $n$, many elliptic curves $E$ such that $X_E^{\alpha}(n)$ had a non-trivial rational point. In contrast, the conjecture states that, given an elliptic curve $E$ and a couple $(n,\alpha)$ with $n$ ``large enough'', the only rational points on $X_E^{\alpha}(n)$ are trivial, i.e. are cusps or come from elliptic curves isogenous to $E$. 

\medskip

Such a framing is of course strongly reminiscent of Mazur's foundational work on rational points of modular curves. Indeed, he shows in \cite{MazurY1} (resp. \cite{FreyMazur}) that, for $p > 13$ (resp. $p > 163$), the only rational points of the modular curve $X_1(p)$ (resp. $X_0(p)$) are cusps. 

Mazur proceeds as follows: he first identifies in \cite{MazurY1} the \emph{Eisenstein quotient} $\tilde{J}$ of the Jacobian $J_0(p)$ of $X_0(p)$, and proves (using some rather delicate algebraic geometry) that $\tilde{J}$ has rank zero. 

In \cite{FreyMazur}, he then introduces the \emph{formal immersion argument}: he shows that the morphism $X_0(p) \rar J_0(p) \rar \tilde{J}$ is a formal immersion at the cuspidal points of finite characteristic. Mazur can then argue that, in the above situation, a non-cuspidal rational point on $X_0(p)$ is never congruent to a cusp modulo some prime $\ell \nmid 2p$. In other words, any elliptic curve $E/\Q$ possessing a rational $p$-isogeny (for a prime $p > 13$) has potentially good reduction away from $2p$. While more work is needed to deduce the main result of \cite{FreyMazur}, this already implies by purely local arguments that no elliptic curve over $\Q$ has a non-trivial $p$-torsion point for $p > 13$. 

\medskip

Mazur's approach has been widely re-used to tackle other problems pertaining to algebraic points on modular curves: the rational points on the modular curve attached to the normaliser of a split Cartan subgroup \cite{BP,BPR}, bounds for the torsion subgroup of any elliptic curve over a number field of given degree \cite{Merel,Parent-Torsion}, and rational points on certain covers of the modular curve attached to the normaliser of a non-split Cartan subgroup, corresponding to an auxiliary level structure \cite{Darmon-Merel} or a smaller Galois image \cite{LFL,lombardo}. 

The formal immersion argument plays a very special role in these proofs, and applying it requires a precise understanding of the geometry of the modular curve of interest, specifically of the special fibres of its compactification over a ring of integers (rather than a number field). 

In the case of classical modular curves, this was done by finding them good definitions as moduli spaces that made sense over more general bases than fields, which was precisely described in the well-known treatises \cite{DeRa} and \cite{KM}. In the latter, the authors are able to leverage their careful study of the curve $X_1(p^r)$ in characteristic $p$ to prove potential good reduction results at $p$ for its Jacobian. \\

Likewise, this article is dedicated to understanding the geometry of the modular curves $X_E^{\alpha}(N)$ over more general bases than characteristic zero fields, as is the case in most of the previously cited works. We also include a discussion of the degeneracy maps yielding the usual Hecke correspondences over a more general base. In particular, our construction provides the following:

\begin{itemize}[noitemsep,label=$-$]
\item proofs for the construction over $\Q$ that do not seem easily accessible otherwise, 
\item upgrades of the equations of $X_E^{\alpha}(n)$ derived over fields to equations valid over any regular excellent Noetherian $\Q$-algebra which is a domain (in Section \ref{subsect:upgrade-eqns}),
\item a study of the mixed characteristic case and the geometry of the special fibres. 
\end{itemize}

The last point does not seem to have been discussed in the literature, in spite of its importance to the formal immersion argument.\footnote{Or other explicit methods to determine the rational points such as the Chabauty--Kim method \cite{BDMTV21}.}

\bigskip

More precisely, our results are the following.\\

\mytheo[general-representability]{(See Propositions \ref{PG-representable-when-good} and Corollary \ref{general-smoothness-equivalence}) Let $G$ be a finite locally free commutative group scheme over a ring $R$ and $N \geq 3$ be an integer. The functor sending a $R$-scheme $S$ to the collection of $S$-isomorphism classes of couples $(F,i)$, where
\begin{itemize}[noitemsep,label=\tiny$\bullet$]
\item $F/S$ is an elliptic curve, 
\item $i: G_S \rar F[N]$ is an isomorphism of group schemes, 
\end{itemize}
is representable by an affine $R$-scheme $Y_G(N)$ of finite presentation. 

If the morphism $Y_G(N) \rightarrow \Sp{R}$ is surjective, it is smooth of relative dimension one.} 

\medskip
\noindent
\textbf{Remarks:}
\begin{itemize}[noitemsep,label=$-$]
\item The smoothness of $Y_G(N)$ does not require any assumption on $N$ or $R$. This contrasts with the situation of classical modular curves: the modular curve $Y(p)$ is smooth over $\Z[\zeta_p]$ at ordinary points, but not at supersingular points (cf. \cite[Theorems 13.7.6, 13.8.4]{KM}). Yet, a closer classical analogue of the $R$-scheme $Y_G(p)$ would be the $\Z$-scheme $Y(p)$, whose fibre at $p$ is nowhere reduced. 
\item Since $Y_G(N)$ is a moduli space of elliptic curves, there is a morphism $j: Y_G(N) \rightarrow \bA^1_R$ sending a couple $(F/S,i)$ as in Theorem \ref{general-representability} to the $j$-invariant of $F$ (as a $S$-point of $\bA^1_R$). 
\item As usual for classical modular curves, one can add extra structure to $Y_G(N)$. For instance, if $m \geq 1$ is an integer, there is a fine moduli space $Y_G(N,\Gamma_0(m))$ for triples $(F/S,i,C)$, where $(F/S,i)$ is as in Theorem \ref{general-representability} and $C$ is a cyclic subgroup of $F$ of degree $m$ (in the sense of \cite[(3.4)]{KM}), which is affine of finite presentation over $\Sp{R}$. 
\item When $m \geq 1$ is prime to $N$ and $d,t \geq 1$ are such that $dt \mid m$, we can construct a \emph{degeneracy map} $D_{d,t}: Y_G(N,\Gamma_0(m)) \rightarrow Y_G(N,\Gamma_0(t))$ corresponding to ``pushing the isomorphism $i$ through the isogeny whose kernel is given by the\footnote{cf. \cite[Theorem 6.7.2]{KM}} degree-$d$ subgroup of $C$'' (cf. Section \ref{gamma0m-structure}). The existence of such degeneracy maps is well-known to the experts in very broad generality, as is the fact that they give rise to Hecke operators. However, we are unaware of a literature reference spelling out that they also exist in the case of $Y_G(N,\Gamma_0(m))$.   
\item Over number rings, the Hecke operator $T_q$ can be defined after inverting $q$ then extended using the N\'eron model property, since we are dealing with Jacobians of curves with good reduction. However, this \emph{ad hoc} construction does not work \emph{a priori} with higher-dimensional families. Hence, we will construct our Hecke operator directly at the level of the $\operatorname{Pic}$ functors using the properties of the $X_G(N,\Gamma_0(q))$ over bases where $q$ is not invertible, thus bypassing the need for the N\'eron model property or Dedekind bases.  
\end{itemize}

\bigskip

\mytheo[etale-representability]{(See Propositions \ref{YG-Gamma0m} and \ref{XG-Gamma0m}) Let $N \geq 3$ be an integer and $R$ be a $\Z[1/N]$-algebra. Let $G$ be a finite locally free commutative group scheme over $R$, \'etale-locally isomorphic to $(\Z/N\Z)^{\oplus 2}$. There is a finite \'etale $(\Z/N\Z)^{\times}$-torsor $\mrm{Pol}(G)$ on $\Sp{R}$ such that for any $R$-scheme $S$, $\mrm{Pol}(G)(S)$ is the collection of bilinear alternating surjective maps $G_S \times G_S \rar (\mu_N)_S$. Let $m \geq 1$ be coprime to $N$. Then
\begin{enumerate}[noitemsep,label=$(\alph*)$]
\item the rule 
\vspace{-5pt}\[(F/S,i,C) \mapsto \left[G_S \times G_S \overset{i \times i}{\rar} F[N] \times F[N] \overset{\mrm{We}}{\rar} \mu_N\right]\vspace{-5pt}\] defines a flat morphism $\det: Y_G(N,\Gamma_0(m)) \rar \mrm{Pol}(G)$, with geometrically connected Cohen--Macaulay fibres that are pure of dimension one,
\item the map $(j,\det): Y_G(N,\Gamma_0(m)) \rar \bA^1_R \times_R \mrm{Pol}(G)$ is finite locally free.
\end{enumerate}

If, for every prime $p \mid m$, either $p^2 \nmid m$ or $p$ is invertible in $R$, then $Y_G(N,\Gamma_0(m)) \rar \mrm{Pol}(G)$ is smooth at any point which is not supersingular of residue characteristic dividing $m$. }

It should be noted that, since isogenies do not preserve the Weil pairing, the degeneracy maps $D_{d,t}$ described above are not morphisms of $\mrm{Pol}(G)$-schemes (see Proposition \ref{degen-AL-Weil}). 

\bigskip

With extra conditions on $N,m$ and the base ring, we can apply the general construction of \cite{KM} to produce compactifications of $Y_G(N,\Gamma_0(m))$. 

\mytheo[etale-compactification]{(See Propositions \ref{XG-Gamma0m} and \ref{XG-Gamma0m-connectedness}) Let $N \geq 3$ be an integer and $R$ be a $\Z[1/N]$-algebra which is \emph{regular} and \emph{excellent}. Let $G$ be a finite \'etale $R$-group scheme, \'etale-locally isomorphic to $(\Z/N\Z)^{\oplus 2}$, and $m \geq 1$ be coprime to $N$ such that, for every prime $p \mid m$, either $p^2 \nmid m$ or $p$ is invertible in $R$. 

Then the relative normalization of $j: Y_G(N,\Gamma_0(m)) \rightarrow \bA^1_R \subset \bP^1_R$ is a finite locally free $\bP^1_R$-scheme $X_G(N,\Gamma_0(m))$ satisfying the following properties: 

\begin{itemize}[noitemsep,label=$-$]
\item $Y_G(N,\Gamma_0(m)) \rar \mrm{Pol}(G)$ extends to a morphism $X_G(N,\Gamma_0(m)) \rightarrow \mrm{Pol}(G)$ which is smooth of relative dimension one at every point which is not supersingular of residue characteristic dividing $m$,
\item the reduced subscheme attached to the closed subspace $\{j=\infty\}$ of $X_G(N,\Gamma_0(m))$ (its \emph{cuspidal subscheme}) is finite \'etale over $R$ and is a relative effective Cartier divisor,  
\item The formation of $X_G(N,\Gamma_0(m))$ and its cuspidal subscheme commutes with any base change.
\end{itemize}
}

\bigskip

The application of interest is the case where the group $G$ is already the $N$-torsion subgroup scheme of some elliptic curve $E/R$, where $R$ is the ring of $S$-integers in some number field. In this case, the Weil pairing induces an isomorphism $(\Z/N\Z)^{\times} \simeq \mrm{Pol}(G)$. 

\medskip

\mycor[etale-representability-ell-curves]{Let $N \geq 3$ be an integer and $\alpha \in (\Z/N\Z)^{\times}$. 
\begin{enumerate}[noitemsep,label=(\alph*),topsep=3pt]
\item For any field $K$ of characteristic not dividing $N$ and any elliptic curve $E/K$, there is a smooth projective geometrically connected curve $X_E^{\alpha}(N)$ over $K$ along with a finite locally free morphism $j: X_E^{\alpha}(N) \rightarrow \bP^1_K$ such that $Y_E^{\alpha}(N) := j^{-1}(\bA^1_K)$ parametrizes couples $(F,i)$ where $F$ is an elliptic curve over some $K$-algebra $A$ and $i: E[N]_A \rightarrow F[N]$ is an isomorphism with determinant $\alpha$. 
\item If $L/K$ is a field extension and $E/K$ is an elliptic curve, $X_E^{\alpha}(N)_L \simeq X_{E_L}^{\alpha}(N)$. 
\item For any prime $q \nmid N$, there is a Hecke correspondence 
\vspace{-7pt}\[T_q: \mrm{Jac}(X_E^{\alpha}(N)) \rar \mrm{Jac}(X_E^{q\alpha}(N)),\vspace{-7pt}\] and the Hecke correspondences commute pairwise. 
\item Let $K$ be a global, or non-Archimedean local, field, and let $v$ be a place of $K$ not dividing $N$. If $E/K$ is an elliptic curve, then $X_E^{\alpha}(N)$ has potentially good reduction at $v$. If $E$ has good reduction at $v$ and $\overline{E}/k$ is the reduction of $E$ modulo $v$, then $X_E^{\alpha}(N)$ has good reduction at $v$ and its fibre at $v$ is $X_{\overline{E}}^{\alpha}(N)$. 
\end{enumerate}
}

\bigskip

We also study the curves $Y_G(N)$ when $N$ is not invertible over the base ring (see Sections \ref{subsect:ordinary} and \ref{subsect:supersingular}). Here is a typical result:    

\medskip

\mytheo[fiber-carac-p]{(See Corollaries \ref{ordinary-moduli} and \ref{fibre-struct-supersing+relrepmoduli}) Let $p$ be a prime and $N \geq 3$ be a multiple of $p$. Let $R$ be a Noetherian local ring of residue characteristic $p$, $\cE$ be an elliptic curve over $R$ and $G$ be the finite locally free $R$-group scheme $\cE[N]$. By Theorem \ref{general-representability}, $Y_G(N)$ is a smooth $R$-scheme of relative dimension one. 
\begin{itemize}[noitemsep,label=$-$]
\item Assume that the special fibre of $\cE$ is ordinary. Then the $j$-invariant map $Y_G(N) \rar \bA^1_{R}$ is quasi-finite flat, and it becomes locally free after any base change $R \rar k$ to a field, but the rank is not the same depending on whether $k$ has characteristic dividing $N$ or not. 
\item Assume that $\cE$ has supersingular special fibre. Then the geometric special fibre of $Y_G(N)$ is the disjoint union of a finite number of copies of $\bA^1$. Furthermore, the $j$-invariant map $Y_G(N) \rar \bA^1_{R}$ is locally constant on the special fibre, so hence it is not quasi-finite or flat. 
\end{itemize}}

\medskip

The main results of this text not figuring in the author's PhD dissertation \cite{Studnia-thesis} are the case of Theorem \ref{general-representability} where $N$ is not invertible in the base ring and Theorem \ref{fiber-carac-p}. 

In \cite{KO,Fisher-13,Frengley-12}, the curve $X_E^{\alpha}(N)$ is usually defined as a Galois twist of the classical modular curve $X(N)$ (or rather its geometrically connected variant defined over $\Q$). It is not immediately apparent that this construction is equivalent to the one above, but we describe precisely in Section \ref{moduli-is-twisting} the equivalence between the two constructions. The result below summarizes its contents. 

\medskip

\mytheo[property-of-twists]{Let $N \geq 3$ be an integer and $k$ be a field of characteristic not dividing $N$ with separable closure $k_s$. 
Let $X(N)$ denote the curve attached by Theorem \ref{etale-compactification} to the constant \'etale group $(\Z/N\Z)^{\oplus 2}_k$; it is endowed with a left action of $\GL{\Z/N\Z}$ given for any noncuspidal point $(E/S,\iota) \in X(N)(S)$ by $M \cdot (E/S,\iota) = (E/S,\iota \circ M^T)$. 

Let $G$ be a finite \'etale $k$-group scheme such that $G(k_s) \simeq (\Z/N\Z)^{\oplus 2}$. Fix a basis $(P,Q)$ of $G(k_s)$. The function $\rho: \mrm{Gal}(k_s/k) \rar \GL{\Z/N\Z}$ defined by 
\vspace{-5pt}\[\forall g \in \mrm{Gal}(k_s/k),\,\rho(g)\begin{pmatrix} P\\Q\end{pmatrix} = \begin{pmatrix}g(P)\\g(Q)\end{pmatrix}\vspace{-5pt}\] is an anti-homomorphism. Let $X(N)_{\rho^{-1}}$ denote the Galois twist (in the sense of Proposition \ref{cocycle-twist}) of $X(N)$ by the homomorphism $\rho^{-1}: \mrm{Gal}(k_s/k) \rar \GL{\Z/N\Z}$, which has a finite locally free morphism $j_{\rho}$ to $\bP^1_k$, and $\tau_{\rho}: X(N)_{k_s} \rar (X(N)_{\rho^{-1}})_{k_s}$ be the $\bP^1_k$-isomorphism given by the twist. 

There is a isomorphism of $\bP^1_k$-schemes $\iota: X_G(N) \rar X(N)_{\rho^{-1}}$ such that, for any non-cuspidal $(E/k_s,\varphi) \in Y_G(N)(k_s)$, one has 
\vspace{-5pt}\[\iota((E/k_s,\varphi)) = \tau_{\rho}(E/k_s,(a,b) \mapsto a\varphi(P)+b\varphi(Q)) \in X(N)_{\rho^{-1}}(k_s).\vspace{-5pt}\]

Suppose finally that $\mrm{Pol}(G)$ is constant (or, equivalently, $\det{\rho}$ is the mod $N$ cyclotomic character). 
Let $\mu_N^{\times}(k_s)$ denote the set of primitive $N$-th roots of unity in $k_s$ and $\underline{\mu_N^{\times}(k_s)}$ be the corresponding constant $k$-scheme. Then $\rho^{-1}$ acts (through its determinant) on $\mrm{Pol}((\Z/N\Z)^{\oplus 2}_k)=(\mu_N^{\times})_k$, and $\det: X(N)_k \rightarrow (\mu_N^{\times})_k$ can be twisted by $\rho^{-1}$ to a morphism $\mrm{We}_{\rho^{-1}}:  X(N)_{\rho^{-1}} \rar \underline{\mu_N^{\times}(k_s)}$. The morphism $\mrm{We}_{\rho^{-1}}$ is proper smooth, its fibres are geometrically connected curves over $k$, and the following diagram commutes:
\vspace{-6pt}\[
\begin{tikzcd}[ampersand replacement=\&]
X_G(N) \arrow{d}{\iota} \arrow{r}{\det} \& \mrm{Pol}(G)\arrow{d}{\mu\mapsto \mu(P,Q)}\\
X(N)_{\rho^{-1}} \arrow{r}{\mrm{We}_{\rho^{-1}}}\& \underline{\mu_N^{\times}(k_s)}
\end{tikzcd}
\vspace{-6pt}\]   
}

\medskip
\noindent
\textbf{Remark}\; A variant of this result holds for the curves $X_G(N,\Gamma_0(m))$, as long as $m$ is coprime to $N$ and is not divisible by the square of the characteristic of $k$. The twist also preserves the degeneracy maps. In particular, the isomorphism $\iota$ respects the Hecke operators. 

\medskip
\noindent

The article is structured as follows. In Section \ref{sec:representable}, we recall the formalism of moduli problems of \cite{KM} and expand on some generalization sketched in this book. We also discuss the properties of the natural degeneracy maps giving rise to Hecke operators. In Section \ref{situation-at-infty}, we generalize the technique of normal compactification introduced in \cite[Chapter 8]{KM}. Its main purpose is to show that, under suitable assumptions, any morphism of moduli schemes extends to a morphism between the compactified moduli schemes. Section \ref{application-moduli} is dedicated to applications to our moduli schemes: we prove Theorems \ref{general-representability} to \ref{fiber-carac-p}. Finally, in Section \ref{moduli-is-twisting}, we show that $X_E^{\alpha}(N)$ can be realized as a connected component of a Galois twist of $X(N)$ and show Theorem \ref{property-of-twists}.  
\bigskip

\nott{
\begin{itemize}[noitemsep,label=\tiny$\bullet$]
\item If $p$ is a prime and $n \geq 1$ is an integer, $v_p(n)$ denotes its $p$-adic valuation. 
\item If $R$ is a ring and $S$ is a set, $S_R$ (or sometimes $\underline{S}_R$) denotes the constant $R$-scheme with underlying set $S$. 
\item If $N \geq 1$ is an integer, $\mu_N$ denotes the finite flat $\Z$-group scheme $\Sp{\Z[X]/(X^N-1)}$, that is, the Cartier dual of the constant group scheme $(\Z/N\Z)_{\Z}$, and $\mu_N^{\times}$ is its ``scheme of $\Z/N\Z$-generators'' (cf. \cite[Corollary 1.12.3]{KM}). Being isomorphic to $\Sp{\Z[X]/\Phi_N(X)}$ (by \cite[Theorem 1.12.9]{KM}), where $\Phi_N(X)$ is the $N$-th cyclotomic polynomial, $\mu_N^{\times}$ is finite flat over $\Z$. Both $\mu_N$ and $\mu_N^{\times}$ are \'etale above $\Z[1/N]$.
\item If $N \geq 1$ is an integer, for every $d \mid N$, there is a natural closed immersion $\mu_d^{\times} \rar \mu_N$ (given on coordinate rings by $X \mapsto X$), and the induced map $\coprod_{d \mid N}{\mu_d^{\times}} \rar \mu_N$ is an isomorphism over $\Z[1/N]$.
\item If $N \geq 1$ is an integer, $\OO_N$ denotes the ring $\Z[1/N,\zeta_N]$, where $\zeta_N$ is a primitive $N$-th root of unity.
\item If $N \geq 1$ is an integer, every $a \in (\Z/N\Z)^{\times}$, $\underline{a}$ denotes the automorphism of $\mu_N$ (or $\mu_N^{\times}$, or $\Sp{\OO_N}$, or their coordinate rings depending on context, which is usually clear) given on $\OO(\mu_N)$ by $X \mapsto X^a$.
\item If $X$ is a scheme, $\OO_X$ (or $\OO$ if it is clear in context) denotes its structure sheaf; at a point $x \in X$, $\mfk{m}_{X,x}$ (or $\mfk{m}_x$ if $X$ is clear) denotes the maximal ideal of the local ring $\OO_{X,x}$, and $\kappa(x)$ denotes the residue field $\OO_{X,x}/\mfk{m}_{X,x}$. The ambient scheme $X$ is omitted from the notation because if $Z \rar X$ is an immersion, then the residue fields of $x \in Z$ and $x \in X$ are naturally identified. 
\end{itemize}

}

\bigskip

\textbf{Acknowledgements:} The inspiration to write this article largely arose in an email exchange with Barry Mazur, who pointed out to me how the definition of the twisted modular curves as moduli spaces made their reduction theory largely easier to understand. This text was largely a part of my PhD dissertation, and I am thus grateful to my supervisor Lo\"ic Merel for his support, help, comments and general encouragement to clarify the literature, as well as my thesis referees Pierre Parent and Henri Darmon for their suggestions and comments. The study of the situation in characteristic dividing $N$ was not in my thesis, and I started thinking about it after Jan Vonk and Vincent Pilloni both asked me this question. I am grateful to Luc Illusie, Adrian Vasiu, Frans Oort and Ofer Gabber for fruitful conversations about some of the algeometric geometry arguments contained in the paper, in particular those in characteristic $p$ and about $p$-divisible groups. I am grateful to the reviewer for their helpful suggestions and comments. 

The author was partly supported by ERC Grant ``GAGARIN'' $101076941$.

\bigskip

\section{Representable moduli problems}
\label{sec:representable}

\subsection{Basic definitions}

We recall the framework introduced in \cite[Chapter 4]{KM} and develop the extension sketched in (4.13) of \emph{loc.cit.}. 

Given a ring $R$, we denote by $\Ell_R$ the category defined as follows. Its objects are relative elliptic curves $E \rar S$ (which we often denote as $E/S$), where $S$ is a $R$-scheme; given two relative elliptic curves $F \rar T$ and $E \rar S$, with $S,T$ being $R$-schemes, a morphism $F/T \rar E/S$ is a Cartesan diagram of $R$-schemes

\[
\begin{tikzcd}[ampersand replacement=\&]
F\arrow{r}{\alpha}\arrow{d} \& E\arrow{d}\\
T \arrow{r} \&  S
\end{tikzcd}
\]

where $\alpha$ preserves the zero section, so that $F \rar E \times_S T$ is an isomorphism \emph{of elliptic curves over $T$} (by \cite[Theorem 2.5.1]{KM}).

\defi{A \emph{moduli problem over $R$} is a contravariant functor $\P: \Ell_R \rar \Set$.  }

\defi{A moduli problem $\P$ over $R$ is 
\begin{itemize}[label=$-$,itemsep=.5pt,topsep=2pt]
\item \emph{rigid} if for any object $E/S$ of $\Ell_R$, the action of $\mrm{Aut}(E/S)$ on $\P(E/S)$ (on the right) is free (i.e. only the identity admits fixed points). 
\item \emph{relatively representable} if for any object $E/S$ of $\Ell_R$, the functor
 \vspace{-7pt}\[T \in \Sch_S \longmapsto \P(E_T/T) \in \Set\vspace{-7pt}\]
  is representable by a $S$-scheme $\P_{E/S}$.
\item \emph{representable} if there exists an object $\cE/\M$ of $\Ell_R$ and $\alpha \in \P(\cE/\M)$ such that, for any object $E/S$ in $\Ell_R$, 
\vspace{-7pt}\[f \in \mrm{Hom}_{\Ell_R}(E/S,\cE/\M) \longmapsto \P(f)\alpha \in \P(E/S)\vspace{-7pt}\] 
is a bijection. We say that $(\cE/\M,\alpha)$ \emph{represents} the moduli problem $\P$, and that $\M$ is the \emph{moduli scheme} of $\P$ (which we will denote by $\MP$). We omit $\alpha$ or $\mathcal{E}$ from the notation when it is irrelevant.
\end{itemize}}

\rem{Let $\P$ be a moduli problem over $R$ and $\cE/\M$ be an object of $\Ell_R$. By Yoneda's lemma, $\P$ is represented by $(\cE/\M,\alpha)$ for some $\alpha \in \P(\cE/\M)$ if it is isomorphic to the functor $\mathrm{Hom}_{\Ell_R}(-,\cE/\M)$.}

\defi{Let $\P$ be a moduli problem over $R$, representable by the elliptic curve $\cE/\M$. The \emph{$j$-invariant map of $\P$} is the morphism $\M \rightarrow \bA^1_R$ given by $j(\cE/\M) \in \bA^1_R(\M)$.} 

\medskip

\defi{Let $\P$ be a moduli problem over $R$ and $R'$ be a $R$-algebra. The \emph{base change} $\P_{R'}$ of $\P$ to $R'$ is the moduli problem on $R'$ obtained by composing the forgetful functor $\Ell_{R'} \rar \Ell_R$ and $\P$.}

One has formally: 

\prop[base-change-representable]{Let $\P$ be a moduli problem over $R$ and let $R'$ be a $R$-algebra. Let $\P'=\P_{R'}$ be the base change of $\P$. 
\begin{itemize}[noitemsep,label=$-$]
\item If $\P$ is rigid, then so is $\P'$. 
\item If $\P$ is relatively representable, then so is $\P'$, and, for any object $E/S$ in $\Ell_{R'}$, $\P'_{E/S}$ and $\P_{E/S}$ are isomorphic $S$-schemes. 
\item Suppose that $\P$ is representable by $(\cE/\M,\alpha)$, and let $p': \cE_{R'}/\M_{R'} \rightarrow \cE/\M$ be the base change. Then $(\cE'/\M',\P(p')\alpha)$ represents $\P'$. In other words, we have a canonical isomorphism $\MPp{\P'} \rar \MP_{R'}$ (we say that the formation of the moduli scheme \emph{commutes with base change}).  
\end{itemize}}

\lem{The formation of the $j$-invariant map commutes with base change of representable moduli problems.}

\lem{Let $T: \P \Rightarrow \P'$ be a natural transformation of relatively representable moduli problem over $R$. Then, for every object $E/S$ in $\Ell_R$, $T$ induces a morphism of $S$-schemes $\P_{E/S} \rar \P'_{E/S}$ whose formation commutes with base change. 
If $\P,\P'$ are representable, then $T$ induces a map $\MP \rar \MPp{\P'}$ whose formation commutes with base change. }

\medskip

Let $\mathscr{C}$ be a property of morphisms of schemes. We say that it is \emph{good} if it satisfies the following conditions:
\begin{itemize}[noitemsep,label=$-$]
\item it is stable under pre- and post-composition by isomorphisms, 
\item it is stable under arbitary base change, 
\item it satisfies fpqc descent in the following sense: if $f: X \rar Y$ is a morphism such that there exists a fpqc cover by maps $Y_i \rar Y$ such that $f_{Y_i}: X \times_Y Y_i \rar Y_i$ satisfies the property $\mathscr{C}$, then $f$ satisfies the property $\mathscr{C}$. 
\end{itemize}

All the properties given in Proposition \ref{fpqc-descent-prop} are good.

\defi[good-property-moduli-problem]{Let $\P$ be a relatively representable moduli problem over $R$ and $\mathscr{C}$ be a good property of morphisms of schemes. We say that $\P$ \emph{satisfies $\mathscr{C}$ over $\Ell_R$} if for every object $E/S$ in $\Ell_R$, then $\P_{E/S} \rar S$ has property $\mathscr{C}$. 

If $Z \subset \bA^1_R$ is a locally closed subscheme, we say that $\P$ \emph{satisfies $\mathscr{C}$ above $Z$} if, for every object $E/S$ in $\Ell_R$ such that $j(E/S) \in \bA^1_R(S)$ factors through $Z$, then $\P_{E/S} \rar S$ satisfies property $\mathscr{C}$. 
}

One checks formally that:

\lem[base-change-good-property]{Let $\mathscr{C}$ be a good property of morphisms of schemes, $\P$ be a relatively representable moduli problem over $R$, and $\P'$ be the base change of $\P$ to an $R$-algebra $R'$. Then 
\begin{itemize}[noitemsep,label=$-$]
\item If $\P$ has property $\mathscr{C}$ over $\Ell_R$, then $\P'$ has property $\mathscr{C}$ over $\Ell_{R'}$,
\item If $R'$ is faithfully flat over $R$ and $\P'$ has property $\mathscr{C}$ over $\Ell_{R'}$, then $\P$ has property $\mathscr{C}$ over $\Ell_R$.
\end{itemize}}

\bigskip

\subsection{Representability results}

\prop[product-rel-rep]{(See \cite[(4.3.4)]{KM}) Let $\P, \P'$ be relatively representable moduli problems over $R$. Then the moduli problem over $R$
\vspace{-7pt}\[\P \times \P': E/S \longmapsto \P(E/S) \times \P'(E/S)\vspace{-7pt}\]
(with its natural action on the morphisms in $\Ell_R$) is relatively representable. Moreover, for every object $E/S$ of $\Ell_R$, one has a canonical isomorphism of $S$-schemes 
\vspace{-7pt}\[(\P \times \P')_{E/S} \rar \P_{E/S} \times_S \P'_{E/S}. \vspace{-7pt}\]
Suppose furthermore that $\P$ is representable by $(\cE/\MP,\alpha)$. Let $p': \cE'/\M' \rightarrow \cE/\M$ denote the pull-back by $\P'_{\cE/\M} \rightarrow \M$ and $\alpha' = \P(p')(\alpha) \in \P(\cE'/\M')$. Let $\beta \in \P'(\cE'/\M')$ correspond to the identity under the identification
\vspace{-7pt}\[ \P'(\cE'/\M') = \mrm{Hom}_{\M}(\M',\P'_{\cE/\M}).\vspace{-7pt}\]
Then $\P \times \P'$ is represented by $(\cE'/\M',(\alpha',\beta))$. 
}

\demo{This is asserted without proof in \cite[(4.3.4)]{KM}, so we give one here. 
The first claim is clear by the universal property of the fibre product of schemes. Let us discuss the second claim: for any object $F/S$ of $\Ell_R$, one has the following sequence of fonctorial bijections
\begin{align*}
\P(F/S) \times \P'(F/S) &\simeq \left\{(\varphi,y) \mid \small\begin{array}{c} \varphi: F/S \rar \cE/\M \\ y \in \P'(F/S) \end{array}\right\} \simeq \left\{(\varphi,y) \mid \small\begin{array}{c} \varphi: F/S \rar \cE/\M \\ y \in \P'(\varphi^{\ast}(\cE/\M)) \end{array}\right\} \\
&\simeq \left\{(\varphi,\psi) \mid \small\begin{array}{c} \varphi: F/S \rar \cE/\M \\ \psi \in \mrm{Hom}_{\M,\varphi}(S,\P'_{\cE/\M}) \end{array}\right\} \\
&\simeq \mrm{Hom}_{\Ell_R}(F/S,(p')^{\ast}\cE/\M) = \mrm{Hom}_{\Ell_R}(F/S,\cE'/\M'), 
\end{align*}
where $\mrm{Hom}_{\M,\varphi}(S,\P'_{\cE/\M})$ denotes the set of morphisms $\psi: S \rightarrow \P'_{\cE/\M}$ that are morphisms of $\M$-schemes, where $\varphi: F/S \rar \cE/\M$ induces a structure of $\M$-scheme on $S$. Hence $\cE'/\M'$ represents $\P \times \P'$. It is easy to check that $(\alpha',\beta) \in \P(\cE'/\M') \times \P'(\cE'/\M')$ is sent by these bijections to $\mathrm{id} \in \mrm{Hom}_{\Ell_R}(\cE'/\M',\cE'/\M')$, which concludes by Yoneda's lemma. }

\smallskip

\prop[relrep-rigid-implies-rep]{Let $R$ be a ring and $\P$ be a moduli problem over $R$. Assume that $\P$ is rigid, relatively representable and affine. Then $\P$ is representable. }

\demo{When $R=\Z$, this is exactly \cite[(4.7.0)]{KM}, but the argument in \emph{loc.cit.} works for any ring. Indeed, the moduli problems used to cover $\Ell_{\Z}$, the ``Legendre'' and ``na\"ive level three'' problems (over $\Z[1/2]$ and $\Z[1/3]$ respectively) are still representable by smooth affine curves over $R[1/2]$ and $R[1/3]$ for any ring $R$, by Proposition \ref{base-change-representable}, and the representability argument in \cite[pp. 112--116]{KM} does not use any property of $\Z$ as base ring.}

\prop[representability-descends]{Let $R \rar R'$ be a faithfully flat \'etale ring homomorphism, and let $\P$ be a moduli problem over $R$. Assume that:
\begin{itemize}[noitemsep,label=$-$]
\item $\P_{R'}$ is representable and affine over $\Ell_{R'}$,
\item For all objects $E/S$ of $\Ell_R$, $T \in \Sch_S \longmapsto \P(E_T/T)$ is a fpqc sheaf. 
\end{itemize}
Then $\P$ is representable.
}

\demo{Let $E/T$ be an object of $\Ell_R$ and $E'/T'$ be its base change by $T \times_R \Sp{R'} \rightarrow T$.  

The pre-sheaf $\mathcal{F}: U \mapsto \P(E_U/U)$ on the category of $T$-schemes is a fpqc sheaf, and its restriction to the category of $T'$-schemes is represented by the affine morphism $(\P_{R'})_{E'/T'} \rightarrow T'$. By Proposition \ref{desc-affine-schemes}, $\mathcal{F}$ is representable by an affine morphism $X \rightarrow T$. Therefore $\P$ is relatively representable and affine. 

Let $u \in \Aut{E/T} \backslash \{\mrm{id}\}$. By faithfully flat descent, $u \in \Aut{E'/T'} \backslash \{\mrm{id}\}$, so it has no fixed points in $\P(E'/T')=\P_{R'}(E'/T')$ because $\P_{R'}$ is rigid. Hence $u$ has no fixed points in $\P(E/T)$. Thus $\P$ is rigid, hence representable by Proposition \ref{relrep-rigid-implies-rep}.}

\bigskip

\subsection{Basic moduli problems and elementary applications}

We now recall the definitions of the moduli problems $[\Gamma(m)]$ and $[\Gamma_0(m)]$, following mostly \cite[Chapter 3]{KM}. 

\defi[gammam-definition]{Let $m \geq 1$ be an integer. The moduli problem $[\Gamma(m)]$ is defined as follows: 
\begin{itemize}[noitemsep,label=$-$]
\item For any elliptic curve $E/S$, $[\Gamma(m)](E/S)$ is the set of couples $(P,Q) \in E[m](S)^2$ such that as Cartier divisors on $E$, one has 
\vspace{-7pt}\[E[m] = \sum_{(a,b) \in (\Z/m\Z)^2}{[aP+bQ]}.\vspace{-7pt}\]  
\item Given a morphism $F/T \rar E/S$ in $\Ell_{\Z}$, the morphism $[\Gamma(m)](E/S) \rar [\Gamma(m)](F/T)$ is induced by the map $E(S) \rar E_T(T) \simeq F(T)$.  
\end{itemize}
}

\rem{When $m \in \OO(S)^{\times}$, the definition describes the classical notion of $\Gamma(m)$-structure by \cite[Lemma 1.8.3, Theorem 1.10.1]{KM}. }

\rem{The group $\GL{\Z/m\Z}$ acts on $[\Gamma(m)]$ on the left by the rule 
\vspace{-7pt}\[\begin{pmatrix} a & b\\c & d\end{pmatrix} \cdot (P,Q) = (aP+bQ, cP+dQ).\vspace{-7pt}\]
Note that our definition of $[\Gamma(m)]$ is not the exact same as \cite[(3.1)]{KM}, but that there is an isomorphism $[\Gamma(m)] \simeq [\Gamma(m)]_{\mathrm{KM}}$ under which the left action of $g \in \GL{\Z/m\Z}$ on $[\Gamma(m)]$ corresponds to the right action of its transpose on $[\Gamma(m)]_{\mathrm{KM}}$ (cf. (10.3) of \emph{op.cit.}).
}

To define $[\Gamma_0(m)]_{\Z}$, we start by defining cyclic subgroups of elliptic curves. 

\defi[cyclic-subgroup]{Let $E$ be an elliptic curve over a scheme $S$ and $m \geq 1$. A \emph{cyclic subgroup of order $m$ of $E$} is a closed subgroup scheme $C$ of $E$, finite locally free over $S$ of rank $m$, and such that, fppf-locally over $S$, there exists a section $P$ of $C$ such that $C$ agrees with $\sum_{i=1}^m{[iP]}$ as Cartier divisors on $E$. The section $P$ is a (local) \emph{generator} of $C$. } 

\defi[gamma0m-definition]{Let $m \geq 1$. The moduli problem $[\Gamma_0(m)]$ is defined as follows:
\begin{itemize}[noitemsep,label=$-$]
\item For any elliptic curve $E/S$, $[\Gamma_0(m)](E/S)$ is the set of cyclic subgroups of order $m$ of $E$. 
\item Given a morphism $F/T \rar E/S$ in $\Ell_{\Z}$ inducing an isomorphism $\iota: F \rightarrow E_T$ of elliptic curves, the morphism $[\Gamma_0(m)](E/S)\rar [\Gamma_0(m)](F/T)$ is given by $C \mapsto \iota^{-1}(C_T)$. 
\end{itemize}
}

\rem{\cite[(3.4)]{KM} also gives another definition for the moduli problem $[\Gamma_0(m)]_R$, as the collection of isomorphism classes of isogenies $E \rar E'$ whose kernel is a cyclic subgroup of order $m$ of $E$. The two definitions are equivalent: if $E/S$ is an elliptic curve and $G$ is a subgroup scheme of $E$ which is finite locally free over $S$, then one can show using \cite[Lemma 07S7]{Stacks}\footnote{I am grateful to MathOverflow user SashaP for pointing out this result to me.} that $E' := E/G$ is an elliptic curve over $S$ and $G$ is the kernel of the obvious isogeny $\pi: E \rar E'$. }

\prop[basic-moduli]{Let $R$ be a ring and $m \geq 1$ be an integer. The moduli problems $[\Gamma_0(m)]_R$ and $[\Gamma(m)]_R$ are relatively representable, finite locally free of constant rank on $\Ell_R$, \'etale above $\Z[1/m]$ (over $\Ell_R$). Moreover, when $m \geq 3$ is invertible in $R$, $[\Gamma(m)]_R$ is representable by a smooth affine $R$-scheme of relative dimension one. }

\demo{The moduli problems are relatively representable with the given properties by \cite[Lemma 4.12.1, Theorem 5.1.1]{KM} and base change. The second part follows from \cite[Corollary 4.7.1]{KM}, where the rigidity assumptions are checked in \cite[(2.7)]{KM}.  
}

This representability result helps us relate properties of representable moduli problems (in the sense of Definition \ref{good-property-moduli-problem}) to properties of their moduli schemes or their $j$-invariant map. 

\lem[decomposition-of-moduli-scheme]{Let $\P$ be a representable moduli problem over $R$ and $\ell$ be an odd prime. Let $\cE/Y(\ell)$ be the elliptic curve representing $[\Gamma(\ell)]_{R[1/\ell]}$. Then $Y(\ell)$ is affine and smooth of relative dimension one over $R[1/\ell]$, and one has a commutative diagram of $R[1/\ell]$-schemes 
\vspace{-6pt}\[
\begin{tikzcd}[ampersand replacement=\&]
\P_{\cE/Y(\ell)} \arrow{r}\arrow{d}{\pi} \& Y(\ell)\arrow{d}{j(\cE)}\\
\MP_{R[1/\ell]} \arrow{r}{j_{\P}} \&  \bA^1_{R[1/\ell]} 
\end{tikzcd}
\vspace{-10pt}\]
where $\pi: \P_{\cE/Y(\ell)} \simeq \MPp{\P_{R[1/\ell]},[\Gamma(\ell)]_{R[1/\ell]}} \rightarrow \MP_{R[1/\ell]}$ is finite \'etale of rank $|\GL{\F_{\ell}}|$ and $j(\cE): Y(\ell) \rar \bA^1_{R[1/\ell]}$ is finite locally free.
}

\demo{The map $\pi$ is well-defined by Proposition \ref{product-rel-rep} and finite \'etale surjective by the properties of the moduli problem $[\Gamma(\ell)]_{\Z[1/\ell]}$ over $\Ell_{\Z[1/\ell]}$. The diagram clearly commutes, and $Y(\ell)$ is smooth over $R[1/\ell]$ of relative dimension one by Proposition \ref{basic-moduli}. To check that $j(\cE)$ is finite locally free, we may assume that $R=\Z[1/\ell]$, so it suffices to check that $j(\cE)$ is finite flat. This morphism is finite by \cite[Proposition 8.2.2]{KM} and between two smooth $\Z[1/\ell]$-schemes of relative dimension one, hence flat (see e.g. \cite[(4.12), Theorem 5.1.1, Notes on Chap. 4]{KM}). }

\cor[over-Ell-concrete]{Let $\P$ be a representable moduli problem over $R$. 
\begin{itemize}[noitemsep,label=$-$]
\item If $\P$ is affine over $\Ell_R$, then $\MP$ is affine.
\item If $\P$ is surjective (resp. locally of finite presentation, locally quasi-finite, finite with $R$ Noetherian, flat) over $\Ell_R$, its $j$-invariant map $\MP \rar \bA^1_R$ is surjective (resp. locally of finite presentation, locally quasi-finite, finite, flat). 
\item If $\P$ is locally of finite presentation over $\Ell_R$ and quasi-finite flat (resp. \'etale) above some open subscheme $U \subset \bA^1_R$, then $f_{U}: j^{-1}(U) \subset \MP \rar \Sp{R}$ is flat (resp. smooth) of relative dimension one with Cohen--Macaulay fibres.
\end{itemize} 
}

\demo{Working Zariski-locally over $R$, we may choose an odd prime $\ell$ invertible in $R$ and consider the diagram of Lemma \ref{decomposition-of-moduli-scheme}. If $\P$ is affine over $\Ell_R$, then $\MPp{\P,[\Gamma(\ell)]}$ is affine over $Y(\ell)_{R[1/\ell]}$, hence affine, so $\MP$ is affine by \cite[Lemma 01ZT]{Stacks}. The argument for the second item is analogous using instead \cite[Lemmas 036K, 036N, 036O]{Stacks}. 

For the third item, base change the diagram of Lemma \ref{decomposition-of-moduli-scheme} by $U \rar \bA^1_R$, and let $V,W$ be the pre-images of $U$ in $\MP$ and $\MPp{\P,[\Gamma(\ell)]_R}$ respectively. As above, $j: V \rar U$ is flat locally of finite presentation with discrete fibres, so $f_U$ is flat, locally of finite presentation, with pure one-dimensional fibres by \cite[Lemma 02NL]{Stacks}.  By \cite[$\text{IV}_2$ (6.3.5), (6.5.2)]{EGA}, it is enough to show that $W$ is Cohen--Macaulay (resp. regular) when $R$ is a field, which follows from \cite[$\text{IV}_2$ (6.8.3)]{EGA} since $W \rar Y(\ell)$ is locally quasi-finite and flat (resp. \'etale) and $Y(\ell)$ is smooth. }

\bigskip

\subsection{Level structures and degeneracy maps}
\label{gamma0m-structure}

We make the following definition inspired by \cite[(4.14)]{KM}. 

\defi[category-strongly-free-gpsch]{Given a ring $R$ and an integer $N \geq 1$, let $\mathbf{FGp}^N_R$ be the following category:
\begin{itemize}[noitemsep,label=$-$]
\item Its objects are commutative $S$-group schemes $G$ killed by $N$, where $S$ is a $R$-scheme, such that, for every $m \mid N$, the subgroup scheme $G[m]$ is finite locally free over $S$, 
\item A morphism $G'/S' \rar G/S$ is a Cartesian diagram    
\vspace{-7pt}\[\begin{tikzcd}[ampersand replacement=\&]G' \arrow{r} \arrow{d} \& G\arrow{d}\\S' \arrow{r}\& S\end{tikzcd}\vspace{-7pt}\]
such that the induced map $f: G' \rar G \times_S S'$ is an isomorphism of $S'$-group schemes.
\end{itemize}}

\defi{Let $N \geq 1$ be an integer and $R$ be a ring. If $\P$ is a moduli problem over $R$, a \emph{level $N$ structure} for $\P$ is a contravariant functor 
\vspace{-7pt}\[F: \mathbf{FGp}^N_R \rar \Set\vspace{-7pt}\]
 such that $\P$ is the composition of $F$ and the obvious ``$N$-torsion functor'' 
\vspace{-7pt}\[\Ell_R \rightarrow \mathbf{FGp}_R^N,\qquad E/S \mapsto E[N]/S.\vspace{-7pt}\] 
Unless it is both relevant and not obvious, when we claim that a moduli problem $\P$ over $R$ is of level $N$, we omit the functor $F$ from the notation.}

\defi{Let $\P,\P'$ be moduli problems over $R$, endowed with level $N$ structures by the functors $F,F'$. A \emph{morphism of moduli problems of level $N$} from $(\P,F)$ to $(\P',F')$ is a natural transformation $f: F \Rightarrow F'$. In particular, $f$ defines a morphism of moduli problems $\P \Rightarrow \P'$.}

\rems{\begin{itemize}[noitemsep,label=$-$]
\item If $G/S$ is an object of $\mathbf{FGp}_R^N$, so is the base change $G_{S'}/S'$ for any $S$-scheme $S'$.
\item As in the previous section, there is a notion of base change for moduli problems of level $N$. 
\item Suppose that $N \mid N'$. The rule $G/S \mapsto G[N]/S$ defines a functor $\mathbf{FGp}_R^{N'} \rightarrow \mathbf{FGp}_R^N$. This lets us endow any moduli problem of level $N$ (say, $(\P,F)$) over $R$ with a natural level $N'$ structure $F_{N \rar N'}$. If $N' \mid N''$, one has $(F_{N \rightarrow N'})_{N' \rightarrow N''} = F_{N \rightarrow N''}$, and this operation clearly commutes with base change. 
\item If $S_0$ is a $R$-scheme, the moduli problem $[S_0]: E/T \mapsto S_0(T)$ over $R$ is relatively representable (with $[S_0]_{E/S} = S_0 \times_R S$), and it has a natural level $N$ structure for any $N$. 
\item Suppose that $(\P,F)$ is a moduli problem of level $N$ over $R$, $E,E'$ are elliptic curves over a $R$-scheme $S$, and $f: E \rar E'$ is an isogeny of degree prime to $N$. Then $f$ induces an isomorphism of $S$-group schemes $E[N] \rar E'[N]$, hence a bijection $\P(E/S) \rar \P(E'/S)$, which we will denote somewhat abusively by $\P f$ or even $f$. When $E=E'$ and $f$ is the multiplication by some integer $d \equiv 1 \mod{N}$, $f$ acts trivially on $E[N]$: this means that $(\Z/N\Z)^{\times}$ acts\footnote{Note that this action and the above isomorphism depend on the choice of $F$. We usually remove this from the notation, because we never consider two different level structures on the same moduli problem.} on $\P$; we will denote this action by $[\cdot]$. In particular, one has a natural action of $\widehat{\Z}^{\times}$ on $(\P,F)$ by reducing modulo $N$ and applying $[\cdot]$. If $N \mid N'$, the constructions are unchanged if we replace $(\P,F)$ by $(\P,F_{N \rightarrow N'})$. 
\item If $(\P,F)$, $(\P',F')$ are moduli problems of respective levels $N,N'$, the moduli problem $\P \times \P'$ has the natural level $M := \mathrm{lcm}(N,N')$ structure $G/S \mapsto F(G[N]) \times F'(G[N'])$. 
\item By \cite[Proposition 1.7.2]{KM}, if $a,b \geq 1$ are coprime integers and $G$ is a finite locally free commutative group scheme over a scheme $S$ killed by $ab$, then $G[a], G[b]$ are finite locally free $S$-group schemes and the addition $G[a] \times_S G[b] \rightarrow G$ is an isomorphism.
\item However, if $G$ is a finite locally free commutative $R$-group scheme killed by $N'$ and $N \mid N'$, it is not true in general that $G[N]$ is $R$-flat\footnote{I am grateful to F. Oort, O. Gabber and P. Bruin for the following constructions. If $R$ is a DVR of mixed characteristic $(0,p)$, $G(\overline{\mrm{Frac}(R)}) \simeq \Z/p^2\Z$ and the special fibre of $G$ is killed by $p$, then $G[p]$ is not flat. Let $E/R$ be an elliptic curve with good supersingular reduction, and let then $G$ be the closure of the subgroup scheme generated by a point of $E$ of order $p^2$. Alternately, the group scheme $\underline{\Z/p^2\Z}$ comes from some $\alpha \in \mrm{Ext}^1(\underline{\Z/p\Z},\underline{\Z/p\Z})$. Let $R=\Z_p[\zeta_p]$ and $f: (\Z/p\Z)_R \rightarrow (\mu_p)_R$ be a generic isomorphism. Let $G$ be the group scheme corresponding to $f_{\ast}(\alpha)$.} 
This is the reason behind the stronger flatness assumption in the definition of $\mathbf{FGp}_R^N$. 
\end{itemize}
}

\prop[basic-moduli-level]{On any ring $R$, for any integer $m \geq 1$, the relatively representable moduli problems $[\Gamma_0(m)]_R$, $[\Gamma(m)]_R$ are endowed with natural level $m$ structures. Moreover, the induced action of $\widehat{\Z}^{\times}$ on $[\Gamma_0(m)]_R$ is trivial.}

\demo{The definition of the level $m$ structure is given by \cite[Proposition 1.10.6]{KM}. If $E/S$ is an object of $\Ell_R$ and $C \in [\Gamma_0(m)](E/S)$, one has $dC=C$ whenever $d \in \Z$ is prime to $m$, so $\widehat{\Z}^{\times}$ acts trivially on $[\Gamma_0(m)]_{R}$.  
}

\rem{Let $m \geq 1$ and $a \in (\Z/m\Z)^{\times}$. Then the action of the matrix $aI_2 \in \GL{\Z/m\Z}$ on the moduli problem $[\Gamma(m)]$ is the same as the action of $[a]$ given by the natural level $m$ structure on $[\Gamma(m)]$.}

\cor[rep-plus-gamma0]{Let $R$ be a ring and $m,N \geq 1$ be integers. Let $\P$ be a representable moduli problem of level $N$ over $R$. The moduli problem $\P \times [\Gamma_0(m)]_R$ is representable and has a natural level $\mathrm{lcm}(N,m)$ structure.} 

\demo{$[\Gamma_0(m)]_R$ is relatively representable of level $m$ by Propositions \ref{basic-moduli} and \ref{basic-moduli-level}, so $\P \times [\Gamma_0(m)]_R$ is representable by Proposition \ref{product-rel-rep} and has a level structure as a product of moduli problems.}

\medskip

In view of our next focus on degeneracy maps, let us recall a few facts on the standard subgroups of cyclic subgroups of elliptic curves, that we will use freely. 

\prop[cyclic-standard]{(\cite[Theorems 6.7.2, 6.7.4]{KM}) Let $C$ be a cyclic subgroup of order $n \geq 1$ of an elliptic curve $E/S$. 
\begin{itemize}[noitemsep,label=\tiny$\bullet$]
\item if $d \mid n$, there is a ``standard'' subgroup scheme $C\{d\}$ of $C$ which is cyclic of degree $d$, uniquely defined by the following property: if $P \in C(T)$ is a local generator of $C$, $\frac{n}{d}P$ generates $C\{d\}$. 
\item if $d \mid d' \mid n$, then $C\{d\} = (C\{d'\})\{d\}$. Moreover, if $f: E \rightarrow E/C\{d\}$ is the obvious isogeny, $f(C)=C/C\{d\}$ is cyclic of degree $\frac{n}{d}$, and $C\{d'\}/C\{d\} = f(C)\left\{\frac{d'}{d}\right\}$. 
\item if $d \mid n$ and $g: E \rightarrow E'$ is an isogeny of degree prime to $n$, then $g(C\{d\}) = g(C)\{d\}$. 
\item if $n = ab$ and $a,b$ are coprime, $C\{a\}$ and $C\{b\}$ are the $a$- and $b$-torsion subschemes $C[a]$ and $C[b]$ of $C$, and the addition $C\{a\} \times C\{b\} \rightarrow C$ is an isomorphism. 
\end{itemize}
} 

\demo{The first and second points are \cite[Theorems 6.7.2, 6.7.4]{KM} respectively. For the third point, $g$ induces an isomorphism $E[n] \rightarrow E'[n]$, so $g(C)$ is cyclic of order $n$, and $g$ maps a generator of $C$ to a generator of $g(C)$ (cf. \cite[(1.10)]{KM}); this implies the conclusion by the definition of the standard cyclic subgroup. The fourth point follows from properties of prime factorization \cite[Corollary 1.7.3, Lemma 3.5.1]{KM}.   }

\defi[degeneracies]{Let $m,d,t,N > 0$ be integers such that $dt \mid m$ and $m$ is coprime to $N$. For any representable moduli problem $\P$ over $R$ of level $N$, we define a \emph{degeneracy map} $D_{d,t}(\P): \MPp{\P,[\Gamma_0(m)]_R} \rar \MPp{\P,[\Gamma_0(t)]_R}$ (denoted $D_{d,t}$ if $\P$ is clear) as follows. Given an object $E/S$ of $\Ell_R$, $C \in [\Gamma_0(m)](E/S)$ and $\alpha \in \P(E/S)$, let $\pi: E \rightarrow E/C\{d\}$ denote the obvious isogeny, and we set
\vspace{-7pt}\[D_{d,t}(\P) \left[(E/S,\alpha,C)\right] = \left[(E/C\{d\},\pi(\alpha),\pi(C)\{t\})\right].\vspace{-7pt}\]
}

\rem{The map $D_{d,t}(\P)$ can be seen as an ``exotic'' morphism of moduli problems in the sense of \cite[(11.2)]{KM}, but we will not take this point of view. }

The following result is formal. 

\lem[elem-degeneracies]{We keep the setting of Definition \ref{degeneracies}.
\begin{itemize}[noitemsep,label=\tiny$\bullet$]
\item The formation of any $D_{d,t}$ is natural with respect to $\P$ (for morphisms of representable moduli problems of level $N$) and commutes with base change and ``level-raising'' (i.e. viewing $\P$ as a moduli problem of level $N'$, where $N'$ is a multiple of $N$ coprime to $m$). 
\item $D_{1,m}(\P)$ is the identity.
\item Let $r,s > 0$ be such that $rs \mid t$. Then $D_{r,s}(\P) \circ D_{d,t}(\P) = D_{dr,s}(\P)$.  
\end{itemize}
}

\lem[degn-compat]{Let $\P$ be a representable moduli problem over $R$ of level $N$ and $m,m' \geq 1$ be such that $N,m,m'$ are pairwise coprime. Let $\P'=\P \times [\Gamma_0(m)]_R$. There is a natural (with respect to representable moduli problems of level $N$) isomorphism of moduli problems (of level $Nmm'$) $\P' \times [\Gamma_0(m')]_R \simeq \P \times [\Gamma_0(mm')]_R$. 

Let $r,s \geq 1$ divide $m'$ such that $rs \mid m'$. Then, under the above isomorphism, the following diagram commutes:
\vspace{-7pt}\[
\begin{tikzcd}[ampersand replacement=\&]
\MPp{\P',[\Gamma_0(m')]_R} \arrow{r}{\sim} \arrow{d}{D_{r,s}(\P')} \& \MPp{\P,[\Gamma_0(mm')]_R} \arrow{d}{D_{r,ms}(\P)}\\
\MPp{\P',[\Gamma_0(s)]} \arrow{r}{\sim} \& \MPp{\P,[\Gamma_0(ms)]_R}
\end{tikzcd}
\vspace{-7pt}\]
}

\demo{By \cite[Lemma 3.5.1]{KM}, one has an isomorphism $[\Gamma_0(mm')] \simeq [\Gamma_0(m)] \times [\Gamma_0(m')]$, from which the conclusion follows using properties of standard cyclic subgroups.}

\medskip

\prop[degn-cartesian]{Let $F: \P \Rightarrow \P'$ be a morphism of representable moduli problems of level $N$ over $R$. Let $m,d,t\geq 1$ be such that $dt \mid m$ and $m$ is coprime to $N$. Then the following diagram is Cartesian:
\vspace{-7pt}\[
\begin{tikzcd}[ampersand replacement=\&]
\MPp{\P,[\Gamma_0(m)]_R} \arrow{r}{D_{d,t}(\P)}\arrow{d}{F} \& \MPp{\P,[\Gamma_0(t)]_R}\arrow{d}{F}\\
\MPp{\P',[\Gamma_0(m)]_R} \arrow{r}{D_{d,t}(\P')} \&  \MPp{\P',[\Gamma_0(t)]_R}
\end{tikzcd}
\vspace{-10pt}\]
}

\demo{The diagram clearly commutes. Let $S$ be a $R$-scheme, and let
\vspace{-7pt}\[\tau_1=(E'/S,\alpha',C') \in \MPp{\P,[\Gamma_0(t)]_R}(S), \tau_2=(E/S,\beta,C) \in \MPp{\P',[\Gamma_0(m)]_R}(S)\vspace{-7pt}\]
 be such that $F(\tau_1)=D_{d,t}(\tau_2) \in \MPp{\P',[\Gamma_0(t)]_R}(S)$. 
Then there exists an isogeny $\pi: E \rar E'$ of degree $d$ such that $\ker{\pi} = C\{d\}$, $C' = \pi(C)\{t\}$, and $F(\alpha')=\pi(\beta)$. 

Let $\alpha \in \P(E/S)$ be the pre-image of $\alpha'$ under $\pi$ and $\tau = (E/S,\alpha,C) \in \MPp{\P,[\Gamma_0(m)]_R}(S)$. Then $\tau$ lies above $\tau_1$ by construction. Moreover, one has $\pi(F(\alpha))=F(\pi(\alpha))=F(\alpha')=\pi(\beta)$, hence $F(\alpha)=\beta$ and $\tau$ lies above $\tau_2$. 

Suppose that $\tau_0 = (E_0/S,\gamma,D) \in \MPp{\P,[\Gamma_0(m)]_R}(S)$ is a point above $(\tau_1,\tau_2)$. Since $F(\tau_0) = \tau_2$, there is an isomorphism $f: E_0 \rar E$ mapping $F(\gamma)$ to $\beta$ and $D$ to $C$. Moreover, since $D_{d,t}(\P)(\tau_0)=\tau_1$, there is an isogeny $g: E_0 \rar E'$ of degree $d$ such that $\alpha'=g(\gamma)$, $\ker{g} = D\{d\}$ and $C' = g(D)\{t\}$. Then $g, \pi \circ f: E_0 \rar E'$ are two isogenies with the same kernel, so $\pi \circ f = \psi \circ g$ for some automorphism $\psi$ of $E'$. Hence
\vspace{-7pt} \[\psi(F(\alpha'))=\psi(F(g(\gamma)))=F(\psi(g(\gamma)))=F(\pi(f(\gamma)))=\pi(f(F(\gamma)))=\pi(\beta)=F(\alpha').\vspace{-7pt}\]

Since $\P'$ is representable, it is rigid, so $\psi$ is the identity. Hence $\pi(\alpha)=\alpha'=\pi(f(\gamma))$, so $f(\gamma)=\alpha$. Therefore $\tau_0=\tau$ and we are done.
}

\lem[rel-coprime-elliptic-curves]{Let $E,F$ be relative elliptic curves over a basis $S$, endowed with cyclic subgroups $C,D$ of coprime degrees $c,d \geq 1$. Suppose we are given an isomorphism $\iota: E/C \rar F/D$ of elliptic curves (thus endowing $F$ with the structure of finite locally free $E/C$-scheme of degree $d$). Then $E \times_{E/C} F$ is an elliptic curve over $S$.}

\demo{$G:=E \times_{E/C} F$ is clearly a proper commutative $S$-group scheme.   
Moreover, the projection $G \rar F$ is the base change of the isogeny $E \rar E/C$ of degree $c$, so it is \'etale above $\Z[1/c]$ and $G \rar S$ is smooth of relative dimension one above $\Z[1/c]$. Similarly, $G \rar S$ is smooth of relative dimension one above $\Z[1/d]$. Since $c,d$ are coprime, $G \rar S$ is smooth of relative dimension one.

To show that the geometric fibres of $G \rar S$ are connected, we may assume that $S = \Sp{k}$ for some field $k$. Let $p: E \rar E/C,\; q: F \rar F/D$ be the projections. In $G$, the multiplication by $c$ is the composition $(p^{\vee}\iota^{-1}q,c \cdot \mrm{id}) \circ \mrm{pr}_F$, where $\mrm{pr}_F: G \rightarrow F$ is the second projection: hence $c$ kills $\pi_0(G_{\overline{k}})$. Similarly, $d$ kills $\pi_0(G_{\overline{k}})$, so $\pi_0(G_{\overline{k}})$ is trivial and $G_{\overline{k}}$ is connected. }

\prop[two-degn-cartesian]{Let $\P$ be a representable moduli problem of level $N$ over $R$. Let $m \geq 1$ be an integer coprime to $N$. Let $s',t' \geq 1$ be two coprime divisors of $m$, and $s,t \geq 1$ be divisors of $s',t'$ respectively. Then the following diagram is Cartesian:
\vspace{-7pt}\[
\begin{tikzcd}[ampersand replacement=\&]
\MPp{\P,[\Gamma_0(m)]_R} \arrow{r}{D_{s,\frac{m}{s'}}}\arrow{d}{D_{t,\frac{m}{t'}}} \& \MPp{\P,[\Gamma_0\left(\frac{m}{s'}\right)]_R}\arrow{d}{D_{t,\frac{m}{s't'}}}\\
\MPp{\P,[\Gamma_0\left(\frac{m}{t'}\right)]_R} \arrow{r}{D_{s,\frac{m}{s't'}}} \&  \MPp{\P,[\Gamma_0\left(\frac{m}{s't'}\right)]_R}
\end{tikzcd}
\vspace{-7pt}\]
}

\demo{The diagram commutes since the diagonal map is $D_{st,m/s't'}$ by Lemma \ref{elem-degeneracies}. Let $s_{\infty},t_{\infty}$ be the greatest divisors of $m$ which have the same prime divisors as $s',t'$ respectively. Then $s_{\infty},t_{\infty},r_{\infty} := \frac{m}{s_{\infty}t_{\infty}}$ are pairwise coprime positive integers. By Lemma \ref{degn-compat}, we may replace $\P$ by $\P \times [\Gamma_0(r_{\infty})]_R$ and assume that $m = s_{\infty}t_{\infty}$. Let us denote $\MPp{\P,[\Gamma_0(d)]_R}$ by $\M_0(d)$ for brevity. 

Let $S$ be a $R$-scheme, $x_0 = (E_0/S, \alpha_0,C_0) \in \M_0(m/s't')(S)$ and
\vspace{-7pt}\[x_1=(E_1/S,\alpha_1,C_1) \in \M_0(m/t')(S), \qquad x_2=(E_2/S,\alpha_2,C_2) \in \M_0(m/s')\vspace{-7pt}\] be such that $D_{s,\frac{m}{s't'}}(x_1)=D_{t,\frac{m}{s't'}}(x_2) = x_0$. Thus, there are cyclic isogenies $f_i: E_i \rar E_0$ with $\deg{f_1}=s$, $\deg{f_2}=t$ such that 
\vspace{-7pt}\[f_i(\alpha_i)=\alpha_0, \quad \ker{f_1}=C_1\{s\}, \quad \ker{f_2}=C_2\{t\},\quad f_1(C_1)\left\{\frac{m}{s't'}\right\}=f_2(C_2)\left\{\frac{m}{s't'}\right\} = C_0.\vspace{-7pt}\]
Then $F = E_1 \times_{E_0} E_2$ is an elliptic curve over $S$ by Lemma \ref{rel-coprime-elliptic-curves}, and let $g_i: F \rar E_i$ be the projection. Since the respective degrees of $g_1,g_2$ are $t,s$, they induce a bijection
\vspace{-7pt}\[[\Gamma_0(m)](F/S) \rar [\Gamma_0(s_{\infty})](E_1/S) \times [\Gamma_0(t_{\infty})](E_2/S),\vspace{-7pt}\] 
so $F/S$ has a cyclic subgroup $C$ of degree $m$ such that $g_2(C\{t_{\infty}\}) = C_2\{t_{\infty}\}, g_1(C\{s_{\infty}\})=C_1\{s_{\infty}\}$. 

Let us check that $C\{t\} = \ker{g_1}$. This is an equality of subgroup schemes of degree $t$, so it can be checked after applying $g_2$, since $\deg{g_2}=s$ is prime to $t_{\infty}$. Now, $g_2: \ker{g_1} \rightarrow \ker{f_2}$ is an isomorphism, and one has $\ker{f_2} = C_2\{t\} = g_2(C\{t_{\infty}\})\{t\} = g_2(C\{t\})$, and the conclusion follows. 

Let us check that $g_1(C)\{m/t'\}=C_1$. By prime factorization \cite[Lemma 3.5.1]{KM} and definition of $C$, it suffices to show that $g_1(C)\left\{\frac{t_{\infty}}{t'}\right\}=C_1\left\{\frac{t_{\infty}}{t'}\right\}$. This is an equality of group schemes killed by $t_{\infty}$, so it suffices to check it after applying the isogeny $f_1$ of degree $s$ prime to $t_{\infty}$. One has 
\vspace{-7pt}\[f_1\left(g_1(C)\left\{\frac{t_{\infty}}{t'}\right\}\right)=f_2(g_2(C\{t_{\infty}\}))\left\{\frac{t_{\infty}}{t'}\right\}  = f_2(C_2\{t_{\infty}\})\left\{\frac{t_{\infty}}{t'}\right\} = C_0\left\{\frac{t_{\infty}}{t'}\right\} =f_1\left(C_1\left\{\frac{t_{\infty}}{t'}\right\}\right),\vspace{-7pt}\]
and the conclusion follows. 

Similarly, one has $C[s]=\ker{g_2}$ and $g_2(C)\{m/s'\}=C_2$. Therefore, if $\alpha \in \P(F/S)$ is the inverse image of $\alpha_1$ by $g_1$, then $f_2(g_2(\alpha))=f_1(g_1(\alpha))=f_1(\alpha_1)=\alpha_0=f_2(\alpha_2)$, hence $g_2(\alpha)=\alpha_2$ and one has 
\vspace{-11pt}\[x := (F/S,\alpha,C) \in \M_0(m)(S), \qquad D_{s,\frac{m}{s'}}(x)=x_2, \qquad D_{t,\frac{m}{t'}}(x)=x_1.\vspace{-9pt}\]

Let $y = (G/S,\beta,D) \in \M_0(m)(S)$ such that 
\vspace{-9pt}\[D_{s,m/s'}(y)=x_2, \qquad D_{t,m/t'}(y)=x_1.\vspace{-7pt}\]
There are isogenies $h_i: G \rightarrow E_i$ such that  
\vspace{-9pt}\[h_i(\beta)=\alpha_i, \quad \ker{h_1}=D\{t\}, \quad \ker{h_2}=D\{s\},\quad h_1(D)\left\{\frac{m}{t'}\right\}=C_1, \qquad h_2(D)\left\{\frac{m}{s'}\right\} = C_2.\vspace{-7pt}\]
By \cite[Proposition 6.7.10]{KM}, $f_i \circ h_i: G \rar E_0$ is an isogeny of kernel $H_i$ of degree $st$, and one has $H_1\{t\}=D\{t\}$, $H_2\{s\}=D\{s\}$. Moreover, $h_1$ induces an isomorphism $D\{s\} \rightarrow C_1\{s\}=\ker{f_1}$, so one has $D\{st\} \subset H_1$. Since $D\{st\},\,H_1$ are finite locally free $S$-group schemes of degree $st$, one has $D\{st\}=H_1$ and similarly $D\{st\}=H_2$. Thus $\ker{f_1\circ h_1}=\ker{f_2\circ h_2}$, so one has $f_2\circ h_2 = \psi \circ f_1 \circ h_1$ for some $\psi \in \Aut{E_0/S}$. Moreover, one has $(f_i\circ h_i)(\beta)=f_i(\alpha_i)=\alpha_0$, so the level $N$ structure shows that $\psi(\alpha_0)=\alpha_0$. Since $\P$ is rigid, one has $\psi=\mrm{id}$, so $h := (h_1,h_2): G \rar F$ is an isogeny. One has $\deg{h}= \frac{\deg{h_1}}{\deg{g_1}}=\frac{t}{t}=1$, so $h$ is an isomorphism. One clearly has $h(\beta)=\alpha$, and $(h_1,h_2)(D)$ is a cyclic subgroup $C'$ of degree $m$ such that 
\vspace{-11pt}\[g_1(C'\{s_{\infty}\})=C_1\{s_{\infty}\}=g_1(C\{s_{\infty}\}),\qquad g_2(C'\{t_{\infty}\})=C_2\{t_{\infty}\}=g_2(C\{t_{\infty}\}),\vspace{-4pt}\]
so $C=C'$. Thus $h: (G/S,\beta,D) \rightarrow (F/S,\alpha,C)$ is an isomorphism and $x=y$.
 
}

\bigskip

\defi[AL-autod]{Let $m,d,N > 0$ be integers such that $d\mid m$, $d$ is coprime to $\frac{m}{d}$ and $N$ is coprime to $m$. Let $\P$ be a representable moduli problem over $R$ of level $N$. The \emph{Atkin--Lehner map} $w_d(\P): \MPp{\P,[\Gamma_0(m)]_R} \rightarrow \MPp{\P,[\Gamma_0(m)]_R}$ (denoted $w_d$ if $\P$ is clear) is defined as follows: if $E/S$ is an object of $\Ell_R$, $\alpha \in \P(E/S)$, $C \in [\Gamma_0(m)](E/S)$, let $\pi: E \rightarrow E'$ be the isogeny such that $\ker{\pi}=C\{d\}=C[d]$, and $C' \in [\Gamma_0(m)](E'/S)$ be such that 
\vspace{-9pt}\[C'[d] = \ker{\pi^{\vee}}, \qquad C'[m/d] = \pi(C[m/d]).\vspace{-6pt}\]
We set $w_d(\P)(E/S,\alpha,C) = (E'/S,\pi(\alpha),C')$.}

This definition makes sense, because, since $\pi$ is a cyclic isogeny of degree $d$, so is its dual $\pi^{\vee}$ by \cite[Corollary 5.5.4 (3)]{KM}\footnote{The map $E[d]/C[d] \rightarrow \ker{\pi^{\vee}}$ induced by $\pi$ is a monomorphism of finite locally free $S$-group schemes of degree $d$, hence an isomorphism by \cite[$\text{IV}_3$ (8.11.5)]{EGA}.}. It is clear from the construction that the formation of $w_d$ is natural with respect to $\P$ (for morphisms of moduli problems of level $N$) and commutes to base change and ``level-raising'' in the sense of Lemma \ref{elem-degeneracies}. The following result is an elementary verification.

\lem[AL-elem]{Let $m,d,N$ be positive integers such that $d \mid m$, $d$ is coprime to $\frac{m}{d}$ and $m,N$ are coprime. Let $\P$ be a representable moduli problem over $R$ of level $N$. 
\begin{itemize}[noitemsep,label=\tiny$\bullet$]
\item If $t \mid \frac{m}{d}$, one has $D_{d,t}(\P) = D_{1,t}(\P) \circ w_d(\P)$. 
\item Let $m' \geq 1$ be coprime to $mN$ and $\P'=\P \times [\Gamma_0(m')]_R$. Under the natural isomorphism $\P' \times [\Gamma_0(m)]_R \rightarrow \P \times [\Gamma_0(mm')]_R$ of Lemma \ref{degn-compat}, the following diagram commutes
\vspace{-7pt}\[
\begin{tikzcd}[ampersand replacement=\&]
\MPp{\P',[\Gamma_0(m)]_R} \arrow{r}{\sim} \arrow{d}{w_d(\P')}\& \MPp{\P,[\Gamma_0(mm')]_R} \arrow{d}{w_d(\P)}\\
\MPp{\P',[\Gamma_0(m)]_R} \arrow{r}{\sim}  \& \MPp{\P,[\Gamma_0(mm')]_R}  
\end{tikzcd}
\vspace{-7pt}\]
\end{itemize}
}

\prop[AL-auto]{Let $m,d,d',N$ be positive integers such that $d \mid m$, $d' \mid m$, $d$ is coprime to $\frac{m}{d}$ and $d'$ is coprime to $\frac{m}{d'}$. Let $\P$ be a representable moduli problem of level $N$ coprime to $m$. One has $w_d(\P) \circ w_{d'}(\P) = [\delta] \circ w_{dd'/\delta^2}(\P)$, where $\delta = \mrm{gcd}(d,d')$. In particular, $w_d$ is an automorphism and $w_d(\P)^2=[d]$. }

\demo{Since $w_d(\P)$ commutes to morphisms of moduli problems of level $N$, it is enough to show the identity when $d'=d$ and when $d,d'$ are coprime. By Lemma \ref{AL-elem}, after replacing $\P$ by $\P \times [\Gamma_0\left(\frac{m}{\mrm{lcm}(d,d')}\right)]_R$, we may assume that either $m=d=d'$ or $m=dd'$ with $d,d'$ coprime. When $d=d'$, if $(E/S,\alpha,C) \in \MPp{\P,[\Gamma_0(m)]}(S)$ and $\pi: E \rightarrow E'$ is the isogeny with $\ker{\pi} = C$, then the conclusion holds because $\pi^{\vee}\circ \pi = d \cdot \mrm{id}$. 

Now assume that $d,d'$ are coprime and $m=dd'$. Let $(E/S,\alpha,C) \in \MPp{\P,[\Gamma_0(m)]_R}(S)$. One has $w_{d'}(E,\alpha,C) = (E_1,\alpha_1,C_1)$, where $\pi_1: E \rightarrow E_1$ is the isogeny with kernel $C[d']$, $\alpha_1=\pi_1(\alpha)$, and $C_1 = \pi(C[d])+\ker{\pi_1^{\vee}}$. One has $w_{d}(E_1,\alpha_1,C_1)=(E_2,\alpha_2,C_2)$ where $\pi_2: E_1 \rightarrow E_2$ is the isogeny with kernel $C_1[d]$, $\alpha_2=\pi_1(\alpha_1)$ and $C_2=\pi_2(C_1[d']) + \ker{\pi_2^{\vee}}$. To check that $(E_2,\alpha_2,C_2) = w_{dd'}(\P)(E,\alpha,C)$, it suffices to show that $C=\ker{\pi_2 \circ \pi_1}$ and $C_2 = \ker{\pi_1^{\vee}\circ \pi_2^{\vee}}$. 

Now, $\pi_2 \circ \pi_1$ is an isogeny of degree $dd'$, and $C[d'] \subset \ker{\pi_1}\subset \ker{\pi_2 \circ \pi_1}$; moreover, $\pi_2(\pi_1(C[d]))=\pi_2(C_1[d])=\{0\}$, so $C \subset \ker{\pi_2 \circ \pi_1}$ is an inclusion of finite locally free $S$-subgroup schemes of degree $dd'$, so it is an equality. Similarly, $\pi_2^{\vee}$ vanishes on $C_2[d]$ by construction, and one has 
\vspace{-13pt}\[\pi_1^{\vee}\circ \pi_2^{\vee}(C_2[d'])=\pi_1^{\vee}\circ \pi_2^{\vee}\circ \pi_2(C_1[d']) = \pi_1^{\vee}(dC_1[d'])=\pi_1^{\vee}(C_1[d'])=\{0\},\vspace{-10pt}\]
and we conclude as previously. 
}

\medskip

\cor[flat-deg]{Let $\P$ be a representable moduli problem over $R$ of level $N$, let $m \geq 1$ be coprime to $N$, and let $d,t \geq 1$ be such that $dt \mid m$. Then $D_{d,t}$ is finite locally free, and \'etale above $\Z[1/m]$.}

\demo{This is Zariski-local over $R$, so we may assume that there is a prime $\ell \nmid 2Nm$ invertible in $R$. By Proposition \ref{degn-cartesian} and \cite[Lemmas 036K, 036W]{Stacks}, it is enough to prove the property for $\P \times [\Gamma(\ell)]_R$ instead of $\P$; by Proposition \ref{degn-cartesian}, we may assume that $R=\Z[1/\ell]$ and $\P=[\Gamma(\ell)]_R$. The schemes $\MPp{\P,[\Gamma_0(s)]}$ for $s \mid m$ are then all affine and regular (and every non-empty open subscheme is of dimension two). 
Write $\M := \MP$. 

When $d=1$, $D_{d,t}$ is the ``standard subgroup'' map $[\Gamma_0(m)]_{\cE/\M} \rightarrow [\Gamma_0(t)]_{\cE/\M}$ of finite $\M$-schemes, so it is finite by \cite[II (6.1.5) (v)]{EGA}. Hence, for arbitrary $d,t$, the composition $D_{1,1}\circ D_{d,t} = D_{1,1} \circ w_d \circ D_{1,d}$ is finite, so $D_{d,t}$ is finite by \emph{loc.cit.}. Since $D_{d,t}$ is finite between regular schemes where every non-empty open subscheme is two-dimensional, the discussion in \cite[Notes on Chap. 4]{KM} shows that it is also flat, thus finite locally free. 

The same argument shows that $D_{d,t}$ is \'etale above $\Z[1/m]$ after replacing ``finite'' and \cite[II (6.1.5) (v)]{EGA} by ``\'etale'' and \cite[$\text{IV}_4$ (17.3.5)]{EGA} respectively. }

\bigskip

\bigskip

\section{Compactifications of moduli schemes}
\label{situation-at-infty}

In this section, unless explicitly specified otherwise, $R$ denotes a base ring which is assumed to be regular and excellent (cf. \cite[$\text{IV}_2$ (7.8.2)]{EGA}), and in particular Noetherian. The goal of this section is to state some properties needed to construct compactified moduli schemes and work with morphisms between them. 

\subsection{The general construction}

\lem[normal-near-infinity-equiv]{Let $X$ be a $R$-scheme of finite type. The following are equivalent:
\begin{itemize}[noitemsep,label=\tiny$\bullet$]
\item there is a closed subscheme $Z \subset X$ finite over $R$ such that $X \backslash Z$ is a normal scheme,
\item the reduced subscheme of $X$ attached to its (closed) non-normal locus is finite over $R$,
\end{itemize}
If $X$ is separated, these conditions hold if and only if there is a morphism $\phi: Z \rar X$ of $R$-schemes with $Z$ finite over $R$ such that $X$ is normal at every point outside $\phi(Z)$. }

\demo{Since $R$ is excellent, the normal locus of $X$ is open by \cite[$\text{IV}_2$ (7.8.3)]{EGA}. The claim formally reduces to the following: if $Z$ is a finite $R$-scheme,  $X$ is separated, and $\phi:Z \rightarrow X$ is a $R$-morphism (which is then finite \cite[II (6.1.5) (v)]{EGA}), then $\phi(Z)$ is the underlying subspace of a closed subscheme $H \subset X$ finite over $R$. Indeed, if $H \subset X$ denotes the scheme-theoretic image of $\phi$, the map $Z \rightarrow H$ induced by $\phi$ is finite surjective by \cite[II (6.1.10); I (9.5.4)]{EGA}, so $H$ is affine by \cite[II (6.7.1)]{EGA}. Since $\OO(H)$ injects in $\OO(Z)$, $H$ is finite over $R$, and $X \backslash H$ is normal.}

\defi[normal-near-infinity]{We say that a $R$-scheme of finite type $X$ is \emph{normal near infinity} if there is a closed subscheme $Z \subset X$ finite over $R$.}

In particular, if $(U_i)_i$ is an affine cover of $\Sp{R}$, $X$ is normal near infinity if and only if $X_{U_i}$ is normal near infinity as a $U_i$-scheme for each $i$. 

\lem[preim-of-closed-are-cofinal]{Let $j: X \rar \bA^1_R$ be a finite map of $R$-schemes. Let $Z \subset X$ be a closed subscheme which is finite over $R$. There exists a monic polynomial $f \in R[t]$ such that $H \subset j^{-1}(V(f))$ as closed subschemes of $X$.}

\demo{Note that $X$ is a $R$-scheme of finite type. Let $x$ denote the coordinate on $\bA^1_R$: the function $j^{\sharp}(x)_{|Z} \in \OO(Z)$ is integral over $R$, so there is a monic $f \in R[t]$ such that $f(j^{\sharp}(x)_{|Z})=0$. Then $Z \rightarrow X \overset{j}{\rightarrow} \bA^1_R$ factors through the closed subscheme $V(f)$. }

\prop[one-compactification]{Let $j: X \rar \mathbb{A}^1_R$ be a finite map, where $X$ is normal near infinity. There is a canonical Cartesian diagram 
\vspace{-5pt}\[
\begin{tikzcd}[ampersand replacement=\&]
X \arrow{r}{j}\arrow{d} \& \mathbb{A}^1_R \arrow{d}\\
\overline{X}\arrow{r}{\overline{j}} \&  \mathbb{P}^1_R
\end{tikzcd}
\vspace{-7pt}\]
such that $\overline{j}$ is finite and $\overline{X}$ is normal at every $x \in \overline{X}\backslash X$. 
The map $X\rar \overline{X}$ is an affine open immersion which is scheme-theoretically dominant in the sense of \cite[$\text{IV}_3$ (11.10.2)]{EGA}. 
}

\demo{The construction of $\overline{j}: \overline{X}\rightarrow \bP^1_R$ follows \cite[(8.6.3)]{KM}: it is the relative normalization of the composition $j': X \overset{j}{\rightarrow} \bA^1_R \subset \bP^1_R$. 
Since $j$ is finite, $X \rightarrow \overline{j}^{-1}(\bA^1_R)$ is an isomorphism, hence $X \rightarrow \overline{X}$ is an affine open immersion. Let $t$ be a coordinate on $\bA^1_R$: by Lemma \ref{preim-of-closed-are-cofinal}, there is a monic $f \in R[t]$ such that the invertible locus $X_f$ of $j^{\sharp}(f)$ in $X$ is affine and normal. Since $\bP^1_R$ is also normal, $\overline{X}_f := \overline{j}^{-1}(\bP^1_R \backslash V(f))$ is a normal open subscheme of $\overline{X}$ containing the subset $\overline{X} \backslash X$. 

Let us check that $\overline{j}$ is finite. It is enough to show that, if $V \subset \bP^1_R \backslash V(f)$ is a connected affine open subscheme and $V_1 = V \cap \bA^1_R$, the relative normalization of $X_1=j^{-1}(V_1)$ in $V$ is a finite $U$-scheme. Now, since $X_1$ is a finite normal $V_1$-scheme, it is the disjoint union of finitely many finite $V_1$-schemes that are normal and integral. Since $V_1 \rightarrow V$ is an affine open immersion of normal integral schemes and $R$ is excellent, we are done by \cite[$\text{IV}_2$ (7.8.3) (vi)]{EGA}.

Finally, $\overline{j}$ is finite, so $\overline{X}_f$ is affine. Since $\OO(\overline{X}_f) \rightarrow \OO(X_f)$ is injective by construction, $X_f$ is scheme-theoretically dense in $\overline{X}_f$. Since $\overline{X}$ has the affine open cover $\{X,\overline{X}_f\}$, $X \rightarrow \overline{X}$ is scheme-theoretically dominant.  }

\rem{With $Y=\bA_R^1$, we see in particular that $\OO(U) \rar \OO(U \cap X)$ is injective. }

\rem{If $W \subset X$ is a closed subscheme, finite over $R$, the composition $W \rightarrow X \rightarrow \overline{X}$ is an immersion of proper $R$-schemes, i.e. a closed immersion by \cite[$\text{IV}_3$ (8.11.5)]{EGA}.}

\cor[extension-to-finite]{Let $j: X \rar \mathbb{A}^1_R$ be a finite map with $X$ normal near infinity. Let $Y$ be a finite $R$-scheme, and let $g: X \rar Y$ be a $R$-morphism. Then $g$ extends uniquely to a map $\overline{g}: \overline{X} \rar Y$. }

\demo{The uniqueness follows from Proposition \ref{one-compactification}. For the existence, let $Z \subset \bA^1_R$ be a closed subscheme, finite over $R$, such that $X \backslash j^{-1}(Z)$ is normal. Then $g$ defines a rational map from the normal Noetherian $R$-scheme $\overline{X} \backslash \overline{j}^{-1}(Z)$ to the finite $R$-scheme $Y$. By \cite[II (6.1.13); I (7.2.2.1), (7.2.3)]{EGA}, $g_{|X \backslash j^{-1}(Z)}$ extends to a morphism $g': \overline{X} \backslash \overline{j}^{-1}(Z) \rightarrow Y$, so we obtain $\overline{g}$ by glueing $g$ and $g'$. }

\cor[partially-functorial-compactification]{Keep the notations of Proposition \ref{one-compactification}. Let $R'$ be a $R$-algebra which is a regular excellent ring and $j': X' \rar \mathbb{A}^1_{R'}$ be a finite map. Assume that there is a closed subscheme $Z'$ of $\mathbb{A}^1_R$, finite over $R$, such that\footnote{Note that this condition is stronger than merely requiring that $X'$ be normal near infinity as a $R'$-scheme.} $X'$ is normal at every point not above $Z'$, so that Proposition \ref{one-compactification} defines a compactification $\overline{j'}: \overline{X'} \rightarrow \bP^1_{R'}$. Then any $\bA^1_R$-morphism $f: X' \rightarrow X$ extends uniquely to a $\bP^1_R$-morphism $\overline{f}: \overline{X'} \rar \overline{X}$. }

\demo{The uniqueness of $\overline{f}$ is a consequence of Proposition \ref{one-compactification} applied to $X'$ with $Y=\overline{X}$. To construct $\overline{f}$, after enlarging $Z'$, we may assume that $X \backslash j^{-1}(Z')$ is normal, so $f$ induces (by e.g. \cite[II (6.3.5)]{EGA}) a morphism $f'$ between the relative normalizations of the horizontal maps in the diagram
\vspace{-6pt}\[
\begin{tikzcd}[ampersand replacement=\&]
X' \backslash (j')^{-1}(Z'_{R'}) \arrow{r}{j'}\arrow{d}{f}\& \bP^1_{R'} \backslash Z'_{R'}\arrow{d}\\
X \backslash j^{-1}(Z') \arrow{r}{j} \&  \bP^1_R \backslash Z'
\end{tikzcd}
\vspace{-6pt}\]
We obtain $\overline{f}$ by glueing $f'$ with $f: X' \rar X$ along $X \backslash j^{-1}(Z')$.  
}

\medskip
\noindent
In Corollary \ref{partially-functorial-compactification}, it is natural to ask about the case where $X'=X_{R'}$: when is it normal near infinity? And when this is the case, taking $f$ as the projection $X_{R'} \rightarrow X$, is $\overline{f}$ an isomorphism?

\prop[base-change-by-normal-morphism]{Suppose that $R'$ is a regular excellent ring and that $R \rightarrow R'$ is flat with geometrically normal fibres\footnote{That is, $\Sp{R'} \rightarrow \Sp{R}$ is normal in the sense of \cite[$\text{IV}_2$ (6.8.1)]{EGA}.} (e.g. $R'$ is the localization of a smooth $R$-algebra). Let $f: X \rightarrow \bA^1_R$ be a finite morphism of $R$-schemes with $X$ normal near infinity. Then $X_{R'}$ is normal near infinity (as a $R'$-scheme), and the projection $X_{R'} \rightarrow X$ extends to an isomorphism $\overline{X_{R'}} \rightarrow \overline{X}_{R'}$ (where the compactifications are taken with respect to $f_{R'}: X_{R'} \rightarrow \Sp{R'}$ and to $f$ respectively).  }

\demo{Let $Z \subset \bA^1_R$ be a closed subscheme finite over $R$ such that $X \backslash f^{-1}(Z)$ is normal. By \cite[$\text{IV}_2$ (6.14.1)]{EGA}, $X_{R'} \backslash (f_{R'})^{-1}(Z_{R'})$ is normal. Moreover, $\overline{X}$ (resp. $\overline{X_{R'}}$) is the relative normalization of the affine map $X \overset{f}{\rightarrow} \bA^1_R \rightarrow \bP^1_R$ (resp. its base change to $R'$); since $\Sp{R'} \rightarrow \Sp{R}$ is normal, the conclusion follows from \cite[$\text{IV}_2$ (6.14.5)]{EGA}.}

\bigskip

We now wish to show that the compactification previously constructed is canonical, and that a much wider class of morphisms can be extended.

\lem[kinda-zmt]{Let $R$ be an excellent ring (not necessarily regular) and $f: Z \rar X$ be a finite morphism of $R$-schemes of finite type. Suppose that $Y \subset X$ is a dense open subscheme such that $f^{-1}(Y) \rar Y$ is an isomorphism, $f^{-1}(Y)$ is dense in $Z$, every $\OO_{Z,z}$ for $z \notin f^{-1}(Y)$ is reduced, and every $\OO_{X,x}$ for $x \notin Y$ is normal. Then $f$ is an isomorphism.}

\demo{The question is Zariski-local over $X$, so we may assume that $X$ is connected. Because $f$ is closed and the non-reduced locus of $Z$ is closed, by \cite[$\text{IV}_2$ (7.8.3) (iv)]{EGA}, we may furthermore assume that $X$ is normal (hence integral) and that $Z$ is reduced. Thus $Y$, hence $f^{-1}(Y)$, and therefore its closure $Z$ are irreducible schemes. Hence $Z$ is an integral scheme and the conclusion follows from \cite[II (6.1.15)]{EGA}.}

\medskip

\prop[unique-compactification]{Let $\alpha, \beta: X \rar \bA^1_R$ be two finite maps, where $X$ is a $R$-scheme of finite type normal near infinity. Let $\overline{\alpha}: X_{\alpha} \rar \bP^1_R$, $\overline{\beta}: X_{\beta}\rar \bP^1_R$ be the finite morphisms attached to $\alpha,\beta$ by Proposition \ref{one-compactification}. There is a unique $R$-morphism $\iota: X_{\alpha} \rar X_{\beta}$ such that the diagram
\vspace{-7pt} \[
\begin{tikzcd}[ampersand replacement=\&]
X  \arrow{r}{\mrm{id}}\arrow{d}{\overline{\alpha}}\& X\arrow{d}{\overline{\beta}}\\
X_{\alpha} \arrow{r}{\iota} \&  X_{\beta}
\end{tikzcd}
\vspace{-7pt}\]
commutes. Moreover, $\iota$ is an isomorphism. 
}

\demo{All schemes are Noetherian, so scheme-theoretic images are well-behaved by \cite[I (9.5)]{EGA}. The uniqueness of $\iota$ follows from Proposition \ref{one-compactification}.  

Let $Z$ be the scheme-theoretic image of the diagonal immersion $\Delta: X \rar X \times_R X \rar X_{\alpha} \times_R X_{\beta}$, and $p_1, p_2$ the two projections of $Z$ on $X_{\alpha}$ and $X_{\beta}$ respectively. We will show that $p_1,p_2$ are isomorphisms, so that we can take $\iota=p_2\circ p_1^{-1}$.  

Since $Z$ is a closed subscheme of the proper $R$-scheme $X_{\alpha} \times X_{\beta}$), both $p_i$ are proper. The intersection $Z \cap [X \times X_{\beta}]$ is the scheme-theoretic image of $X \rightarrow X \times X_{\beta}$, which is a section of the first projection $X \times_R X_{\beta} \rightarrow X$ and is therefore a closed immersion. Thus $p_1$ is an isomorphism above the open subscheme $X$, and $p_1^{-1}(X)$ is dense in $Z$. By symmetry, the same claims hold for $p_2$. Therefore, as subspaces of $X_{\alpha} \times X_{\beta}$, one has $Z \backslash \Delta(X) \subset (X_{\alpha} \backslash X) \times (X_{\beta}\backslash X)$. Therefore, both $p_i$ are quasi-finite, hence finite by \cite[$\text{IV}_3$ (8.11.1)]{EGA}. 

Finally, let $W \subset X$ be a closed subscheme, finite over $R$, such that $X \backslash W$ is normal. Then $Z \backslash p_i^{-1}(W)$ is the scheme-theoretic image of $X \backslash W \rightarrow (X_{\alpha} \backslash W) \times_R (X_{\beta} \backslash W)$. Since $X \backslash W$ is reduced, so is $Z \backslash p_i^{-1}(W)$. We conclude by applying Lemma \ref{kinda-zmt} to $p_i$. }

\bigskip

\prop[functorial-compactification]{Let $X,Y$ be two $R$-schemes normal near infinity and admitting finite maps $\alpha,\beta$ to $\bA^1_R$ respectively. Let $\overline{\alpha}: X_{\alpha} \rar \bP^1_R$, $\overline{\beta}: Y_{\beta} \rar \bP^1_R$ be the maps given by Proposition \ref{one-compactification}, where $X$ (resp. $Y$) is an open subscheme of $X_{\alpha}$ (resp. $Y_{\beta}$). Let $f: X \rar Y$ be any finite map of $R$-schemes. Then there is a unique map $\overline{f}: X_{\alpha} \rar Y_{\beta}$ making the following diagram commute:
 \vspace{-7pt}\[
\begin{tikzcd}[ampersand replacement=\&]
X  \arrow{r}{f}\arrow{d}\& Y\arrow{d}\\
X_{\alpha} \arrow{r}{\overline{f}} \&  Y_{\beta}
\end{tikzcd}
\vspace{-7pt}\]
Moreover, $\overline{f}$ is finite, and maps any $x \in X_{\alpha} \backslash X$ into $Y_{\beta} \backslash Y$. 
}  

\demo{By Proposition \ref{one-compactification}, the finite map $\alpha': X \overset{f}{\rar} Y \overset{\beta}{\rar} \bA^1_R$ defines a compactification $X_{\alpha'}$ of $X$. By Corollary \ref{partially-functorial-compactification}, $f$ extends to $\overline{f_{\alpha'}}: X_{\alpha'} \rar Y_{\beta}$, which is a map of finite $\mathbb{P}^1_R$-schemes, so it is finite and maps $X_{\alpha'}\backslash X$ into $Y_{\beta}\backslash Y$. Let $\iota: X_{\alpha} \rar X_{\alpha'}$ be the isomorphism preserving $X$ (thus mapping $X_{\alpha} \backslash X$ into $X_{\alpha'}\backslash X$) given by Proposition \ref{unique-compactification}. Thus $\overline{f}=\overline{f_{\alpha'}}\circ \iota$ satisfies all the requested conditions.  
The uniqueness follows from Proposition \ref{one-compactification} since $Y_{\beta}$ is separated. }

\medskip

As a consequence of the previous results in this section, we obtain: 

\prop[compactification-functor]{Let $\mathbf{ACrv}_R$ be the category whose objects are $R$-schemes of finite type, normal near infinity and admitting a finite map to $\mathbb{A}^1_R$, and whose morphisms are finite maps of $R$-schemes. Let $\mathbf{PCrv}_R$ be the category of proper $R$-schemes that are normal near infinity, its morphisms being finite $R$-maps. There is a \emph{normal compactification functor} $\mathbf{ACrv}_R \rar \mathbf{PCrv}_R$, denoted by $X \longmapsto \overline{X}$, $f \longmapsto \overline{f}$, satisfying the following properties:
\begin{itemize}[noitemsep,label=$-$]
\item Every object $X$ of $\mathbf{ACrv}_R$ has a natural open immersion of $R$-schemes $X \rightarrow \overline{X}$ with dense image,
\item For every object $X$ of $\mathbf{ACrv}_R$, every separated scheme $Y$ and every open subscheme $U \subset \overline{X}$, $Y(U) \rar Y(X \cap U)$ is injective,
\item For any morphism $f: X \rar Y$ in $\mathbf{ACrv}_R$, $\overline{f}$ sends any point in $\overline{X} \backslash X$ into $\overline{Y} \backslash Y$,
\item $\overline{X}$ is normal at every point $x \notin X$, 
\item $\overline{\bA^1_R}=\bP^1_R$. 
\end{itemize} }

\lem[compact-normal]{Let $f: X \rar Y$ be a morphism in $\mathbf{ACrv}_R$. Then $\overline{f}: \overline{X} \rar \overline{Y}$ is the relative normalization of $X \rar Y \subset \overline{Y}$.}

\demo{Let $g: Y \rightarrow \bA^1_R$ be a finite morphism, $\iota_X: X \rightarrow \overline{X}, \iota_Y: Y \rightarrow \overline{Y}$ the compactifications, and apply the final claim in \cite[Lemma 035I]{Stacks} to the diagram
\vspace{-8pt}\[
\begin{tikzcd}[ampersand replacement=\&]
X \arrow{d}{\iota_X} \arrow{r}{\iota_Y \circ f} \& \overline{Y} \arrow{d}{\overline{g}} \\
\overline{X} \arrow{r}{\overline{g\circ f}} \arrow{ur}{\overline{f}}\&\bP^1_R
\end{tikzcd}
\vspace{-10pt}\]} 

\cor[compactification-is-dense]{Let $X$ be an object of $\mathbf{ACrv}_R$ and $f: X \rar \bA^1_R$ be finite. Assume that the compactification $\overline{f}: \overline{X} \rar \bP^1_R$ of $f$ is flat above some open subscheme $V \subset \bP^1_R$ containing $\bP^1_R \backslash \bA^1_R$. Then $X \rightarrow \overline{X}$ is universally scheme-theoretically dominant relatively to $\Sp{R}$ (in the sense of \cite[$\text{IV}_3$ (11.10.8)]{EGA}).} 

\demo{Since $\overline{X} = X \cup \overline{f}^{-1}(V)$ is an open cover, it suffices to show that for any $R$-scheme $S$, the open immersion $(\overline{f}_S)^{-1}(V_S \cap \bA^1_S) \rightarrow (\overline{f}_S)^{-1}(V_S)$ is scheme-theoretically dominant. This immersion is the base change by the flat map $(\overline{f}_S)^{-1}(V_S) \rightarrow \bP^1_S$ of the scheme-theoretically dominant affine open immersion $\bA^1_S \rightarrow \bP^1_S$, so we conclude by \cite[$\text{IV}_3$ (11.10.5)]{EGA}.}

We conclude this section by giving a flatness criterion for a compactified map. 

\lem[dimension-of-fibres]{Let $X$ be an object of $\mathbf{ACrv}_R$ which is nowhere quasi-finite over $S := \Sp{R}$. Then the fibres of $\pi: \overline{X} \rightarrow S$ are pure of dimension one, and their generic points lie in $X$. }

\demo{By \cite[$\text{IV}_3$ (13.1.4)]{EGA}, $\pi$ is not quasi-finite at any point in the (topological) closure of $X$, so $\pi$ is nowhere quasi-finite. Thus, for any $s \in S$, $\overline{X}_s$ has no isolated point. It is finite over $\bP^1_{\kappa(s)}$, so it is pure of dimension one. Since the open subspace $X_s$ of $\overline{X}_s$ is cofinite, it is dense.}

\prop[flatness-cond]{Let $f: X \rightarrow Y$ be a morphism in $\mathbf{ACrv}_R$, where $X,Y$ are nowhere quasi-finite over $\Sp{R}$. Assume that the flatness loci of $X \rightarrow \Sp{R}$ and $Y \rightarrow \Sp{R}$ are dense in every fibre. Assume that $\overline{Y}$ is regular at every point in $\overline{Y} \backslash Y$ and that $\overline{X}$ is Cohen--Macaulay at every point in $\overline{X} \backslash X$. Then $\overline{f}$ is flat above an open neighborhood of $\overline{Y} \backslash Y$.}

\demo{Because $\overline{f}$ is finite hence closed, $\overline{f}^{-1}(\overline{Y} \backslash Y)=\overline{X} \backslash X$, and flat loci are open (\cite[$\text{IV}_3$ (11.1.1)]{EGA}), it is enough to show that $\overline{f}$ is flat at any $x \in \overline{X} \backslash X$. Indeed, let $y=f(x)$ and $s \in \Sp{R}$ be the point above which they lie. By assumption, $p: X \rightarrow S,\, q: Y \rightarrow S$ are flat in a neighborhood of the generic points of $X_s,\,Y_s$, hence $p,\,q$ are open at these points and one has $\dim{\OO_{\overline{X},x}}=\dim{\OO_{\overline{Y},y}}=\dim{\OO_{S,s}}+1$ by \cite[$\text{IV}_3$ (15.4.2)]{EGA} and Lemma \ref{dimension-of-fibres}. The conclusion follows by applying \cite[$\text{IV}_3$ (15.4.2)]{EGA} again to $\overline{f}$ at $x$. 
}

\bigskip
\bigskip

\subsection{Smoothness at infinity}

\defi[smooth-at-infty-gnl]{Let $R$ be a regular excellent ring. An object $Y$ of $\mathbf{ACrv}_R$ is \emph{smooth at infinity} if there exists a closed subscheme $Z \subset Y$, finite over $R$, such that $\overline{Y} \backslash Z \rar \Sp{R}$ is smooth of relative dimension one, and the scheme $(\overline{Y}\backslash Y)^{red}$ is finite \'etale over $\Sp{R}$.}

\lem[smooth-at-infty-dense]{Let $Y$ be a $R$-scheme of finite type, where $R$ is an excellent ring. Let $Z,W \subset Y$ be disjoint closed subschemes, finite over $R$, such that $Y \backslash Z$ is smooth over $R$ of relative dimension one. The quasi-compact open immersion $Y \backslash W \rightarrow Y$ is universally scheme-theoretically dense relatively to $\Sp{R}$.}

\demo{Working Zariski-locally over $Y$, we may assume by \cite[$\text{IV}_3$ (11.10.10)]{EGA} that $R$ is a field $k$ over which $Y$ is smooth of relative dimension one. Hence every connected component of $Y$ is integral and is not contained in $W$, so we can apply \cite[$\text{IV}_3$ (11.10.4)]{EGA}. }

\prop[smooth-at-infty-base-change]{Let $R \rar S$ be a morphism of regular excellent rings, and $Y$ be an object of $\mathbf{ACrv}_R$ which is smooth at infinity. Then $Y_S$ is an object of $\mathbf{ACrv}_S$ which is smooth at infinity, and the canonical morphism $\overline{Y_S} \rar \overline{Y} \times_R \Sp{S}$ is an isomorphism, which also induces an isomorphism $(\overline{Y_S} \backslash Y_S)^{red} \rar \left[(\overline{Y}\backslash Y)^{red}\right]_S$. 
Furthermore, $(\overline{Y} \backslash Y)^{red}$ is a relative effective Cartier divisor on $\overline{Y} \rar \Sp{R}$ in the sense of \cite[Definition 062T]{Stacks}.}

\demo{Let $Z \subset Y$ be a closed subscheme, finite over $R$, such that $\overline{Y} \backslash Z$ is smooth over $R$ of relative dimension one. If $f: Y \rightarrow \bA^1_R$ is finite, then $f_S: Y_S \rightarrow \bA^1_S$ is finite, and $Y_S \backslash Z_S$ is smooth over $S$ hence normal by \cite[$\text{IV}_4$ (17.5.7)]{EGA}. Hence $Y_S$ is an object of $\mathbf{ACrv}_S$. The morphism $c: \overline{Y_S} \rar \overline{Y}_S$ given by Corollary \ref{partially-functorial-compactification} is an isomorphism above the open subscheme $Y_S \subset \overline{Y}_S$. To show that $c$ is an isomorphism, we apply Lemma \ref{kinda-zmt}. Indeed, $c^{-1}(Y_S) \subset \overline{Y_S}$ is a dense open subscheme by construction, and $\overline{Y_S},\overline{Y}_S$ are normal at every point above $\overline{Y} \backslash Y$. Finally, $Y_S \subset \overline{Y}_S$ is a dense open subscheme by Lemma \ref{smooth-at-infty-dense}. 

Therefore, $\overline{Y_S} \backslash Z_S$ is smooth over $S$ of relative dimension one and $c$ induces a surjective closed immersion $c': (\overline{Y_S} \backslash Y_S)^{red} \rar \left[(\overline{Y} \backslash Y)^{red}\right]_S$. Since $Y$ is smooth at infinity, $\left[(\overline{Y} \backslash Y)^{red}\right]_S$ is \'etale over $S$, thus reduced, and $c'$ is an isomorphism. 

For the assertion about relative effective Cartier divisors, since $(\overline{Y}\backslash Y)^{red} \rar \Sp{R}$ is finite \'etale, it is enough by \cite[Lemma 062Y]{Stacks} to prove that for any $R$-algebra $k$ which is a field, $\left[(\overline{Y}\backslash Y)^{red}\right]_k$ is an effective Cartier divisor in $\overline{Y}_k$. By the previous paragraph, we may assume that $R=k$. The claim is then formal by \cite[Lemma 0B8Y]{Stacks}, since the finite \'etale $k$-scheme $(\overline{Y} \backslash Y)^{red}$ is a closed subscheme of the smooth $k$-scheme $\overline{Y} \backslash Z$ of relative dimension one.}

\prop[abhyankar-compactification]{Let $R$ be a regular excellent ring and $Y$ be an object in $\mathbf{ACrv}_R$, smooth at infinity. Let $Z \subset Y$ be a closed subscheme, finite over $\Sp{R}$, such that $Y \rar \Sp{R}$ is smooth outside $Z$. Let $f: X \rar Y$ be a finite morphism of $R$-schemes which is \'etale above $Y \backslash Z$, so that $X$ is an object in $\mathbf{ACrv}_R$.  
If $f$ is tamely ramified in codimension one\footnote{in the sense of \cite[Section 0BSE]{Stacks}.} over $\overline{Y} \backslash Z$, then $\overline{X} \backslash f^{-1}(Z) \rar \Sp{R}$ is smooth of relative dimension one, and $X$ is smooth at infinity. Moreover, $\overline{X} \backslash f^{-1}(Z)\rar \overline{Y} \backslash Z$ is finite locally free. }

\demo{For the first part, since $f^{-1}(Z)$ is finite over $\Sp{R}$, and $X \backslash f^{-1}(Z) \rar Y \backslash Z$ is finite \'etale, $X \backslash f^{-1}(Z)$ is normal, so $X$ is an object of $\mathbf{ACrv}_R$.

By Lemma \ref{compact-normal}, $\overline{f}: \overline{X} \rar \overline{Y}$ is the relative normalization of $X \overset{f}{\rar} Y \rar \overline{Y}$. We now apply Abhyankar's lemma \cite[Lemmas 0EYG, 0EYH]{Stacks}: \'etale-locally on $\overline{Y}\backslash Z$, the map $\overline{f}_Z: \overline{X} \backslash f^{-1}(Z) \rar \overline{Y}\backslash Z$ is given by finitely many copies of maps of the form $A \rar A[x]/(x^e-\alpha)$, where $e$ is an integer invertible in $A$ and $\alpha \in A$ cuts out exactly the inverse image in $A$ of $(\overline{Y}\backslash Y)^{red}$: in particular, $A$ is smooth of relative dimension one over $R$ and $A/\alpha A$ is \'etale over $R$. \footnote{These facts is not stated as such in the reference, but follow directly from the given proof.} Thus $\overline{f}_Z$ is finite locally free by \'etale descent. 

Now, let $A$ be a smooth $R$-algebra of relative dimension one, $\alpha \in A$ and $e \geq 1$ an integer which is invertible in $A$. Let us show that if $A/\alpha A$ is \'etale over $R$, then $A' := A[x]/(x^e-\alpha)$ is smooth of relative dimension one over $R$. Since $A'$ is finite free over $A$, we may assume by \cite[$\text{IV}_2$ (6.8.6)]{EGA} that $R$ is a field $k$. For each $\mathfrak{m} \in \mathrm{Max}(A)$ containing $\alpha$, $A/\alpha A$ is reduced, so $\alpha$ is a uniformizer at $\mathfrak{m}$ and $A'_{\mathfrak{m}}$ is a discrete valuation ring by \cite[Lemma 09EV]{Stacks}; $A'[\alpha^{-1}]$ is \'etale over $A[\alpha^{-1}]$ so is normal, and we are done.  

As a consequence, \'etale descent (e.g. \cite[Lemma 036U]{Stacks}) shows that $\overline{X} \backslash f^{-1}(Z)$ is smooth over $R$. 
The local description also implies that \'etale-locally on $\overline{Y} \backslash Z$, $(\overline{X} \backslash X)^{red}$ is given by copies of $A[x]/\sqrt{(x^e-\alpha,\alpha A)}=A[x]/(x) \simeq A$. Thus, $(\overline{X} \backslash X)^{red}$ is smooth over $\Sp{R}$ by descent. It is also finite over $\Sp{R}$, hence it is finite \'etale.
}

\rem{The tame ramification assumption in Proposition \ref{abhyankar-compactification} is automatically verified if $R$ is a product of domains of characteristic zero.} 

\subsection{Compactified moduli schemes}

Recall that the base rings are always assumed to be regular excellent, unless explicitly specified otherwise. 

\defi[normal-near-infinity-mp]{A finite representable moduli problem $\P$ over $R$ is \emph{normal near infinity} if $\MP$ is normal near infinity in the sense of Definition \ref{normal-near-infinity}. 

In this case, we define the \emph{compactified moduli scheme} $\CMP$ as the compactification of $\MP$ in the sense of Proposition \ref{compactification-functor}. It is endowed with a finite morphism $\overline{j}: \CMP \rar \bP^1_R$ extending the usual $j$-invariant. Its \emph{cuspidal subscheme} $\CP$ is the reduced scheme attached to the closed subspace $\{\overline{j}=\infty\}$ of $\CMP$. }

\rem{If $\P$ is a finite representable moduli problem over $R$, its $j$-invariant map $j: \MP\rightarrow \bA^1_R$ is finite by \cite[Proposition 8.2.2]{KM}, so its compactification makes sense. Moreover, by Lemma \ref{preim-of-closed-are-cofinal} and the proof of Proposition \ref{compactification-functor}, Definition \ref{normal-near-infinity-mp} agrees with the notions developed in \cite[(8.6.2),(8.6.3)]{KM}.}

\medskip

By Corollary \ref{extension-to-finite}, one has:
	
\cor[finite-extension-mp]{Let $\P$ be a finite representable moduli problem over $R$, normal near infinity. Let $Y$ be a finite $R$-scheme. Any morphism $\MP \rar Y$ of $R$-schemes extends uniquely to $\CMP \rar Y$.}

\smallskip

By Proposition \ref{compactification-functor}, one has:

\prop[finite-maps-extend]{Let $\P,\P'$ be finite representable moduli problems over $R$ normal near infinity. Let $f: \MP \rar \MPp{\P'}$ be any finite map of $R$-schemes. Then $f$ extends to a finite map $\CMP \rar \CMPp{\P'}$.}

\medskip

For later reference, we give a criterion to detect when such extensions are flat in a neighborhood of the cuspidal subscheme.

\prop[flatness-at-cusps]{In the situation of Proposition \ref{finite-maps-extend}, assume that $\P,\P'$ are flat over $R$ and that $\CMP$, $\CMPp{\P'}$ are respectively Cohen--Macaulay and regular at every point of the cuspidal subscheme. Then $\overline{f}$ is flat on an open neighborhood of the cuspidal subscheme of $\CMP$, and the compactified $j$-invariants of $\P,\P'$ are finite locally free. 
}

\demo{By Corollary \ref{over-Ell-concrete}, the moduli schemes have finite flat $j$-invariants, so they are flat and nowhere quasi-finite over $R$, and we can apply Proposition \ref{flatness-cond}.}

\medskip

\defi[normal-near-infinity-relative]{If $R'$ is a $R$-algebra (which is a regular excellent ring) and $\P$ is a finite representable moduli problem over $R'$ with $j$-invariant $j: \MP \rightarrow \bA^1_{R'}$, we say that $\P$ is \emph{normal near infinity over $R$} if there exists a closed subscheme $Z$ of $\bA^1_R$, finite over $R$, such that $\MP$ is normal outside $j^{-1}(Z_{R'})$.}

\cor[base-changed-morphisms]{Let $R'$ be a $R$-algebra (which is a regular excellent ring) and let $\P,\P'$ be finite representable moduli problems over $R,R'$ respectively that are both normal near infinity over $R$. Let $f:\MPp{\P'} \rar \MP_{R'}$ be induced by a morphism of moduli problems $\P' \Rightarrow \P_{R'}$. Then $f$ extends to a $\bP^1_{R'}$-morphism $\CMPp{\P'} \rar \CMP_{R'}$, which sends $\CPp{\P'}$ to $ (\CP)_{R'}$. }

\demo{This follows from Corollary \ref{partially-functorial-compactification}, the finite maps to $\bA^1$ being the $j$-invariants.}

\defi[formation-commutes-base-change]{Suppose that we can take $\P'=\P_{R'}$ in Corollary \ref{base-changed-morphisms}. The morphisms $\overline{m}: \CMPp{\P'} \rar \CMP_{R'}$ and $c: \CPp{\P'} \rar (\CP)_{R'}$ constructed by this Corollary are called the \emph{base change morphisms} for $\P$ and $R \rar R'$. We then say that the formation of the compactified moduli scheme (resp. the cuspidal subscheme) for $\P$ \emph{commutes with the base change $R \rar R'$} if $\overline{m}$ (resp. $c$) is an isomorphism.}

\rem{In the situation of Definition \ref{formation-commutes-base-change}, $\overline{m}$ is finite and an isomorphism above $\overline{j}^{-1}(\mathbb{A}^1_{R'})$.}

\bigskip

\prop[normal-near-infinity-ascends]{(cf. \cite[Proposition 8.6.7]{KM}) Let $R \rar R'$ be a morphism of $R_0$-algebras (where $R_0,R,R'$ are regular excellent) and $\P$ be a moduli problem over $R$ normal near infinity over $R_0$. Suppose that $\Sp{R'} \rightarrow \Sp{R}$ is normal. Then $\P_{R'}$ is normal near infinity over $R_0$. Moreover, the formation of the compactified moduli scheme and the cuspidal subscheme commutes with the base change $R \rar R'$.}

\demo{By Proposition \ref{base-change-by-normal-morphism}, $\P_{R'}$ is normal near infinity over $R_0$ and the base change map $\overline{m}: \CMPp{\P_{R'}} \rightarrow \CMP_{R'}$ is an isomorphism of $\bP^1_{R'}$-schemes. Since $\CP$ is reduced and $R \rightarrow R'$ is flat with reduced fibres, $(\CP)_{R'}$ is reduced by \cite[$\text{IV}_2$ (6.8.2), (3.3.5)]{EGA}. Since the base change map $c: \CPp{\P_{R'}} \rar (\CP)_{R'}$ is a surjective closed immersion, it is an isomorphism.}

\medskip

\prop[normal-near-infinity-descends]{Let $R \rar R'$ be a faithfully flat morphism of $R_0$-algebras (with $R_0,R,R'$ regular excellent) and $\P$ a finite representable moduli problem over $R$. If $\P_{R'}$ is normal near infinity over $R_0$, then $\P$ is normal near infinity over $R_0$. If $R \rightarrow R'$ is also finite and $\P_{R'}$ is normal near infinity, so is $\P$. }

\demo{Let $Z \subset \bA^1_{R_0}$ be a closed subscheme, finite over $R_0$, such that $\MPp{\P_{R'}} \backslash j^{-1}(Z_{R'})$ is normal. Then $\MPp{\P_{R'}} \backslash j^{-1}(Z_{R'}) \rar \MP \backslash j^{-1}(Z_R)$ is faithfully flat from a normal scheme, hence its target is normal by \cite[$\text{IV}_2$ (6.5.4)]{EGA}. Thus, $\P$ is normal near infinity over $R_0$. 

Suppose now that $\P_{R'}$ is normal near infinity and $R'$ is finite faithfully flat over $R$: it suffices by the previous paragraph to show that $\P_{R'}$ is normal near infinity over $R$. Since $R'$ is finite over $R$, the $j$-invariant map $\MPp{\P_{R'}} \rightarrow \bA^1_R$ is finite, and we are done by Lemma \ref{preim-of-closed-are-cofinal}. } 

\medskip

Using these results, we can compactify certain constructions that we previously made for moduli problems with level structure.

\prop[normal-near-inf-level]{Let $\P$ be a finite representable moduli problem over a regular excellent ring $R$ endowed with a level $N$ structure for some $N \geq 1$. Assume that $\P$ is normal near infinity. Then the action of $(\Z/N\Z)^{\times}$ on $\MP$ extends to an action of $(\Z/N\Z)^{\times}$ on $\CMP$, which is natural with respect to $\P$ (equipped with its level structure).}  

\demo{This is a direct application of Proposition \ref{compactification-functor}.}

\lem[normal-near-inf-gamma0]{Let $\P$ be a finite representable moduli problem over a regular excellent ring $R$ and let $m \geq 1$. Suppose that $\P \times [\Gamma_0(m)]_R$ is normal near infinity. Then, for any $d \mid m$, $\P \times [\Gamma_0(d)]_R$ is normal near infinity.}

\demo{Let $Z$ be a closed subscheme of $X := \MPp{\P,[\Gamma_0(m)]_R}$, finite over $R$, outside which $X$ is normal. If $d \mid m$, the map $D_{1,d}: X \rar Y := \MPp{\P,[\Gamma_0(d)]_R}$ is finite locally free, so $Y$ is normal outside the image of $D_{1,d}: Z \rightarrow Y$. }

By Proposition \ref{compactification-functor}, this implies:

\prop[degn-AL-infty]{Let $\P$ be a finite representable moduli problem over a regular excellent ring $R$ endowed with a level $N$ structure. Let $m \geq 1$ be coprime to $N$ such that $\P \times [\Gamma_0(m)]_R$ is normal near infinity. Then:
\begin{itemize}[noitemsep,label=$-$]
\item for any $d,t \geq 1$ such that $dt \mid m$, the morphism $D_{d,t}$ of Definition \ref{degeneracies} extends to a finite morphism $\CMPp{\P,[\Gamma_0(m)]_R} \rar \CMPp{\P,[\Gamma_0(t)]}$ of $R$-schemes, also denoted $D_{d,t}$. 
\item for any $d \geq 1$ dividing $m$ and coprime to $\frac{m}{d}$, the automorphism $w_d$ of $\MPp{\P,[\Gamma_0(m)]_R}$ of Definition \ref{AL-autod} extends to an automorphism of $R$-schemes of $\CMPp{\P,[\Gamma_0(m)]_R}$, also denoted $w_d$.
\end{itemize}
}

\rem[degn-AL-infty-natural]{In the situation of Proposition \ref{degn-AL-infty}, the formation of the $D_{d,t}$ and $w_d$ is natural with respect to $\P$ in the obvious sense and respects ``level-raising''. Moreover, the relations between the various $w_d$, $D_{s,t}$, $[n]$ (for $n \in (\Z/N\Z)^{\times}$) and canonical isomorphisms proved throughout Section \ref{gamma0m-structure} still hold.\footnote{Namely, this applies to Lemmas \ref{elem-degeneracies}, \ref{degn-compat} and \ref{AL-elem}, Proposition \ref{AL-auto}, and the fact that the diagrams in Propositions \ref{degn-cartesian} and \ref{two-degn-cartesian} commute. } }

\cor[flat-deg-infty]{Consider the situation of Proposition \ref{degn-AL-infty}, where we furthermore assume that $\P$ is locally free over $\Ell_R$ and that, for any $s \mid m$ and every point $x \in X:=\CMPp{\P,[\Gamma_0(s)]_R}$ lying in the cuspidal subscheme, the local ring $\OO_{X,x}$ is regular. Then the $D_{d,t}$ defined in Proposition \ref{degn-AL-infty} are finite locally free.}

\demo{The map $D_{d,t}$ is finite, and it is flat at every point outside the cuspidal subscheme by Corollary \ref{flat-deg}. Since the moduli problems $\P$ and $[\Gamma_0(s)]_R$ are finite locally free, the regularity assumption and Proposition \ref{flatness-at-cusps} imply that $D_{d,t}$ is flat.}

\bigskip

In order to apply Proposition \ref{degn-AL-infty}, we need concrete criteria in order to decide when the moduli problem $\P \times [\Gamma_0(m)]_R$ is normal near infinity and find out whether, for instance, the compactified moduli scheme is regular at points in the cuspidal subscheme. It will be convenient to work with the following strengthening of this notion:

\defi[smooth-at-infty-moduli]{Let $R'$ be a $R$-algebra which is a regular excellent ring. A finite representable moduli problem $\P$ over $R'$ is \emph{smooth at infinity over $R$} if the following properties hold:
\begin{itemize}[noitemsep,label=$-$]
\item $\P$ is normal near infinity.
\item There exists a closed subscheme $Z \subset \bA^1_R$, finite over $R$, such that $\CMP \rar \Sp{R'}$ is smooth outside $j^{-1}(Z_{R'})$.  
\item $\CP$ is finite \'etale over $R'$.
\end{itemize}
When $R=R'$, we say that $\P$ is smooth at infinity. 
}

\rem{Suppose that $R=R'$. Then Definition \ref{smooth-at-infty-moduli} is the same as \cite[Definition 8.6.5]{KM}. Since the $j$-invariant map $\MP \rar \bA^1_R$ is finite, it is also the same as asking that $\MP \rar \Sp{R}$ is smooth at infinity in the sense of Definition \ref{smooth-at-infty-gnl}. }

\medskip

Since a scheme which is smooth over a regular base is regular, hence normal, one has:

\lem[smooth-at-infinity-implies-normal]{Let $R'$ be a $R$-algebra. If $\P$ is a moduli problem over $R'$ which is smooth at infinity over $R$, then $\P$ is normal near infinity over $R$.}

\cor[flatness-at-cusps-smooth]{Let $\P,\P'$ be finite flat representable moduli problems over $R$, smooth at infinity, and $f: \MP \rightarrow \MPp{\P'}$ be a finite flat morphism. Then the extension $\overline{f}: \CMP \rightarrow \CMPp{\P'}$ is finite locally free, and so are the compactified $j$-invariants of $\P,\P'$.}

\demo{This follows from Proposition \ref{flatness-at-cusps}.}

\medskip

\prop[smooth-at-infinity-base-change]{(cf. \cite[Proposition 8.6.6]{KM}) Let $R \rar R'$ be a morphism of $R_0$-algebras (with $R_0,R,R'$ regular excellent). Let $\P$ be a finite representable moduli problem over $R$ which is smooth at infinity over $R_0$. Then $\P_{R'}$ is smooth at infinity over $R_0$, and the formation of the compactified moduli scheme and the cuspidal subscheme for $\P$ commutes with the base change $R \rar R'$. Moreover, the cuspidal subscheme of $\P$ is a relative effective Cartier divisor in $\CMP$.} 

\demo{This is an immediate consequence of Proposition \ref{smooth-at-infty-base-change}.}

\medskip

\prop[smooth-at-infinity-descends]{Let $R \rar R'$ be a morphism of $R_0$-algebras (with $R_0,R,R'$ regular excellent) such that $\Sp{R'} \rightarrow \Sp{R}$ is normal and surjective. Let $\P$ be a finite representable moduli problem over $R$ such that $\P_{R'}$ is smooth at infinity over $R_0$. Then $\P$ is smooth at infinity over $R_0$ and the formation of the compactified moduli scheme and the cuspidal subscheme for $\P$ commutes with the base change $R \rar R'$.}

\demo{Let $Z \subset \bA^1_{R_0}$ be a closed subscheme, finite over $R_0$, such that $\CMPp{\P_{R'}} \backslash j^{-1}(Z_{R'})$ is smooth of relative dimension one over $\Sp{R'}$. Now, $\P$ is normal near infinity over $R_0$ by Proposition \ref{normal-near-infinity-descends}, and the formation of the compactified moduli scheme and the cuspidal subscheme commutes with base change by Proposition \ref{normal-near-infinity-ascends}. In particular, after making the faithfully flat base change $R \rar R'$, $\CMP \backslash j^{-1}(Z_R) \rar \Sp{R}$ becomes smooth of relative dimension one and $\CP \rar \Sp{R}$ becomes finite \'etale. By descent, $\CMP \backslash j^{-1}(Z_R)$ is smooth over $R$ of relative dimension one and $\CP$ is finite \'etale over $R$.}

\medskip

\prop[struct-gamma0m-smooth]{Let $\P$ be a finite representable moduli problem over $R$, where $R$ is a Noetherian ring (without any extra assumptions). Let $m_0,m_1 \geq 1$ be coprime integers such that $m_0$ is square-free and $m_1$ is invertible in $R$. Let $m=m_0m_1$ and $Z \subset \bA^1_R$ be a closed subscheme finite over $R$ such that $\P$ is \'etale over $\Ell_R$ above $U := \bA^1_R \backslash Z$. Then:
\begin{itemize}[noitemsep,label=$-$]
\item The morphism $f: j^{-1}(U) \subset \MPp{\P,[\Gamma_0(m)]_R]} \rar \Sp{R}$ is flat, and its fibres are Cohen--Macaulay, pure one-dimensional and geometrically reduced. Moreover, $f$ is smooth of relative dimension one, except at supersingular points in residue characteristic $p \mid m_0$. 
\item Assume furthermore that $R$ is a regular excellent domain of characteristic zero. Then $\P \times [\Gamma_0(m)]_R$ is smooth at infinity and $\CMPp{\P,[\Gamma_0(m)]_R} \rar \CMP$ is finite locally free. 
\item If furthermore $R'$ is a $R$-algebra which is a regular excellent ring, the conclusions of the previous point hold after replacing $(R,\P)$ with $(R',\P_{R'})$. 
\end{itemize}}

\demo{The last point follows from the second one using Proposition \ref{smooth-at-infinity-base-change}. We start with the first point. 

The morphism $j^{-1}(U) \subset \MPp{\P,[\Gamma_0(m)]_R]} \rar \Sp{R}$ is flat, smooth above $R[1/m]$, and its fibres are Cohen--Macaulay and pure one-dimensional by Corollary \ref{over-Ell-concrete}. To prove the rest of the claim, we may assume by \cite[$\text{IV}_2$ (6.8.6)]{EGA} that $R$ is an algebraically closed field $k$ of prime characteristic $p \mid m_0$. Replacing $(\P,m)$ with $(\P \times [\Gamma_0(\frac{m}{p})]_k,p)$, we may assume that $m=m_0=p$. As in the proof of Corollary \ref{flat-deg}, one reduces to the case where $\P=[\Gamma(\ell)]_k$ for some odd prime $\ell$ invertible in $k$, which is well known by \cite[(13.1.4), Theorem 13.4.7, (13.5.6)]{KM}. 

We are now interested in compactifications: we assume from now on that $R$ is a regular excellent domain of characteristic zero and $m=m_0m_1$ with $m_0$ square-free and $m_1$ coprime to $m_0$ and invertible in $R$. If we can show that $\P \times [\Gamma_0(m)]_R$ is smooth at infinity, then in particular the proof shows that $\P$ is smooth at infinity. Since the forgetful map $\MPp{\P,[\Gamma_0(m)]_R} \rightarrow \MP$ is flat, its compactification is flat by Corollary \ref{flatness-at-cusps-smooth}. 

Replacing $(\P,m_0,m_1)$ by $(\P \times [\Gamma_0(m_1)]_R,m_0,1)$, we may assume that $m=m_0$. The claim is Zariski-local over $R$, so we may assume that there is a prime $\ell \nmid 2m$ which is invertible in $R$. Let $Y_0(m)^R$ be the coarse moduli scheme over $R$ attached to $[\Gamma_0(m)]_R$ (cf. \cite[(8.1.1)]{KM}). After extending $Z$ if needed, we may assume that it is defined by a monic polynomial $f \in R[t]$ divisible by $t(t-1728)\prod_{p \mid m}{f_p(t)}$, where $f_p(t) \in \Z[t]$ is a monic lift of the supersingular polynomial in $\F_p[t]$. 

By \cite[Theorem 10.10.3 (5)]{KM}, $Y_0(m)^{\Z}$ is smooth at infinity in the sense of Definition \ref{smooth-at-infty-gnl} (by \cite[Theorem 7.4.2 (4)]{KM}, for $p \mid m$, $[\Gamma_0(p)]$ is the quotient of $[\Gamma(p)]$ by the balanced subgroup of lower-triangular\footnote{Recall that we use a different convention for the $\GL{\Z/N\Z}$-action on $[\Gamma(N)]$ than \emph{op.cit.}.} matrices in $\GL{\F_p}$). Since $R$ is flat over $\Z$, one has $Y_0(m)^R = (Y_0(m)^{\Z})_R$ by \cite[Proposition 8.1.6]{KM}, so $Y_0(m)^R$ is smooth at infinity by Proposition \ref{smooth-at-infty-base-change}. 

By Proposition \ref{abhyankar-compactification}, it suffices to show that the map $\pi: X := \MPp{\P \times [\Gamma_0(m)]_R} \rightarrow Y_0(m)^R$ (see \cite[(8.1.4)]{KM}) is finite \'etale above $j^{-1}(U)$. Arguing again as in the proof of Corollary \ref{flat-deg}, it is enough to show this when $\P=[\Gamma(\ell)]_R$. In this case, $\pi$ is the quotient by $\GL{\F_{\ell}}$, so it is enough to show by \cite[Theorem A7.1.1]{KM} that $-\mathrm{id}$ acts trivially on $X$ and that for any $g \in \GL{\F_{\ell}} \backslash \{\pm \mathrm{id}\}$, the stabilizer of $g$ in $Y(\ell)_R=\MP$ does not meet $j^{-1}(U)$. Both are easy to check, because the function $j(j-1728)$ is invertible on $U$, so elliptic curves corresponding to points in $j^{-1}(U) \subset \MP$ have no non-trivial automorphisms by \cite[(8.4.2)]{KM}. }

\cor[smooth-at-infinity-gamma0m]{Let $\P$ be a finite \'etale representable moduli problem over a regular excellent domain $R$. Let $m_0,m_1 \geq 1$ be coprime, where $m_0$ is square-free and $m_1 \in R^{\times}$. Then $\P \times [\Gamma_0(m_0m_1)]_R$ is normal near infinity over $\Z$. If $R$ has characteristic zero, it is smooth at infinity over $\Z$. }

\section{Application to our moduli problems}
\label{application-moduli}

\subsection{Representability results}

\defi{Let $G$ be a finite locally free commutative group scheme over a ring $R$ and $N \geq 1$ be an integer. The moduli problem $\P_G(N)$ over $R$ is defined as follows: 
\begin{itemize}[noitemsep,label=$-$]
\item if $E/S$ is an object of $\Ell_R$, $\P_G(N)(E/S)$ is the set of isomorphisms $\iota: G \times_R S \rar E[N]$ of $S$-group schemes, 
\item for a morphism $F/T \rar E/S$ in $\Ell_R$, the morphism $\P_G(N)(E/S) \rar \P_G(N)(F/T)$ is induced by the isomorphism $F[N] \rar E[N] \times_S T$. 
\end{itemize}
It is endowed with a natural structure of moduli problem of level $N$. In particular, for any $m \geq 1$ coprime to $N$, the moduli problem $\P_G(N,\Gamma_0(m)) := \P_G(N) \times [\Gamma_0(m)]_R$ is endowed with a natural level structure of level $Nm$. 
}

\rem{Let $R$ be a $\Z[1/N]$-algebra and $G$ be the constant group scheme $(\Z/N\Z)^{\oplus 2}_R$. Let $e_1,e_2 \in G(R)$ be the standard basis of sections, then, for any object $E/S$ of $\Ell_R$, 
\vspace{-7pt}\[f \in \P_G(N)(E/S) \mapsto (f(e_1),f(e_2)) \in [\Gamma(N)]_R(E/S)\vspace{-7pt}\] 
defines an isomorphism of functors $\P_G(N) \rar [\Gamma(N)]_R$. 
}

\medskip

We will use the following classical result:

\lem[hom-finite-free-is-representable]{Let $S$ be a scheme and $G, H$ be finite locally free $S$-schemes. Then the functors $\Sch_S \rightarrow \Set$
\vspace{-7pt}\[\underline{\mathrm{Hom}}(G,H): T \mapsto \mrm{Hom}_{\Sch_T}(G_T,H_T),\qquad \underline{\mathrm{Isom}}(G,H): T \mapsto \mrm{Isom}_{\Sch_T}(G_T,H_T)\vspace{-7pt}\] are representable by affine $S$-schemes of finite presentation. If moreover $G,H$ are $S$-group schemes, the functors $\Sch_S \rightarrow \Set$
\vspace{-7pt}\[\underline{\mathrm{Hom}}_{\Gp}(G,H): T  \mapsto \mrm{Hom}_{\GpSch_T}(G_T,H_T),\qquad \underline{\mathrm{Isom}}_{\Gp}(G,H): T \mapsto \mrm{Isom}_{\GpSch_T}(G_T,H_T)\vspace{-7pt}\] are also representable by affine $S$-schemes of finite presentation.
}

\demo{If $T$ is a $S$-scheme, one has 
\vspace{-7pt}\[\mrm{Hom}_T(G_T,H_T) \simeq \mrm{Hom}_G(G \times_S T, H \times_S G) \simeq \mrm{Hom}_S(T,\mrm{Res}_{G/S}(H \times_S G)),\vspace{-7pt}\]
where $\mrm{Res}_{G/S}$ denotes the Weil restriction along the finite locally free map $G \rightarrow S$, which exists for $H \times_S G \rightarrow G$ by \cite[\S 7.6 Theorem 4]{BLR}. By the proof of \emph{loc.cit.} and \cite[\S 7.6 Proposition 5]{BLR}, $\mrm{Res}_{G/S}(H \times_S G)$ is an affine $S$-scheme of finite presentation. The three other functors can be constructed using $\underline{\mathrm{Hom}}$ and fibre products (which proves that they can be represented by affine $S$-schemes of finite presentation) as follows. For any $S$-scheme $T$, there is a natural bijection
\vspace{-7pt}\[\mrm{Isom}_T(G_T,H_T) \simeq \{(f,g) \in \underline{\mathrm{Hom}}(G,H)(T) \times \underline{\mathrm{Hom}}(H,G)(T) \; \mid \; f \circ g = \mathrm{id}_{H_T},\, g \circ f = \mathrm{id}_{G_T}\}. \vspace{-7pt}\] 
If $G,H$ are group schemes with multiplications $m_G,m_H$, $\underline{\mathrm{Hom}}_{\Gp}(G,H)$ is the equalizer of the two maps $\underline{\mathrm{Hom}}(G,H) \rightarrow \underline{\mathrm{Hom}}(G \times G,H)$ given by 
\vspace{-7pt}\[f \mapsto f \circ m_G,\qquad f \mapsto m_H \circ (f \times f),\vspace{-7pt}\] 
and one has clearly $\underline{\mathrm{Isom}}_{\Gp}(G,H) \simeq \underline{\mathrm{Hom}}_{\Gp}(G,H) \times_{\underline{\mathrm{Hom}}(G,H)} \underline{\mathrm{Isom}}(G,H)$.}

\cor[iso-over-finite-extension]{Let $R$ be a ring and $G,H$ be finite locally free commutative $R$-group schemes. Assume that $\underline{\mathrm{Aut}}(H) := \underline{\mathrm{Isom}}_{\Gp}(H,H)$ is a flat $R$-scheme. The following are equivalent:
\begin{enumerate}[noitemsep,label=(\roman*)]
\item \label{iofe-1} for every prime ideal $\mfk{p}$ of $R$, there is a faithfully flat $R_{\mfk{p}}$-algebra $R'$ such that $G_{R'}$ and $H_{R'}$ are isomorphic, 
\item \label{iofe-2} there is a faithfully flat quasi-finite $R$-algebra of finite presentation $R'$ such that $G_{R'}$ and $H_{R'}$ are isomorphic,
\item \label{iofe-3} $\underline{\mathrm{Isom}}_{\Gp}(G,H)$ is a faithfully flat $R$-scheme. 
\end{enumerate}
}

\demo{Clearly, \ref{iofe-2} implies \ref{iofe-1}. Since $\underline{\mathrm{Isom}}_{\Gp}(G,H)$ is an affine $S$-scheme of finite presentation, \ref{iofe-3} implies \ref{iofe-2} by \cite[$\text{IV}_3$ (17.16.2)]{EGA}. Now assume \ref{iofe-1} and let $\mathfrak{p} \subset R$ be a prime ideal. There is a faithfully flat $R_{\mfk{p}}$-algebra $R'$ such that $G_{R'} \simeq H_{R'}$, so $\underline{\mathrm{Isom}}_{\Gp}(G,H)_{R'} \simeq \underline{\mathrm{Aut}}(H)_{R'}$ is faithfully flat over $R'$. By descent, $\underline{\mathrm{Isom}}_{\Gp}(G,H)_{R_{\mfk{p}}}$ is faithfully flat over $R_{\mfk{p}}$, so \ref{iofe-3} holds. 
}

\prop[PG-representable-when-good]{Let $G$ be a finite locally free commutative group scheme over a ring $R$ and $N$ be a positive integer. Then the moduli problem $\P_G(N)$ over $R$ is relatively representable and affine of finite presentation. If $N \geq 3$, then $\P_G(N)$ is representable.}

\demo{By Lemma \ref{hom-finite-free-is-representable}, $\P_G(N)$ is relatively representable and affine of finite presentation. When $N \geq 3$, it is rigid by \cite[Corollary 2.7.2]{KM}, hence representable by Proposition \ref{relrep-rigid-implies-rep}.}

\rem{By \cite[Propostion 1.7.2]{KM}, if $G$ is a finite locally free commutative $R$-group scheme killed by $ab$, where $a,b \geq 1$ are coprime integers, then $\P_G(ab)$ identifies (as a moduli problem of level $ab$) with the product $\P_{G[a]}(a) \times \P_{G[b]}(b)$, and $G[a],G[b]$ are finite locally free commutative $R$-group schemes. 

If $G/\Sp{R}$ is in fact an object in the category $\mathbf{FGp}^{N}_{\Z}$ of Definition \ref{category-strongly-free-gpsch}, then, for any divisor $d \mid N$, there is a forgetful morphism $\P_G(N) \rightarrow \P_{G[d]}(d)$ of moduli problems of level $N$.  }

\subsection{The \'etale case}
\label{subsect:etale-case}

If $R$ is a $\Z[1/N]$-algebra and $G$ is a finite locally free $\Z[1/N]$-group scheme, then $G$ is \'etale. Let us assume for the sake of simplicity that $R$ is connected. Then, there is a finite abelian group $A$ such that \'etale-locally on $R$, $G$ is isomorphic to the constant group scheme given by $A$. If there is a $R$-scheme $S$ and an elliptic curve $E/S$ such that $E[N] \simeq G_S$, it follows that $A$ is isomorphic to $(\Z/N\Z)^{\oplus 2}$. This explains the following definition. 

\defi[torsion-group]{Let $N \geq 1$ be an integer, $S$ be a $\Z[1/N]$-scheme. An \emph{\'etale $N$-torsion group over $S$} is a group scheme $G$ over $S$ which is \'etale-locally isomorphic to $(\Z/N\Z)^{\oplus 2}$. In particular, a $N$-torsion group is commutative, finite locally free and \'etale over $S$.}

Suppose that $G$ is an \'etale $N$-torsion group over some $\Z[1/N]$-scheme $S$ and assume that we have an isomorphism $\varphi: G \rar E[N]$ of $S$-group schemes for some elliptic curve $E/S$. Then $G \times_S G \overset{\varphi \times \varphi}{\rar} E[N] \times E[N] \overset{\mrm{We}}{\rar} (\mu_N)_S$ is a surjective bilinear alternating morphism between \'etale $S$-schemes. 

\lem[pol-torsor-exists]{Let $S$ be a scheme and $G$ be an \'etale $N$-torsion group over $S$. The functor mapping $T \in \Sch_S$ to the collection of bilinear alternating surjective pairings $G_T \times G_T \rar (\mu_N)_T$ is representable by an \'etale $(\Z/N\Z)^{\times}$-torsor over $S$. In particular, this torsor is finite locally free and \'etale over $S$.}

\demo{The last sentence follows from the previous claims by descent. Let us prove the rest of the claim. First, note that the functor carries an action of $(\Z/N\Z)^{\times}$. 

Second, when restricted to the subcategory of $S_1$-schemes for some suitable \'etale faithfully flat $S$-scheme $S_1$, $G$ becomes isomorphic to the constant group scheme $(\Z/N\Z)^{\oplus 2}$, and it is then direct to check that the restricted functor is represented by the \'etale $(\Z/N\Z)^{\times}$-torsor $(\mu_N^{\times})_{S_1}$ over $S_1$.    
So, by Proposition \ref{desc-affine-schemes}, it is enough to show that the functor is an fpqc sheaf. Indeed, it is not hard to see that this functor is a subsheaf of the fpqc sheaf $\underline{\mrm{Hom}}(G \times G, \mu_N)$.}

\defi[weil-pairings-group]{Let $G$ be an \'etale $N$-torsion group over a $\Z[1/N]$-scheme. The torsor constructed in Lemma \ref{pol-torsor-exists} is denoted $\mrm{Pol}(G)$. A section of $\mrm{Pol}(G)$ is a \emph{Weil pairing} on $G$. If $G$ admits a Weil pairing, we say it is \emph{polarized}.}

\rem{Let $E$ be an elliptic curve over a $\Z[1/N]$-scheme $S$. Then $E[N]$ is an \'etale $N$-torsion group over $S$, and the usual Weil pairing is a canonical section of $\mrm{Pol}(E[N])$.} 

\defi{Let $G$ be an \'etale $N$-torsion group over a $\Z[1/N]$-algebra $R$. There is a morphism of moduli problems $\det: \P_G(N) \rightarrow [\mathrm{Pol}(G)]$ defined as follows. For any object $E/S$ of $\Ell_R$ and any $\varphi \in \P_G(N)(E/S)$, we let $\det{\varphi}$ be the bilinear alternating surjective pairing 
\vspace{-7pt}\[\det{\varphi}: G_S \times G_S \overset{\varphi\times\varphi}{\rar} E[N] \times E[N] \overset{\mrm{We}}{\rar} (\mu_N)_S.\vspace{-7pt}\]}

\defi{Let $G$ be a polarized \'etale $N$-torsion group over a $\Z[1/N]$-algebra $R$. For $u \in \mrm{Pol}(G)(R)$, let $\P_G^u$ denote the subfunctor of $\P_G(N)$ parametrizing isomorphisms $\iota$ with $\det{\iota}=u$. For $m \geq 1$ coprime to $N$, we define $\P_G^u(\Gamma_0(m)) = \P_G^u \times [\Gamma_0(m)]_R$. }

\rem{If $G$ is a polarized \'etale $N$-torsion group over a connected $\Z[1/N]$-algebra $R$, $\mrm{Pol}(G)$ is indexed by its $R$-points, so $\P_G(N)$ is the coproduct of the $\P_G^u$ over $u \in \mrm{Pol}(G)(R)$. A fix that works for general $\Z[1/N]$-algebras $R$ is to fix a Weil pairing $W$ on $G$, and then $\P_G(N)$ is the coproduct of the $\P_G^{d \cdot u}$ over $d \in (\Z/N\Z)^{\times}$.  }

\rem{When $N \geq 3$, $\P_G(N) \rightarrow [\mrm{Pol}(G)]$ is not a morphism of moduli problems of level $N$, and $\P_G^u$ does not have an obvious level $N$ structure. Indeed, we will see that the degeneracy map $D_{d,t}$ sends the moduli scheme of $\P_G^u(\Gamma_0(m))$ to itself to that of $\P_G^{d \cdot u}(\Gamma_0(t))$.}

\prop[PG-elementary]{Let $N \geq 1$ be an integer and $G$ be an \'etale $N$-torsion group over a $\Z[1/N]$-algebra $R$. Then the moduli problem $\P_G(N)$ is finite \'etale of constant rank $r_N := |\GL{\Z/N\Z}|$ over $\Ell_R$. If $N \geq 3$, it satisfies the following properties: 
\begin{enumerate}[noitemsep,label=(\roman*)]
\item\label{PG-elem-1} It is representable by an elliptic curve $\mathcal{E}_G(N)/Y_G(N)$, with $Y_G(N) \rar \Sp{R}$ smooth affine of relative dimension one,
\item\label{PG-elem-2} Its $j$-invariant is finite locally free of constant rank $\frac{r_N}{2}$, \'etale on the locus where $j(j-1728)$ is invertible,
\item\label{PG-elem-3} If $R$ is regular excellent, then $\P_G(N)$ is smooth at infinity over $\Z[1/N]$: the formation of its compactified moduli scheme $X_G(N)$ and cuspidal subscheme commutes with any base change,
\item\label{PG-elem-4} If $R$ is regular excellent, the compactified moduli scheme $X_G(N)$ attached to $\P_G(N)$ is smooth proper over $\Sp{R}$ of relative dimension one, and its compactified $j$-invariant is finite locally free of rank $\frac{r_N}{2}$,
\item\label{PG-elem-6} If $N \geq 3$, the level $N$ structure on $\P_G(N)$ defines a natural action of $(\Z/N\Z)^{\times}$ on $Y_G(N)$ which preserves the $j$-invariant. When $R$ is regular excellent Noetherian, it extends to an action of $(\Z/N\Z)^{\times}$ on $X_G(N)$ whose formation commutes with base change.
\end{enumerate} }

\demo{We already know that $\P_G(N)$ is relatively representable, affine of finite presentation, and representable if $N \geq 3$. Suppose that we can show that it is finite \'etale of constant rank $|\GL{\Z/N\Z}|$ over $\Ell_R$ and that it is smooth at infinity. Then, apart from the rank calculation for the $j$-invariant, \ref{PG-elem-1} and \ref{PG-elem-2} follow from Corollary \ref{over-Ell-concrete}, they imply \ref{PG-elem-3} and \ref{PG-elem-4} by Proposition \ref{smooth-at-infinity-base-change} and Corollary \ref{flatness-at-cusps-smooth} (apart from the \'etaleness and rank calculation for the $j$-invariant), and \ref{PG-elem-6} is a consequence of Proposition \ref{normal-near-inf-level}.

In summary, we need to show that $\P_G(N)$ is finite \'etale of rank $r_N$ and that, if $N \geq 3$, it is smooth at infinity over $\Z[1/N]$, and that its $j$-invariant has rank $\frac{r_N}{2}$. When $G$ is constant, we may assume by Proposition \ref{smooth-at-infinity-base-change} that $R=\Z[1/N]$, so that $\P_G(N) \simeq [\Gamma(N)]_{\Z[1/N]}$. Apart from the rank calculation, the claims follow from \cite[Theorem 5.1.1, Corollary 8.4.5, Proposition 8.6.8]{KM}. To compute the rank of $\P_G(N)$, note that if $E$ is an elliptic curve over an algebraically closed field $k$ of characteristic not dividing $N$, $\P_G(N)(E/k)$ has $r_N$ elements, and, if $N \geq 3$ and $j(E) \notin \{0,1728\}$, $(P,Q),(P',Q') \in \P_G(N)(E/k)$ produce equal $k$-points in the moduli scheme if and only if $(P,Q) \in \{\pm (P',Q')\}$. 

In general, $G$ is \'etale-locally constant, so by descent $\P_G(N)$ is finite \'etale of rank $r_N$ with (if $N \geq 3$) $j$-invariant of rank $\frac{r_N}{2}$. Moreover, $\P_G(N)$ is \'etale-locally over $R$ smooth at infinity over $\Z[1/N]$, so and the result follows by descent and Proposition \ref{smooth-at-infinity-descends}.
}

\nott{We denote by $\mathcal{E}(N)/Y(N)$ the elliptic curve $\mathcal{E}_{(\Z/N\Z)^{\oplus 2}}(N)/Y_{(\Z/N\Z)^{\oplus 2}}(N)$ (over the base ring $\Z[1/N]$). The compactification of $Y(N)$ is called $X(N)$. Since $\GL{\Z/N\Z}$ acts on the left on $(\Z/N\Z)^{\oplus 2}$, it acts on the right on $\P := \P_{(\Z/N\Z)^{\oplus 2}}(N)$; this moduli problem with its right action of $\GL{\Z/N\Z}$ is exactly the base change to $\Z[1/N]$ of the moduli problem $[\Gamma(N)]$ defined in \cite[(3.1)]{KM}.     
This procedure identifies the action of $a \in (\Z/N\Z)^{\times}$ (coming from the level structure) with the action of the matrix $aI_2 \in \GL{\Z/N\Z}$. 
}

\nott{We recall a notation introduced at the beginning: for an integer $N \geq 1$, $\OO_N$ is the ring $\Z[1/N,\zeta_N]$, where $\zeta_N$ is a primitive $N$-th root of unity. For $a \in (\Z/N\Z)^{\times}$, $\underline{a}$ is the automorphism of $\OO_N$ (or $\Sp{\OO_N}$) given by $\zeta_N \mapsto \zeta_N^a$. }

\prop[XN-Weil]{Let $\phi^{\mrm{univ}}: (\Z/N\Z)^{\oplus 2}_{Y(N)} \rightarrow \cE(N)[N]$ be the universal isomorphism of finite \'etale $Y(N)$-group schemes. The Weil pairing $\langle \phi^{\mrm{univ}}(1,0) ,\,\phi^{\mrm{univ}}(0,1)\rangle_{\mathcal{E}(N)}$ defines a map $\mrm{We}: Y(N) \rar \Sp{\OO_N}$, which is smooth of relative dimension one with geometrically connected fibres. It extends to a smooth proper map $\mrm{We}: X(N) \rar \Sp{\OO_N}$ with geometrically connected fibres of dimension one. 

Furthermore, the natural right action of $\GL{\Z/N\Z}$ on $Y(N)$ preserves the $j$-invariant and extends uniquely to $X(N)$. Moreover, the following diagram commutes, for any $M \in \GL{\Z/N\Z}$:
 \[
\begin{tikzcd}[ampersand replacement=\&]
X(N)  \arrow{r}{M}\arrow{d}{\mrm{We}}\& X(N)\arrow{d}{\mrm{We}}\\
\Sp{\OO_N} \arrow{r}{\underline{\det{M}}} \& \Sp{\OO_N}  
\end{tikzcd}
\]
}

\demo{Since $\OO_N$ is finite \'etale over $\Z[1/N]$, the extension of $\mrm{We}$ follows from Corollary \ref{extension-to-finite}. Since $X(N)$ is smooth proper of relative dimension one over $\Z[1/N]$ and $\OO_N$ is finite \'etale over $\Z[1/N]$, $\mrm{We}: X(N) \rar \Sp{\OO_N}$ is smooth proper of relative dimension one. It is clear that $\GL{\Z/N\Z}$ acts on $Y(N)$ by automorphisms of $\bA^1_{\Z[1/N]}$-schemes, so this action extends to an action on $X(N)$ by automorphisms of $\bP^1_{\Z[1/N]}$-schemes by Proposition \ref{compactification-functor}. Checking that the diagram commutes can also be done after restricting to $Y(N)$ (by Proposition \ref{compactification-functor}), where it is an easy computation. 

Note that $\mrm{We}: Y(N) \rar \Sp{\OO_N}$ is the moduli scheme of the moduli problem $[\Gamma(N)]^{\Z[\zeta_N]-\text{can}}_{\OO_N}$ (see \cite[Propositions 9.3.1, (9.4.3.1)]{KM}). Hence, by \cite[Corollary 10.9.2]{KM}, the geometric fibres of $\mrm{We}: X(N) \rar \Sp{\OO_N}$ are irreducible. In particular, $\mrm{We}: Y(N) \rar \Sp{\OO_N}$ also has geometrically connected fibres. }

\prop[XG-to-polG]{Suppose that $G$ is an \'etale $N$-torsion group over some $\Z[1/N]$-algebra $R$ with $N \geq 3$. Then the rule $(E/S,\iota) \longmapsto \det{\iota} \in \mrm{Pol}(G)(S)$ defines an affine morphism $\det: Y_G(N) \rar \mrm{Pol}(G)$, which extends to a proper morphism $\det: X_G(N) \rar \mrm{Pol}(G)$ if $R$ is regular and excellent. In either case, $\det$ is smooth of relative dimension one with geometrically integral fibres and one has $\det \circ [u] = u^2 \cdot \det$ for $u \in (\Z/N\Z)^{\times}$. }

\demo{It is clear that $\det$ is well-defined and affine, since it is between affine schemes. By Corollary \ref{extension-to-finite}, it extends to a proper map $X_G(N) \rightarrow \mrm{Pol}(G)$ if $R$ is regular excellent. The compatibility with the action of $(\Z/N\Z)^{\times}$ can be checked on $Y_G(N)$, where it is formal. 

To check that that $\det$ is smooth and the geometric fibres of $\det$ are connected, we can work \'etale-locally over $R$ and assume that $G$ is constant. It is then not hard to see that $\det$ identifies with the base change by $\Z[1/N] \rightarrow R$ of the (affine) morphism $\mathrm{We}$ of Proposition \ref{XN-Weil}, so $\det$ is indeed smooth with geometrically connected fibres of dimension one. The same argument applies for $X_G(N)$ and the extension of $\det$ if $R$ is regular and excellent.}

\cor[counting-connected-components-XG]{In the situation of Proposition \ref{XG-to-polG}, the geometric fibres of $Y_G(N) \rar \Sp{R}$ (and $X_G(N) \rar \Sp{R}$ if $R$ is regular and excellent) have $\varphi(N)$ connected components.}

\demo{We may assume that $R$ is an algebraically closed field $k$. Then $\mrm{Pol}(G)$ consists of $\varphi(N)$ copies of $\Sp{k}$ and we conclude by Proposition \ref{XG-to-polG} since the geometric fibres of $\det$ are geometrically connected.}

\cor[enriched-jinv-free]{In the situation of Proposition \ref{XG-to-polG}, the map $(j,\det): Y_G(N) \rar \bA^1_R \times \mrm{Pol}(G)$ (or $(\overline{j},\det): X_G(N) \rar \bP^1_R \times \mrm{Pol}(G)$ if $R$ is regular excellent) is finite locally free of rank $\frac{|\SL{\Z/N\Z}|}{2}$. }

\demo{Since $\mrm{Pol}(G)$ is finite \'etale over $R$, this map is finite locally free because the $j$-invariant is. To compute the rank, it is enough to deal with the uncompactified case because $\bA^1_{\mrm{Pol}(G)} \subset \bP^1_{\mrm{Pol}(G)}$ is open and dense. Working \'etale-locally over $G$, we may assume that $G$ is constant; by base change we are reduced to the case $R=\Z[1/N]$. Then $\bA^1_R \times \mrm{Pol}(G)$ is connected, so $(j,\det)$ has constant rank $\rho$. Because $\mrm{Pol}(G)$ is finite free over $\Sp{R}$ of rank $\varphi(N)$ ($\varphi$ being Euler's totient function) one has $\varphi(N)\rho=\frac{|\GL{\Z/N\Z}|}{2}$.}

\medskip

\prop[forgetful-maps]{Let $N,d \geq 3$ be integers with $d \mid N$. Let $G$ be an \'etale $N$-torsion group over a $\Z[1/N]$-algebra $R$. Then $G[d]$ is an \'etale $d$-torsion group over $R$, and there is an obvious forgetful map $\mrm{Pol}(G) \rightarrow \mrm{Pol}(G[d])$ of $R$-schemes which respects the action of $(\Z/N\Z)^{\times}$. The forgetful map $Y_G(N) \rightarrow Y_{G[d]}(d)$ of $\bA^1_R \times \mrm{Pol}(G[d])$-schemes is finite \'etale. If furthermore $R$ is regular and excellent, it extends to a finite locally free morphism $X_G(N) \rightarrow X_{G[d]}(d)$ of $\bP^1_R \times \mrm{Pol}(G[d])$-schemes.}

\demo{The group scheme $G[d]$ is clearly an \'etale $d$-torsion group, and the morphism 
\vspace{-7pt}\[\varphi: G \times_R G \rightarrow G[d] \times_R G[d], \qquad (x,y) \mapsto \left(\frac{N}{d}x,\frac{N}{d}y\right)\vspace{-7pt}\] is finite \'etale surjective. Let $S$ be a $R$-scheme and $\pi \in\mrm{Pol}(G)(S)$. Then, if $T$ is a $S$-scheme and $P,Q \in G_S(T)$, the quantity $\frac{N}{d}\pi(P,Q)$ only depends on the couple $\varphi(P,Q)$. Since $\varphi$ is finite \'etale surjective, by descent, there is a unique $\pi[d]: G[d]_S \times_S G[d]_S \rightarrow (\mu_d)_S$ such that $\pi[d]\circ \varphi = \frac{N}{d}\pi$. By construction (and because $\varphi$ is faithfully flat), $\pi[d]$ is bilinear and alternating. Since $\pi$ and $\frac{N}{d}: G \rightarrow G[d]$ are surjective, $\pi$ is surjective. Thus $\pi \mapsto \pi[d]$ defines a morphism $\mrm{Pol}(G) \rightarrow \mrm{Pol}(G[d])$. 

For any object $E/S$ of $\Ell_R$, the rule 
\vspace{-7pt}\[(\iota: G_S \rightarrow E[N]) \in \P_G(N)(E/S) \mapsto (\iota[d]: G[d]_S \rightarrow E[d]) \in \P_{G[d]}(d)(E/S)\vspace{-7pt}\]
 defines a morphism of moduli problems $f: \P_G(N) \Rightarrow \P_{G[d]}(d)$ and induces a map $f'$ between the moduli schemes. By construction, $f'$ is a morphism of $\bA^1_R$-schemes, and the properties of the Weil pairing show that it also preserves the maps to $\mrm{Pol}(G[d])$. As a morphism of finite $\bA^1_R$-schemes, $f'$ is finite, so let us check that it is \'etale. This is Zariski-local over $R$, so we may assume that some odd prime $\ell \nmid N$ is invertible in $R$. For any object $E/S$ of $\Ell_R$, $f_{E/S}: \P_G(N)_{E/S} \rightarrow \P_{G[d]}(d)_{E/S}$ is finite \'etale as a morphism of finite \'etale $S$-schemes. In particular, $f: \MPp{[\Gamma(\ell)]_R,\P_G(N)} \rightarrow \MPp{[\Gamma(\ell)]_R,\P_{G[d]}(d)}$ is finite \'etale, hence $\MPp{[\Gamma(\ell)]_R,\P_G(N)} \rightarrow Y_G(d)$ is finite \'etale. Since $[\Gamma(\ell)]_R$ is \'etale over $\Ell_R$, $f'$ is \'etale by \cite[$\text{IV}_4$ (17.7.7)]{EGA}. The last claim then follows by Proposition \ref{PG-elementary} and Corollary \ref{flatness-at-cusps-smooth}. }
 
\medskip

\prop[XG-pol]{Let $N \geq 3$ and let $G$ be a polarized \'etale $N$-torsion group over a connected $\Z[1/N]$-algebra $R$. Let, for each $u \in \mrm{Pol}(G)(R)$, $Y_G^u(N)$ (resp. $X_G^u(N)$ if $R$ is regular excellent) be the inverse image of $u$ under $\det: Y_G(N) \rar \mrm{Pol}(G)$ (resp. $\det: X_G(N) \rar \mrm{Pol}(G)$).

The moduli problem $\P_G^u$ over $R$ is represented by the elliptic curve 
\vspace{-7pt}\[\cE_G(N) \times_{Y_G(N)} Y_G^u(N)/Y_G^u(N),\vspace{-7pt}\] and it is finite \'etale of rank $|\SL{\Z/N\Z}|$ over $\Ell_R$.  

If $R$ is regular excellent, $\P_G^u$ is smooth at infinity and its compactified moduli scheme identifies with $X_G^u(N)$ as a $\CMPp{\P_G(N)}$-scheme. }

\demo{If $(\mathcal{E}_G(N)/Y_G(N),\alpha)$ represents $\P_G(N)$, the couple $(\mathscr{E}_G^u(N)/Y_G^u(N), \alpha_{|\det^{-1}(u)})$ tautologically represents $\P_G^u$.

For any object $E/S$ of $\Ell_R$, $(\P_G^u)_{E/S} \rightarrow \P_{G}(N)_{E/S}$ is a base change of the closed open immersion $u: \Sp{R} \rar \mrm{Pol}(G)$, so $(\P_G^u)_{E/S}$ is a finite \'etale $S$-scheme. Thus $\P_G^u$ is finite \'etale over $\Ell_R$. If $k$ is an algebraically closed field and a $R$-algebra and $E/k$ is an elliptic curve, one checks easily that $\P_G^u(E/k)$ has $|\SL{\Z/N\Z}|$ elements. Hence $\P_G^u$ is finite \'etale over $\Ell_R$ with rank $|\SL{\Z/N\Z}|$.  

Suppose now that $R$ is regular and excellent. Then $Y_G(N)$ is the disjoint union of the $\MPp{\P_G^u}$, so all $\P_G^u$ are normal near infinity and by \cite[Lemma 03GO]{Stacks} $X_G(N)$ is the disjoint union of the $\CMPp{\P_G^u}$. The composition $\CMPp{\P_G^u} \overset{\det}{\rightarrow} \CMPp{\P_G(N)} \rightarrow \mrm{Pol}(G)$ is equal to $u$ on the scheme-theoretically dense open subscheme $\MPp{\P_G^u}$, so it is equal to $u$ on $\CMPp{\P_G^u}$. Since $X_G(N)$ is the disjoint union of the $\CMPp{\P_G^u}$, $\CMPp{\P_G^u}$ agree with $X_G^u(N)$ as a $X_G(N)$-scheme.} 

\rem{Let $G$ be a polarized \'etale $N$-torsion group scheme over a connected $\Z[1/N]$-algebra $R$ (with $N \geq 3$). In $Y_G(N)$, the equality $\det \circ [v]=v^2 \cdot \det$ holds for any $v \in (\Z/N\Z)^{\times}$, so $[v]$ induces an isomorphism $Y_G^u(N) \rar Y_G^{v^2 \cdot u}(N)$ of $\bA^1_R$-schemes. When $R$ is regular and excellent, it extends by Proposition \ref{compactification-functor} to an isomorphism $X_G^u(N) \rar X_G^{v^2 \cdot u}(N)$ of $\bP^1_R$-schemes. }

\medskip

\cor[XG-pol-prop]{Let $G$ be a polarized \'etale $N$-torsion group over a connected $\Z[1/N]$-algebra $R$ with $N \geq 3$. Then, for any $u \in \mrm{Pol}(G)(R)$:
\begin{enumerate}[noitemsep,label=$(\roman*)$]
\item\label{PGu-1} $Y_G^u(N) \rar \Sp{R}$ (resp. $X_G^u(N) \rar \Sp{R}$ if $R$ is regular excellent) is affine (resp. proper) and smooth of relative dimension one with connected geometric fibres. 
\item\label{PGu-2} The $j$-invariant $Y_G^u(N) \rar \bA^1_R$ (resp. $X_G^u(N) \rar \bP^1_R$ if $R$ is regular excellent) is finite locally free of rank $\frac{|\SL{\Z/N\Z}|}{2}$. 
\item\label{PGu-3} If $R$ is regular excellent, then $\P_G^u$ is smooth at infinity over $\Z[1/N]$ and the formation of its compactified moduli scheme and cuspidal subscheme commute with base change. 
\end{enumerate}
}

\medskip

\nott{Let $N \geq 3$ be an integer and $E$ be an elliptic curve over a $\Z[1/N]$-algebra $R$. We write $Y_E(N) := Y_{E[N]}(N)$. If $\alpha \in (\Z/N\Z)^{\times}$, we write $Y_E^{\alpha}(N)$ for the closed subscheme $\det^{-1}(\alpha \cdot \mrm{We})$ of $Y_E(N)$, where $\mrm{We} \in \mrm{Pol}(E[N])(R)$ is the (usual) Weil pairing. 

If furthermore $R$ is a regular excellent domain, we write 
\vspace{-7pt}\[X_E(N) := X_{E[N]}(N),\qquad X_E^{\alpha}(N) := X_{E[N]}^{\alpha \cdot \mrm{We}}(N).\vspace{-7pt}\]}

\subsection{Some general remarks about the bad fibres}

One can summarize the results of Section \ref{subsect:etale-case} as follows. 

\prop[potentially-good-red]{Let $K$ be a global field or local non-Archimedean field and $N \geq 6$ be an integer invertible in $K$. Let $E/K$ be an elliptic curve, $a \in (\Z/N\Z)^{\times}$, and $v \nmid N$ be a (finite) place of $K$. The smooth projective geometrically connected curve $X_E^a(N)$ over $K$ (with genus at least $2$) has potentially good reduction at $v$. If $E$ has furthermore good reduction at $v$, so does $X_E^a(N)$, and the special fibre of $X_E^a(N)$ at $v$ is $X_{\overline{E}}^a(N)$, where $\overline{E}$ is the reduction of $E$ modulo $v$. }

\demo{The ring $\OO_v$ of $v$-integers in $K$ is a regular excellent domain (by \cite[$\text{IV}_2$ (7.8.3)]{EGA}). Let $k$ be its residue field. If the \'etale $N$-torsion group $E[N]$ over $K$ extends to a finite \'etale $\OO_v$-group scheme $G$, then $G$ is an \'etale $N$-torsion group and the Weil pairing of $E$ extends to a Weil pairing $W$ on $G$ (cf. Definition \ref{weil-pairings-group}). Then one has $X_E^a(N) = X_G^{a \cdot W}(N)_K$, where $X_G^{a \cdot W}(N)$ is a smooth proper $\OO_v$-scheme of relative dimension one with special fibre $X_{G_k}^{a \cdot W}(N)$. 

If $E$ is the generic fibre of an elliptic curve $\cE/\OO_v$, we can simply take $G=\cE[N]$, and $X_{G_k}^{a \cdot W}(N)=X_{\cE_k}^a(N)$. In general, the $N$-torsion of $E$ is defined over a finite separable extension $L/K(\zeta_N)$. So $E[N]_L \simeq (\Z/N\Z)^{\oplus 2}_L$ clearly extends to the valuation ring of any place $w$ above $v$, and the conclusion follows.}

There are three natural (and not mutually exclusive) directions in which one can go in order to extend Proposition \ref{potentially-good-red}:
\begin{enumerate}[label=(\arabic*),itemsep=1pt]
\item\label{add-hecke} We can study the geometry of the curves $X_E^a(N,\Gamma_0(m))$, where $m$ is square-free and coprime to $N$, at places dividing $m$ of good reduction for $E$.    
\item\label{bad-reduction} We can study the geometry of the curves $X_E^a(N)$ at places $v \nmid N$ of bad reduction for $E$.  
\item\label{bad-places} We can study the geometry of the curves $X_E^a(N)$ at places $v \mid N$ of good reduction for $E$. 
\end{enumerate} 

The direction \ref{add-hecke} will be pursued in Section \ref{subsect:etale+gamma0} in order to show that the $X_G^a(N)$ admit well-behaved Hecke operators. Because its claims fit naturally in the set-up of \cite[Chapters 12--13]{KM}, it is the more approachable of the three. 

The directions \ref{bad-reduction} and \ref{bad-places} seem more complex in comparison. The immediate issue is that our moduli problems do not have good enough formal properties any more: while the definition of $\P_G(N)$ carries over in the case \ref{bad-places}, the moduli problem is not finite any more, which makes it delicate to discuss its compactification. In the case \ref{bad-reduction}, the moduli problem $\P_G(N)$ simply does not make sense any more. A satisfactory resolution of these problems would likely involve a better definition of $\P_G(N)$, an achievement similar to the moduli-theoretic description of $X_0(N), X_1(N), X(N)$ over $\Z$ given in \cite{KM}. 

Our aim in this section is less ambitious and we will not attempt to formulate a more subtle definition of $\P_G(N)$ that adapts to the case \ref{bad-reduction}; we will merely point out some more accessible facts about the geometry of $X_E^a(N)$ in this case. The case of \ref{bad-places} will be studied starting from Section \ref{subsect:smoothness-general}.

\rem{Let $E/\Q$ be an elliptic curve and $N \geq 3$ be an integer. It is not true that the smooth projective geometrically connected curve $X_E^1(N)$ over $\Q$ has good reduction at places not dividing $N$. 

To give an explicit example, we can for instance consider the case of the curve $X_E^1(7)$, whose explicit equation is determined by \cite[Theorem 2.2]{HalbKraus}. Following \cite[\S 6]{HalbKraus}, take for $E/\Q$ the elliptic curve with minimal Weierstrass equation 
\vspace{-7pt}\[y^2 + xy + y = x^3 - 8x - 7.\vspace{-7pt}\]
It has conductor $105$, rational $2$-torsion, and LMFDB label \cite{lmfdb} $105.a2$ (in Cremona's tables \cite{Cremona-tables}, $E/\Q$ is the curve $105A2$). 

By the LMFDB, $E$ has the reduced Weierstrass equation $y^2=(x+33)(x+78)(x-111)$, so that, by \cite[Theorem 2.2]{HalbKraus}, $X_E^1(7)$ admits the following projective plane quartic 
\vspace{-7pt}\[144\left[(y+z)x^3-3x^2yz\right]-189\left[(z+x)y^3-3y^2xz\right]+45\left[(x+y)z^3-3z^2xy\right]=0\vspace{-7pt}\] as an equation.

A direct computation with Magma \cite{magma} shows that this smooth projective plane quartic is isomorphic over $\Q$ to the one defined by the equation 
\vspace{-7pt}\[x^3y - 3x^2y^2 - 3x^2yz + 3x^2z^2 + 4xy^3 + 6xy^2z - 9xyz^2 + 7xz^3 - 4y^4 + 9y^3z - 3y^2z^2 - 2yz^3 - 6z^4 = 0, \vspace{-7pt}\] which is minimal. 

Another computation with Magma shows that the curve defined by this ternary quartic form is smooth over $\Q$ and every $\F_p$ with $p \nmid 105$, and that it is not smooth at $p=3,5,7$. }

\rem{In the setting of Proposition \ref{potentially-good-red} with $N > 5$ prime, the general calculation of the Tate module of the Jacobian of $X_{E}^a(N)$ done in the author's PhD dissertation \cite{Studnia-thesis} shows that, when this Jacobian has good reduction at some place $v$ of $K$ (not dividing $N$), then the action of the inertia subgroup at $v$ of $\mrm{Gal}(K_s/K)$ on $E[N]$ is given by $\{\pm \mrm{id}\}$. Thus, after replacing $E$ if neeeded with a quadratic twist (which does not change $X_{E}^a(N)$, by Proposition \ref{quadratic-twist-equiv} and the main results of Section \ref{moduli-is-twisting}), $E[N]$ extends to a finite \'etale group scheme over the ring of $v$-integers of $K$ by the results of \cite[Exp. V, Th. 4.1, \S\S 7--8]{SGA1}.}

\subsection{Equations of $X_E^{\alpha}(N)$ over a regular excellent base}
\label{subsect:upgrade-eqns}

In this section, we analyze to which extent the equations found for various $X_E^{\alpha}(N)$ over a field (if needed of characteristic zero) extend to any basis.

\prop[generic-iso-of-smooth-extends]{Let $R$ be a normal Noetherian domain with generic point $\eta \in \Sp{R}$. Let $X,Y,Z$ be three smooth $R$-schemes. Let $f: X\rar Z,\, g: Y \rar Z$ be finite dominant morphisms of $R$-schemes and $h_{\eta}: X_{\eta} \rightarrow Y_{\eta}$ be an isomorphism of $Z_{\eta}$-schemes. Then $h_{\eta}$ extends to an isomorphism $X \rightarrow Y$ of $Z$-schemes.  
}

\demo{Since $X \rightarrow S$ is smooth, the relative normalization of $X_{\eta} \rightarrow X$ is $X$ by \cite[$\text{IV}_2$ (6.14.5)]{EGA}. Since $X$ is finite over $Z$, the decomposition $X_{\eta} \rightarrow X \rightarrow Z$ is the relative normalization by \cite[Lemma 035I]{Stacks}. The same applies after replacing $X$ by $Y$. Hence, by functoriality, $h_{\eta}, h_{\eta}^{-1}$ extend to morphisms $\alpha: X \rightarrow Y$, $\beta: Y \rightarrow X$ of $Z$-schemes. Then $\beta \circ \alpha: X \rightarrow X$ is a $Z$-morphism that agrees with the identity on $X_{\eta}$. Since $X_{\eta} \rightarrow X$ is scheme-theoretically dominant by \cite[$\text{IV}_3$ (11.10.5)]{EGA} and $X$ is separated over $Z$, one has $\beta\circ\alpha=\mrm{id}_X$. Similarly, one has $\alpha \circ \beta=\mrm{id}_Y$, so $\alpha$ is an isomorphism. }

\lem[reduce-to-generic-fibre]{Let $S$ be a Noetherian integral scheme and $X$ be a dominant $S$-scheme of finite type whose nonempty fibres are pure of dimension $d$. The following are equivalent:
\begin{itemize}[noitemsep,label=$-$]
\item $X$ has dense generic fibre,
\item the set-theoretic closure of the generic fibre of $X$ meets the interior of every irreducible component of every non-empty fibre of $X \rar S$.
\end{itemize}}

\demo{The forward direction is clear, let us show the reverse one. Let $\eta \in S$ be the generic point and $j: H \rightarrow X$ be the scheme-theoretic image of $X_{\eta} \rightarrow X$. By \cite[I (9.5.4)]{EGA}, $j(H)$ meets the interior of every irreducible component of every fibre of $X \rightarrow S$ and we need to show that $j$ is surjective.  

The irreducible components of $j(H)_{\eta}=X_{\eta}$ have dimension $d$. Let $s \in S$. Since $H$ has dense generic fibre, the irreducible components of $j(H)_s$ have dimension at least $d$ by \cite[$\text{IV}_3$ (13.1.3)]{EGA}. Since $X_s$ is pure of dimension $d$, $j(H)_s$ is a reunion of irreducible components of $X_s$. Now, $j(H)_s$ meets the interior of any irreducible component of $X_s$, so $j(H)_s=X_s$ and $j$ is surjective.}

\prop[miracle-flatness-bis]{Let $R$ be a normal Noetherian domain and $X$ be a $R$-scheme of finite type. Assume that
\begin{itemize}[noitemsep,label=$-$]
\item the non-empty fibres $X \rar \Sp{R}$ are geometrically reduced and pure of dimension $d \geq 1$,
\item the set-theoretic closure of the generic fibre of $X$ meets the interior of every irreducible component of every non-empty fibre of $X \rar \Sp{R}$,
\item the generic fibre of $X$ is irreducible,
\end{itemize}
Then $X$ is flat over $R$.}

\demo{I am grateful to O. Gabber for telling me about the following proof.  

By Lemma \ref{reduce-to-generic-fibre}, $X$ has dense generic fibre, so $X$ is irreducible. Then $X$ is equidimensional over $\Sp{R}$ in the sense of \cite[$\text{IV}_3$ (13.2.2)]{EGA}. Because $R$ is a normal Noetherian domain, $\Sp{R}$ is geometrically unibranch (cf. \cite[(\textbf{0}, 23.2.1); $\text{IV}_2$ (6.15.1)]{EGA}). Thus $X \rar \Sp{R}$ is universally open by \cite[$\text{IV}_3$ (14.4.4)]{EGA}, hence flat by \cite[$\text{IV}_3$ (15.2.3)]{EGA}. }

\prop[models-of-XEalpha]{Let $N \geq 3$ be an integer and $R$ be a regular excellent domain where $N$ is invertible. Let $E/R$ be an elliptic curve and fix $\alpha \in (\Z/N\Z)^{\times}$. Let $X' := X_E^{\alpha}(N)$ (which is a finite flat $\bP^1_R$-scheme by the $j$-invariant). Suppose given a $\bP^1_R$-scheme $S$, projective over $R$, such that, for every morphism $R \rar k$ with $k$ a field, there is an isomorphism $\iota_k: S_k \rightarrow X'_k$ of $\bP^1_k$-schemes. 

Then $\iota_{\eta}: S_{\eta} \rar X'_{\eta}$ extends to an isomorphism $\iota: S \rar X'$ of $\bP^1_R$-schemes. }

\demo{The fibres of $S \rar \Sp{R}$ are isomorphic to those of $X'$, so they are smooth connected of dimension one by Corollary \ref{XG-pol-prop}, hence they are geometrically integral of dimension one. Let $Z \subset S$ be the set-theoretic closure of its generic fibre. Since $S$ is projective, the image of $Z$ in $\Sp{R}$ is closed and contains the generic point, hence $Z$ meets every fibre of $S \rar \Sp{R}$. By Proposition \ref{miracle-flatness-bis}, $S$ is flat over $R$, thus smooth over $R$. The assumptions imply that the structure map $j': S \rightarrow \bP^1_R$ is quasi-finite (because this can be checked after a base change $R \rightarrow k$); both schemes are proper over $R$ so $j'$ is finite. We can then use Proposition \ref{generic-iso-of-smooth-extends}. }

\cor{Let $R$ be a regular excellent domain which is a $\Q$-algebra. Let $n \in \{7,8,9,11,13\}$ and $\alpha \in (\Z/n\Z)^{\times}$. Let $E/R$ be an elliptic curve with Weierstrass equation $y^2=x^3+ax+b$ and assume that, if $n \geq 11$, $j(E)(j(E)-1728)$ is invertible in $R$. Then the given equations for $X_E^{\alpha}(n)$ (over $\mrm{Frac}(R)$) hold over $R$.}

\demo{\emph{Case where $n=7$:} We use the equations of \cite[Theorem 1.1]{Fisher-711}. We must show that the $j$-map defined by Fisher in Remark 4.3 of \emph{op.cit.} extends to a map $S \rar \bP^1_R$ (where $S$ denotes the projective $R$-scheme given by the equation $\mathcal{F}=0$ of Theorem 1.1 of \emph{op.cit.}). We need to show that if $R'$ is a $R$-algebra, $x_0,y_0,z_0 \in R'$ generate the unit ideal, and one has $\mathcal{F}(x_0,y_0,z_0)=0$, then $c_4(\mathcal{F})(x_0,y_0,z_0)$ and $H(\mathcal{F})(x_0,y_0,z_0)$ generate the unit ideal of $R'$ (this is because one has $\psi(\mathcal{F}) \in (R')^{\times}$, cf. Definition 4.1 of \emph{op.cit.}). We can assume that $R'$ is an algebraically closed field $k$ of characteristic zero, in which case $\mathcal{F}$ is deduced from the equation $F := x^3y+y^3z+z^3x$ by post-composing with a linear action of $\PSL{\F_7}$ on $kx \oplus ky \oplus kz$: thus, it suffices to show that if $k$ is an algebraically closed field of characteristic zero, the common vanishing locus of $F,H,c_4$ (as defined in pp. 11--12 of \emph{op.cit.}) in $\bA^3$ is reduced to zero. This can be checked explicitly in Magma, so we are done.  

\medskip

\emph{Case where $n=8$:}
We use the equations of \cite[Theorem 1.7.2, 1.7.3, 1.7.4, 1.7.5]{Chen-PhD} with the forgetful map to $X_E^{\alpha}(4)$ always given by Corollary 6.1.3 of \emph{op.cit.}, so the $j$-invariant is always given by $[1728\tilde{A_4}(t,s):\tilde{A_4}(t,s)^3-\tilde{B_4}(t,s)^2]$, where $\tilde{A_4}(t,s)=s^8A_4(t/s), \tilde{B_4}(t,s)=s^{12}B_4(t/s)$, and $A_4,B_4$ are defined in Theorem 3.3.1 of \emph{op.cit.}. To check that this map is well-defined over $R$, we need to check 
\begin{itemize}[noitemsep,label=\tiny$\bullet$]
\item that over the base ring $R_0 := \Q[c_4,c_6,(c_4^3-c_6^2)^{-1}]$, both $t$ and $s$ are in the radical of the ideal $(\tilde{A_4},\tilde{B_4})$,
\item that over the base ring $R_0 := \Q[a,b,(4a^3+27b^2)^{-1}]$, for $i \in \{1,5,7\}$, the radical of the homogeneous ideal of $R_0[x_0,x_1,x_2,x_3,x_4]$ generated by $x_0,3x_4,f_i,g_i,h_i$ is $(x_0,x_1,x_2,x_3,x_4)$, where $f_1,g_1,h_1$ (resp. $f_5,g_5,h_5$, resp. $f_7,g_7,h_7$) are defined in Theorem 1.7.2 (resp. 1.7.4, resp. 1.7.5) of \emph{op.cit.}, and $3x_4$ comes from the coordinate substitution between Corollary 6.1.3 and the end of the proof of Theorem 1.7.2 (resp. from the description of the forgetful map to $X_E^{1}(4)$ in Theorem 6.3.3, resp. the coordinate substitution at the end of the proof of Theorem 1.7.5),
\item that over the same $R_0$, the radical of the homogeneous ideal of $R_0[t,u_0,u_1,u_2,s]$ generated by $t,s,F_3,G_3,H_3$ (as in the proof of Theorem 1.7.2 of \emph{op.cit.} which specifies the forgetful map to $X_E^3(4)$) is $(t,u_0,u_1,u_2,s)$.
\end{itemize}

All of this can easily be checked in Magma. 

\medskip

\emph{Case where $n=9$:}
We use the equations of \cite{Fisher-9} over a regular excellent $\Q[a,b,(4a^3+27b^2)^{-1}]$-algebra $R$ which is a domain and verify that the $j$-invariant is defined over $R$ in a ``two-step'' strategy similar to the case $n=8$. Indeed, we can first directly check in Magma that, if $A^{\pm},B^{\pm}$ are the polynomials in $\Q[c_4,c_6,(c_4^3-c_6^2)^{-1},r,s]$ described in Theorem 1.1 of \emph{op.cit.}, then for each $\epsilon \in \{\pm\}$, one has $\sqrt{(A^{\epsilon},B^{\epsilon})}\supset (r,s)$. 

Fix a sign $\epsilon\in \{\pm 1\}$. Then we need to show that the map $\varphi_{\epsilon}: X_{\mrm{Frac}(R)}^{\epsilon} \rar \mathbb{P}^1_{\mrm{Frac}(R)}$ defined in Theorem 1.3 of \emph{op.cit} comes from a morphism $X^{\epsilon} \rar \mathbb{P}^1_R$, where $X^{\epsilon}$ is the subscheme of $\mathbb{P}^3_R$ (with homogeneous coordinates $x,y,z,t$) cut out by the ideal $(F_1^{\epsilon},F_2^{\epsilon})$ and $F_i^{\epsilon}(x,y,z,t)$ denotes the polynomial in $4$ variables with coefficients in $R$ defined in Theorem 1.2 of \emph{op.cit.}). 

This does not directly follow from the main results of \emph{op.cit.}, because while the equation $[\gamma_2^{\epsilon}:3\gamma_1^{\epsilon}]$ given in Theorem 1.3 of \emph{op.cit.} can be explicitly computed, it is not globally defined. Indeed, it can be checked that $\gamma_2^{\epsilon},\gamma_1^{\epsilon}$ vanish whenever $F_1^{\epsilon},F_2^{\epsilon}$ and the $D_{\epsilon}$ at the bottom of p. 11 of \emph{op.cit.} simultaneously vanish, which does happen.

However, it is possible, using the results of \emph{op.cit.}, to work out the following formulas that Fisher was happy to send me upon request \cite{Fisher-perso}: one has $[\gamma_2^{\epsilon}:3\gamma_1^{\epsilon}]=[c_2^{\epsilon}:9c_1^{\epsilon}]$, where the $c_i^{\epsilon} \in \Z[a,b,x,y,z,t]$ are the following polynomials: 

\begin{align*}
c_1^{+1} &= &&3x^2y - 9bxz^2 + 6a^2xzt + 9abxt^2 - 9ay^3 - 54by^2z\\
         &   &&+ 18a^2y^2t + 9a^2yz^2 + 54abyzt - 9a^3yt^2 + 2a(a^3 + 9b^2)t^3,\\
c_2^{+1} &= &&-x^3 - 6ax^2z - 9axy^2 + 108bxyz - 9a^2xz^2 - 54abxzt\\
         &   &&- 3a^3xt^2 + 108by^3 - 162aby^2t + 27abyz^2 - 18a^3yzt \\ 
         &   &&+ 81a^2byt^2 - 6(a^3 + 9b^2)z^3 - 24a^3bt^3,\\
c_1^{-1} &= &&9bx^3 + 6a^2x^2z + 9abx^2t + 6(4a^3 + 27b^2)xyt - 9abxz^2 \\
         &   &&+ 6(2a^3 + 9b^2)xzt + 3a^2bxt^2 + 6(4a^3 + 27b^2)y^3 + 9(4a^3 + 27b^2)y^2z \\
         &   &&+ 6(4a^3 + 27b^2)yz^2 + 3(2a^3 + 15b^2)z^3 - 3a^2bz^2t \\
         &   &&+ a(2a^3 + 9b^2)zt^2 + 3b(a^3 + 6b^2)t^3,\\  
c_2^{-1} &= &&6a^2x^3 - 27abx^2z - 18(a^3 + 9b^2)x^2t - 27(4a^3 + 27b^2)xy^2 \\
         &   &&- 18(4a^3 + 27b^2)xyz - 9(2a^3 + 9b^2)xz^2 - 18a^2bxzt + 3a(2a^3 + 9b^2)xt^2 \\
         &   &&+ 9a(4a^3 + 27b^2)y^2t - 18b(4a^3 + 27b^2)yt^2 + 3a^2bz^3 - 3a(2a^3 + 9b^2)z^2t\\
         &   &&- 27b(a^3 + 6b^2)zt^2 - a^2(2a^3 + 15b^2)t^3.
\end{align*}

One checks furthermore in Magma that in the ring $\Q[a,b,(4a^3+27b^2)^{-1},x,y,z,t]$, one has $\sqrt{(c_1^{\epsilon},c_2^{\epsilon},F_1^{\epsilon},F_2^{\epsilon})} \supset (x,y,z,t)$. 

\medskip

\emph{Case where $n=11$:}
We use the equations of \cite[Theorem 1.2]{Fisher-711} and argue like in the case $n=7$: all we need to show is that if $F$ is the cubic form $v^2w + w^2x + x^2y + y^2z + z^2v$ and $H$ is the determinant of its matrix of second derivatives, the radical of the ideal of $k[v,w,x,y,z]$ generated by $F$, the partial derivatives of $H$, and the function $c_4$ (defined on p. 16 of \emph{op.cit.}) is exactly the ideal $(v,w,x,y,z)$, which is an easy computation with Magma. 

\medskip

\emph{Case where $n=13$:}
We argue as in the case $n=7$, and use the equations of \cite[Theorem 5.10]{Fisher-13} (which are defined over $\Q[a,b]$ for the reduced Weierstrass equation $y^2=x^3+ax+b$, see p. 21 of \emph{op.cit.}) with $\alpha=1$ and those of \cite[Theorem 5.14]{Fisher-13} where $\alpha$ is not a square modulo $13$ (also defined over $\Q[a,b,(4a^3+27b^2)^{-1}]$, see p. 24). What remains to be done is checking that in either case, the radical of the ideal of $R[x_0,x_1,x_2,x_3,x_4,x_5,x_6]$ generated by the cubic forms $f_i$ (for $1 \leq i \leq 7$), $\mathcal{Q}$, $\mfk{c}_6$, $\mathcal{F}$ is exactly $(x_0,x_1,x_2,x_3,x_4,x_5,x_6)$. We can thus assume that $R$ is an algebraically closed field. 

These quantities are obtained by pre-composing an equation of $X(13) \subset \mathbb{P}^6$ (with its $j$-invariant) by a matrix in $\GL{k}$ (see proofs of Theorem 5.10 and 5.14 of \emph{op.cit.}). It is therefore enough to show that, in the notation of \S 5.1 of \emph{op.cit.}, the common vanishing locus of the $f_i$, of $Q$, $F$ and $c_6$ is trivial. This can be directly checked in Magma.\footnote{On the other hand, we have not been able to show the same result directly with the twisted forms using Magma.}  
}

\rem{
We expect that the approach of Proposition \ref{models-of-XEalpha} works for $n=6$ once suitable equations for $X_E^{\alpha}(6)$ and its $j$-invariant map are known. However, the forgetful map $X_E^1(6) \rar X_E^1(3)\simeq \bP^1$ described in \cite[Theorem 4.2.2]{Chen-PhD} is not globally defined, since, even over the base field $\Q(a,b)$, the numerator and the denominator of the fraction have common zeros on $X_E^1(6)$. For $(n,\alpha)=(6,-1)$, \cite[Theorem 1.7.1]{Chen-PhD} only gives the equation of a curve which, while birational to $X_E^{-1}(6)$, is not smooth in general (e.g. for $a=b=1$).}

\subsection{\'Etale group schemes and degeneracy maps} 
\label{subsect:etale+gamma0}

\nott{In this section, we call $\psi$ the multiplicative function defined as follows: if $q > 1$ is a power of the prime $p$, then $\psi(q)=q\left(1+\frac{1}{p}\right)$, and we reserve the notation $\varphi$ for Euler's totient function. Recall that $[\Gamma_0(n)]_{\Z}$ is finite locally free over $\Ell_{\Z}$ of rank $\psi(n)$.}

\prop[YG-Gamma0m]{Let $m,N \geq 1$ be pairwise coprime integers and $G$ be an \'etale $N$-torsion group over some $\Z[1/N]$-algebra $R$. The moduli problem $\P_G(N,\Gamma_0(m))$ is relatively representable of level $Nm$, and finite locally free of rank $r_{N,m}:=|\GL{\Z/N\Z}|\psi(m)$ over $\Ell_R$. If furthermore $N \geq 3$, it is representable by a $R$-scheme $Y_G(N,\Gamma_0(m))$ and the following properties hold:
\begin{enumerate}[noitemsep,label=(\roman*)]
\item\label{PGYm-1} the forgetful map $Y_G(N,\Gamma_0(m)) \rightarrow Y_G(N)$ is finite locally free; post-composing it with $(j,\det)$ defines a finite locally free morphism $ Y_G(N,\Gamma_0(m)) \rightarrow \bA^1_R \times \mrm{Pol}(G)$, also denoted $(j,\det)$, which has rank $r'_{N,m} := \frac{r_{N,m}}{2\varphi(N)}$,
\item\label{PGYm-2} $\det: Y_G(N,\Gamma_0(m)) \rightarrow \mrm{Pol}(G)$ is flat, its geometric fibres are pure of dimension one and Cohen--Macaulay. It is smooth of relative dimension one except at points of prime residue characteristic $p$ such that either $p^2 \mid m$ or $p \mid m$ and their $j$-invariant is supersingular.  
\item\label{PGYm-2b} the same claim holds for the structure map $Y_G(N,\Gamma_0(m)) \rar \Sp{R}$.
\end{enumerate}}

\demo{The moduli problem is relatively representable and finite locally free of rank $r_{N,m}$ over $\Ell_R$, \'etale above $\Z[1/m]$, by Proposition \ref{product-rel-rep}, using \cite[Lemma 4.12.1, Theorem 5.1.1]{KM} and Proposition \ref{PG-elementary}. When $N \geq 3$, it is representable by Proposition \ref{product-rel-rep}. The claims \ref{PGYm-1}, \ref{PGYm-2}, \ref{PGYm-2b} follow from Propositions \ref{PG-elementary} and \ref{struct-gamma0m-smooth} and Corollary \ref{over-Ell-concrete} because the geometric fibres of $\det$ are closed open subschemes of the geometric fibres of the structure map. }

\cor[YG-Gamma0m-pol]{Keep the notations of Proposition \ref{YG-Gamma0m}. Assume that $N \geq 3$ and that $G$ is polarized. Let $u \in \mrm{Pol}(G)(R)$. The moduli problem $\P_G^u(\Gamma_0(m))$ over $R$ is relatively representable, finite locally free of rank $|\SL{\Z/N\Z}|\psi(m)$ over $\Ell_R$, \'etale above $R[1/m]$. It is represented by the closed open subscheme $Y_G^u(N,\Gamma_0(m)) := \det^{-1}(u)$ of $Y_G(N,\Gamma_0(m))$.}

\demo{Since $\P_G^u(\Gamma_0(m))=\P_G^u \times [\Gamma_0(m)]_R$, this follows from Propositions \ref{product-rel-rep} and \ref{XG-pol}.}

\prop[XG-Gamma0m]{Keep the notations of Proposition \ref{YG-Gamma0m}. Assume furthermore that $N \geq 3$ and that $m=m_0m_1$, where $m_0, m_1$ are coprime, $m_0$ is square-free, and $m_1$ is invertible in $R$. Then:
\begin{enumerate}[noitemsep,label=(\roman*)]
\item\label{PGm-2} $\det: Y_G(N,\Gamma_0(m)) \rightarrow \mrm{Pol}(G)$ and the structure map $Y_G(N,\Gamma_0(m)) \rightarrow \Sp{R}$ have reduced geometric fibres. They are smooth of relative dimension one except at points with supersingular $j$-invariant in residue characteristic $p \mid m_0$.
\item\label{PGm-3} if $R$ is regular excellent, $\P_G(N,\Gamma_0(m))$ is smooth at infinity over $\Z[(Nm_1)^{-1}]$, and the formation of its compactified moduli scheme $X_G(N,\Gamma_0(m))$ and cuspidal subscheme commutes with base change.   
\item\label{PGm-4} if $R$ is regular excellent, the map $Y_G(N,\Gamma_0(m)) \rightarrow Y_G(N)$ extends to a finite locally free map $X_G(N,\Gamma_0(m)) \rightarrow X_G(N)$; in particular, $(j,\det)$ also extends to a finite locally free morphism $(\overline{j},\det): X_G(N,\Gamma_0(m)) \rightarrow \bP^1_R \times \mrm{Pol}(G)$ of rank $r'_{N,m}$. The claim of \ref{PGm-2} holds verbatim after replacing $Y_G$ with $X_G$. 
\end{enumerate}}

\rem{When $G$ is constant, we write $Y(N,\Gamma_0(m))$ (resp. $X(N,\Gamma_0(m))$) for $Y_G(N,\Gamma_0(m))$ (resp. $X_G(N,\Gamma_0(m))$).}

\demo{The claim \ref{PGm-2} is simply a rewriting of the claims \ref{PGYm-2} and \ref{PGYm-2b} of Proposition \ref{YG-Gamma0m}. So assume from now on that $R$ is regular excellent and let us show that $\P_G(N,\Gamma_0(m))$ is smooth at infinity over $\Z[(Nm_1)^{-1}]$. By Proposition \ref{smooth-at-infinity-descends}, this can be checked \'etale-locally over $\Sp{R}$, so we can assume that $G$ is constant. Then $\P_G(N,\Gamma_0(m))$ is the base change of $([\Gamma(N)] \times [\Gamma_0(m)])_{\Z\left[\frac{1}{Nm_1}\right]}$, it is smooth at infinity over $\Z[1/Nm_1]$ by Proposition \ref{struct-gamma0m-smooth}. Then \ref{PGm-3} is Proposition \ref{smooth-at-infinity-base-change} and \ref{PGm-4} follows formally from \ref{PGm-2} and Corollary \ref{flatness-at-cusps}.}

\cor[XG-Gamma0m-pol]{Keep the notations of Proposition \ref{XG-Gamma0m}. Assume that $R$ is regular and excellent and that $G$ is polarized. Let $u \in \mrm{Pol}(G)(R)$. Then $\P_G^u(\Gamma_0(m))$ is smooth at infinity over $\Z[1/Nm]$ and the closed open subscheme $X_G^u(N,\Gamma_0(m))=\det^{-1}(u)$ of $X_G(N,\Gamma_0(m))$ identifies with the compactified moduli scheme $\CMPp{\P_G^u(\Gamma_0(m))}$ as a $\CMPp{\P_G(N,\Gamma_0(m))}$-scheme. }

\demo{This is proved as Proposition \ref{XG-pol}.}

\prop[XG-Gamma0m-connectedness]{Keep the notations of Proposition \ref{XG-Gamma0m}. Then $\det: Y_G(N,\Gamma_0(m)) \rightarrow \mrm{Pol}(G)$ has geometrically connected fibres, and the geometric fibres of the structure map of $Y_G(N,\Gamma_0(m))$ to $\Sp{R}$ have $\varphi(N)$ connected components. When $R$ is regular excellent, the same claims hold after replacing $Y_G$ by $X_G$.}

\demo{When $G$ is constant and $R=\Z[1/N]$, $\det: X_G(N,\Gamma_0(m)) \rightarrow \mrm{Pol}(G)$ has geometrically connected fibres above $\Z[1/Nm]$ by \cite[Corollary 10.9.5 (2)]{KM}, hence so it has geometrically connected fibres by e.g. \cite[Lemma 0ENE]{Stacks}. The case of $Y_G(N,\Gamma_0(m))$ follows since the cuspidal points are non-isolated and normal in the geometric fibres of $X_G(N,\Gamma_0(m)) \rightarrow \mrm{Pol}(G)$. By base change, the result holds when $G$ is constant for any $R$. In general, $G$ is \'etale-locally constant, so the conclusion follows by descent. 

The case of the structure maps to $\Sp{R}$ is treated as in the proof of Corollary \ref{counting-connected-components-XG}. }

\cor[YG-Gamma0m-connectedness]{Let $N \geq 3, m \geq 1$ be coprime integers. Let $G$ be an \'etale $N$-torsion group over a $\Z[1/N]$-algebra $R$. Then $\det: Y_G(N,\Gamma_0(m)) \rightarrow \mrm{Pol}(G)$ has connected geometric fibres. In particular, if $G$ is polarized and $u \in \mrm{Pol}(G)(R)$, then $Y_G^u(N,\Gamma_0(m)) \rightarrow \Sp{R}$ has geometrically connected fibres. }

\demo{We may assume that $R$ is an algebraically closed field $k$ of characteristic $p$. If $p=0$ or $p^2 \nmid m$, this is Proposition \ref{XG-Gamma0m-connectedness}, so we may assume that $p$ is prime and $m=p^rm'$ with $r \geq 2$ and $p \nmid m'$. Let $u \in \mrm{Pol}(G)(k)$. It suffices to show that $\M := \MPp{\P_G^u \times [\Gamma_0(m)]_k}$ is connected. Since $\P' := \P_G^u \times [\Gamma_0(m')]_k$ is finite \'etale over $\Ell_k$, \cite[Theorem 13.4.7, (13.5.6)]{KM} implies that $\M$ is the union of closed subschemes homeomorphic to $\MPp{\P'}$ (hence irreducible by Proposition \ref{XG-Gamma0m-connectedness} and Corollary \ref{over-Ell-concrete}) that share a common point, so we are done. }

\cor[flat-deg-constant]{In Corollaries \ref{flat-deg} and \ref{flat-deg-infty}, the map $D_{d,t}$ is finite locally free of constant rank $\frac{\psi(m)}{\psi(t)}$.}

\demo{It is enough to treat the case of Corollary \ref{flat-deg}, since the moduli scheme is dense in its compactification. As in the proof of Corollary \ref{flat-deg}, we may assume that the auxiliary moduli problem $\P$ is $[\Gamma(\ell)]_{\Z[1/\ell]}$ for some prime $\ell \nmid 2m$. Then, for every $d \mid m$, $Y(\ell,\Gamma_0(d)) \rightarrow \Sp{\Z[\frac{1}{\ell},\zeta_{\ell}]}$ is flat with geometrically connected fibres, so $Y(\ell,\Gamma_0(d))$ is connected and every $D_{d,t}$ has constant rank $r(m,d,t,\ell)$. When $d=1$, $D_{d,t}$ is a morphism of finite locally free $Y(\ell)$-schemes, the domain of rank $\psi(m)$ and the target of rank $\psi(t)$, so $r(m,1,t,\ell)=\frac{\psi(m)}{\psi(t)}$. In general, the equation $D_{1,1} \circ D_{d,t} = D_{1,1} \circ w_d \circ D_{1,d}$ implies $\psi(t)r(m,d,t,\ell)=\psi(d)\frac{\psi(m)}{\psi(d)}$ and we are done.}  

\prop[ext-degn]{Let $N \geq 3$, $R$ be a $\Z[1/N]$-algebra and $G$ be an \'etale $N$-torsion group over $R$. The collection of $Y_G(N,\Gamma_0(m))$, taken over integers $m=m_0m_1$ coprime to $N$ with $m_0,m_1 \geq 1$ coprime, $m_0$ square-free and $m_1 \in R^{\times}$, is endowed with finite locally free degeneracy maps $D_{d,t}$ (for $dt \mid m$) of constant rank $\frac{\psi(m)}{\psi(t)}$ and Atkin--Lehner automorphisms as in Definitions \ref{degeneracies} and \ref{AL-autod}, and an action of $\widehat{\Z}^{\times}$ which factors through $(\Z/N\Z)^{\times}$. The action of $\widehat{\Z}^{\times}$ commutes to the $D_{d,t}$ and to the $w_d$. The relations proved throughout Section \ref{gamma0m-structure} hold\footnote{This applies to Lemma \ref{elem-degeneracies}, the first item of Lemma \ref{AL-elem}, Proposition \ref{AL-auto}, and the fact that the diagrams in Propositions \ref{degn-cartesian} and \ref{two-degn-cartesian} commute.}.

If $R$ is furthermore regular excellent, those maps extend to the compactified moduli schemes, and the same claims hold for their extensions.}

\demo{This is a direct translation of the results of Section \ref{gamma0m-structure} (note Proposition \ref{basic-moduli-level} for the $\widehat{\Z}^{\times}$-action), Corollaries \ref{flat-deg-infty}, \ref{flat-deg-constant}, as well as Propositions \ref{XG-Gamma0m} and \ref{compactification-functor}. }

\prop[degen-AL-Weil]{Let $m_0,m_1,m,N,R,G$ be as in Proposition \ref{XG-Gamma0m}.
\begin{itemize}[noitemsep,label=$-$]
\item Let $d,t \geq 1$ be such that $dt \mid m$. Then $\det \circ D_{d,t} = \underline{d} \circ \det$.
\item Let $d \geq 1$ be a divisor of $m$ coprime to $m/d$. Then $\det \circ w_d = \underline{d} \circ \det$.
\item When $R$ is regular excellent, both identities hold over $X_G(N,\Gamma_0(m))$.
\end{itemize}}

\demo{The first two items are computations, the last one follows by Proposition \ref{compactification-functor}. }

\prop[degen-AL-GL2]{Keep the notations of Proposition \ref{XG-Gamma0m}; assume that $P,Q \in G(R)$ form a basis of $G$. Through the isomorphism of moduli problems of level $N$ 
\vspace{-8pt}\[\P_G(N) \Rightarrow [\Gamma(N)]_R,\qquad (\iota : G_S \overset{\sim}{\rightarrow} E[N]) \mapsto (\iota(P),\iota(Q)),\vspace{-8pt}\] $\GL{\Z/N\Z}$ acts on the left on the level $N$ structure of the moduli problem $\P_G(N)$. Thus:
\begin{itemize}[noitemsep,label=\tiny$\bullet$]
\item $\GL{\Z/N\Z}$ acts on every $Y(N,\Gamma_0(m))$ and commutes to the degeneracy maps and Atkin--Lehner automorphisms. If $a \in \widehat{\Z}^{\times}$, $[a]$ and $\overline{a} \cdot I_2 \in \GL{\Z/N\Z}$ agree in $\mathrm{Aut}(Y(N,\Gamma_0(m)))$. 
\item if moreover $R$ is regular excellent, this action extends to $X_G(N,\Gamma_0(m))$. The claims in the previous item hold verbatim after replacing moduli schemes by their compactifications.  
\end{itemize}}

\demo{The rule is clearly an isomorphism of moduli problems of level $N$. The first item is then a consequence of Lemma \ref{elem-degeneracies} and the remarks following Definition \ref{AL-autod}, and the action of $\widehat{\Z}^{\times}$ is an explicit computation. The second item follows from Proposition \ref{compactification-functor}. }

\bigskip

\prop[degen-forgetful-maps]{Let $N,e \geq 3$ be integers with $e \mid N$ and $G$ be an \'etale $N$-torsion group over a $\Z[1/N]$-algebra $R$. For every $m_0,m_1 \geq 1$ such that $N,m_0,m_1$ are pairwise coprime, $m_1 \in R^{\times}$ and $m_0$ is square-free, the forgetful map $\P_G(N) \rightarrow \P_{G[e]}(e)$ induces a finite \'etale forgetful map $Y_G(N,\Gamma_0(m)) \rightarrow Y_{G[e]}(e,\Gamma_0(m))$, where $m=m_0m_1$. This forgetful map commutes to all degeneracy maps, Atkin--Lehner automorphisms and to the action of $\widehat{\Z}^{\times}$. 

If $R$ is furthermore regular excellent, this forgetful map extends to a finite locally free morphism $X_G(N,\Gamma_0(m)) \rightarrow X_{G[e]}(e,\Gamma_0(m))$ which commutes to the compactifications of the degeneracy maps, Atkin--Lehner automorphisms and of the action of $\widehat{\Z}^{\times}$.  
}

\demo{That the forgetful map is finite \'etale is proved as in Proposition \ref{forgetful-maps}, and its compactification is locally free by Proposition \ref{flatness-at-cusps}. The forgetful morphism $\P_G(N) \Rightarrow \P_{G[e]}(e)$ respects the level $N$ structure, so the remaining claims are proved as in Proposition \ref{degen-AL-GL2}. }

\subsection{Hecke correspondences for \'etale polarized group schemes}

In this section, $R$ denotes a regular excellent domain, $N \geq 3$ is an integer invertible in $R$, and $G$ is a polarized\footnote{It is possible to remove this assumption by slightly extending the notion of Jacobian to non-geometrically connected curves (cf. e.g. \cite[Section 1.5.2]{Studnia-thesis} or \cite[Appendix A]{RSZB}), but this is not needed here.} \'etale $N$-torsion group over $R$. Let us first recall elementary facts about pushforward functoriality on Picard groups. 

\defi{(cf. \cite[II (6.5.2), (6.5.5); $\text{IV}_4$ (21.5.7.2)]{EGA}) If $f: X \rightarrow Y$ is a finite locally free morphism of schemes, $f$ induces a norm homomorphism $f_{\ast}: \operatorname{Pic}(X) \rightarrow \operatorname{Pic}(Y)$. The rule $X \mapsto \operatorname{Pic}(X), f \mapsto f_{\ast}$ is a functor from the category of schemes with finite locally free morphisms to the category of abelian groups. }

\rem{This definition recovers the intuitive notion of direct image on cycles of codimension $1$, see \cite[$\text{IV}_4$ (21.10.17)]{EGA}.}

\rem{Let $X,Y$ be proper connected one-dimensional schemes of finite type over a field $k$, so line bundles on $X,Y$ have a well-defined degree (e.g. \cite[Definition 0AYR]{Stacks}). Let $f: X \rightarrow Y$ be a finite locally free morphism of $k$-schemes, then $f_{\ast}$ sends $\operatorname{Pic}^0(X)$ to $\operatorname{Pic}^0(Y)$.}

\prop[base-change-pseudo-cartesian]{Consider the following commutative diagram of schemes, where all the maps are finite locally free:
\vspace{-7pt}\[
\begin{tikzcd}[ampersand replacement=\&]
Y \arrow{r}{f} \arrow{d}{v} \& X \arrow{d}{u}\\
Y' \arrow{r}{f'} \& X'
\end{tikzcd}\]
The two morphisms 
\vspace{-4pt}\[u^{\ast}\circ (f')_{\ast},\; f_{\ast} \circ v^{\ast}: \operatorname{Pic}(Y') \rar \operatorname{Pic}(X)\vspace{-4pt}\] are equal if, for any affine open subscheme $U \subset X'$ and any $b \in \OO_{Y'}((f')^{-1}(U))$, one has $u^{\sharp}(N_{Y'/X'}(b))=N_{Y/X}(v^{\sharp}(b))$. This assumption is verified if the diagram is Cartesian above a scheme-theoretically dense open subscheme $U \subset X'$. }

It is not hard to find a proof using the ideas of \cite[II (6.5.8)]{EGA}; the reader can also find a self-contained argument in e.g. \cite[Proposition A.4.10]{Studnia-thesis}.

\prop[hecke-operator-defined]{Let $\alpha \in \mrm{Pol}(G)(R)$ and $q \nmid N$ be a prime. Recall the degeneracy maps 
\vspace{-7pt}\[D_{1,1}: X_G^{\alpha}(N,\Gamma_0(q)) \rar X_G^{\alpha}(N),\qquad D_{q,1}:  X_G^{\alpha}(N,\Gamma_0(q)) \rar X_G^{q\cdot \alpha}(N).\vspace{-5pt}\]
The rule $(D_{q,1})_{\ast}D_{1,1}^{\ast}$ defines a morphism of abelian pre-sheaves on $\Sch_R$
\vspace{-7pt}\[\left[S \mapsto \operatorname{Pic}(X_G^{\alpha}(N)\times_R S)\right] \rar \left[S \mapsto \operatorname{Pic}(X_G^{q \cdot \alpha}(N)\times_R S)\right]\vspace{-7pt}\] which preserves the subfunctors $\operatorname{Pic}^0$. Thus, its sheafification induces a morphism of Abelian schemes 
\vspace{-7pt}\[T_q: \mrm{Jac}(X_G^{\alpha}(N)) \rar \mrm{Jac}(X_G^{q\cdot \alpha}(N)).\vspace{-7pt}\] }

\demo{By Corollary \ref{XG-pol-prop}, $X_G^{\alpha}(N), X_G^{q\cdot \alpha}(N)$ are smooth proper $R$-schemes of relative dimension one with geometrically connected fibres, so their $\operatorname{Pic}^0$ functors are representable by Abelian schemes by \cite[\S 9.4 Proposition 4]{BLR}. By Propositions \ref{ext-degn} and \ref{degen-AL-Weil}, $D_{1,1},D_{q,1}$ are (well-defined) finite locally free morphisms, so the morphisms of pre-sheaves on $\Sch_R$
\vspace{-7pt}
\begin{align*}
(D_{1,1})^{\ast}:\;& \operatorname{Pic}\big(X_G^{\alpha}(N)\times_R -\big) \; \longrightarrow \; \operatorname{Pic}\big(X_G^{\alpha}(N,\Gamma_0(q))\times_R -\big)\\
(D_{q,1})_{\ast}:\;& \operatorname{Pic}\big(X_G^{\alpha}(N,\Gamma_0(q))\times_R -\big) \; \longrightarrow \; \operatorname{Pic}\big(X_G^{q \cdot \alpha}(N)\times_R -\big)
\end{align*}
\vspace{-18pt}

\noindent are well-defined. Moreover, since $X_G^{\alpha}(N,\Gamma_0(q)) \rar \Sp{R}$ also has geometrically connected fibres by Proposition \ref{XG-Gamma0m-connectedness}, these morphisms preserve the degree-zero subfunctors by the above and the results of \cite[Section 0AQV]{Stacks}. 
}

\defi{The operator $T_q$ of Proposition \ref{hecke-operator-defined} is the \emph{Hecke operator of index $q$}. It is a morphism of Abelian schemes. }

\prop[hecke-operators-commute]{Let $\alpha \in \mrm{Pol}(G)(R)$, and $q,r \nmid N$ be two distinct primes. The two operators 
\vspace{-7pt}\[T_q\circ T_r, T_r\circ T_q: \mrm{Jac}(X_G^{\alpha}(N)) \rar \mrm{Jac}(X_G^{qr\cdot \alpha}(N))\vspace{-7pt}\]
are equal. }

\demo{Consider the following commutative diagram for any $R$-scheme $S$, where all maps are finite locally free:
\vspace{-7pt}\[
\begin{tikzcd}[ampersand replacement=\&]
X_G^{\alpha}(N,\Gamma_0(qr))_S \arrow{r}{D_{1,r}} \arrow{d}{D_{r,q}}\arrow{rd}{D_{r,1}}\& X_G^{\alpha}(N,\Gamma_0(r))_S\arrow{d}{D_{r,1}}\\
X_G^{r\cdot \alpha}(N,\Gamma_0(q))_S \arrow{r}{D_{1,1}} \& X_G^{r \cdot \alpha}(N)_S
\end{tikzcd}
\vspace{-4pt}\]
Let us show that the two maps 
\vspace{-7pt}\[(D_{1,1})^{\ast}(D_{r,1})_{\ast}, (D_{r,q})_{\ast}(D_{1,r})^{\ast}: \operatorname{Pic}(X_G^{\alpha}(N,\Gamma_0(r))_S) \rar \operatorname{Pic}(X_G^{r\cdot \alpha}(N,\Gamma_0(q))_S)\vspace{-7pt}\] agree. Indeed, we know that this diagram is Cartesian above $Y_G^{\alpha}(N)$ by Proposition \ref{two-degn-cartesian}, so by Proposition \ref{base-change-pseudo-cartesian} it is enough to see that $Y_G^{\alpha}(N) \rightarrow X_G^{\alpha}(N)$ is universally scheme-theoretically dominant (over $\Sp{R}$). This follows from Corollaries \ref{compactification-is-dense} and \ref{XG-pol-prop}. 

Thus, as morphisms between presheaves of Picard groups, one has 
\vspace{-7pt}
\begin{align*}
(D_{q,1})_{\ast}(D_{1,1})^{\ast}(D_{r,1})_{\ast}(D_{1,1})^{\ast} &= (D_{q,1})_{\ast}(D_{r,q})_{\ast}(D_{1,r})^{\ast}(D_{1,1})^{\ast}\\
&= (D_{q,1}\circ D_{r,q})_{\ast}(D_{1,1}\circ D_{1,r})^{\ast} = (D_{qr,1})_{\ast}(D_{1,1})^{\ast},
\end{align*}
\vspace{-20pt}

\noindent where the final $D_{1,1}$ is the (forgetful) degeneracy map $X_G^{\alpha}(N,\Gamma_0(qr)) \rar X_G^{\alpha}(N)$. In particular, this map does not depend on the order of $q,r$, and the conclusion follows. }

\rem{To show that these Hecke operators commute, it is necessary to work with the ``norm-commuting diagrams'' of Proposition \ref{base-change-pseudo-cartesian} instead of genuine Cartesian diagrams. Indeed, let $N \geq 3$ be an integer and $p,q$ be distinct primes not dividing $N$. Even for the constant group scheme $G=(\Z/N\Z)^{\oplus 2}$ over $\C$ (equipped with some choice of polarization which we omit from the notation), the diagram 
\vspace{-8pt}\[
\begin{tikzcd}[ampersand replacement=\&]
X_G(N,\Gamma_0(pq)) \arrow{r}{D_{1,p}} \arrow{d}{D_{1,q}} \& X_G(N,\Gamma_0(p)) \arrow{d}{D_{1,1}}\\
X_G(N,\Gamma_0(q)) \arrow{r}{f'} \& X_G(N)
\end{tikzcd}
\vspace{-7pt}\]
is only Cartesian above $Y_G(N)$. Indeed, the rightmost and bottom maps are both ramified above the cusp $0$, so their Cartesian product cannot be normal at a point above $0$, but $X_G(N,\Gamma_0(pq))$ is normal.} 

\rem{The same arguments can be used to show that the Hecke operators commute to the push-forwards and the pull-backs of the forgetful maps defined in Proposition \ref{forgetful-maps}. }

\bigskip

\subsection{Smoothness of $Y_G(N)$ when $N$ is not invertible}
\label{subsect:smoothness-general}

\defi[N-tors-group]{Let $N \geq 1$ be an integer and $S$ be a scheme. A finite locally free commutative $S$-group scheme $G$ killed by $N$ is a \emph{$N$-torsion group} if, for every $p \mid N$, the finite locally free $S$-group scheme $G[p^{v_p(N)}]$ is a $\text{BTT}_{v_p(N)}$ of height $2$ in the sense of \cite[D\'efinition 1.1]{Illusie}.}

\rem{Let $S$ be a scheme and $G$ be a finite commutative $S$-group scheme killed by $N$. 
\begin{itemize}[noitemsep,label=\tiny$\bullet$]
\item If $G$ is a $N$-torsion group over $S$, for any $d \mid N$, $G[d]$ is finite locally free over $S$ of rank $d^2$. 
\item If $S'$ is a $S$-scheme and $G$ is a $N$-torsion group over $S$, $G_{S'}$ is a $N$-torsion group over $S'$.
\item If $N$ is invertible in $S$, $G$ is a $N$-torsion group over $S$ if and only if it is \'etale-locally isomorphic to $(\Z/N\Z)^{\oplus 2}$. 
\item If $E/S$ is an elliptic curve, $E[N]$ is a $N$-torsion group over $S$ (by \cite[Exemple 1.2]{Illusie}).
\end{itemize}}

Let us give an equivalent to Definition \ref{N-tors-group}. We will need the following lemma. 

\lem[reform-FV]{Let $S$ be a scheme, $f: G \rightarrow H$ be a morphism of finite $S$-group schemes of finite presentation. Suppose that $G$ is locally free over $S$. Then $f$ is surjective in the fppf topology if and only if it is faithfully flat. In this case, $H$ is finite locally free.}

\demo{The ``if'' is \cite[$\text{IV}_4$ (17.16.2) (i)]{EGA}. Conversely, assume that $G \rightarrow H$ is a surjective morphism of abelian fppf-sheaves. For any $s \in S$, $G_s \rightarrow H_s$ is a surjective morphism of abelian fppf sheaves (represented by finite $\kappa(s)$-group schemes), so it is faithfully flat by \cite[(2.3.3.2)]{Anantharaman}. Hence $f$ is flat by the fibre criterion for flatness.  }

\cor[alternate-def-N-torsion]{Let $G$ be a finite commutative $S$-group scheme of finite presentation killed by $N \geq 1$. Then $G$ is a $N$-torsion group if and only if the following two conditions hold:
\begin{itemize}[noitemsep,label=\tiny$\bullet$]
\item For any $d \mid N$, $G \overset{\cdot d}{\rightarrow} G[N/d]$ is finite locally free of rank $d^2$,
\item Let $p \mid N$ be prime such that $p^2 \nmid N$. The morphism $\Fr: G[p]_{R/(p)} \rar \ker{\mrm{Ver}}$ is faithfully flat, where   
\vspace{-11pt}\[\Fr: G[p]_{R/(p)} \rar G[p]_{R/(p)}^{(p)}, \qquad \mrm{Ver}: G[p]_{R/(p)}^{(p)} \rar G[p]_{R/(p)}\vspace{-9pt}\]
 are the relative Frobenius and Verschiebung respectively (cf. e.g. \cite[Exp. $\mrm{VII}_{\mrm{A}}$, \S 4]{SGA3}).
\end{itemize}
}

\cor[weN-descent]{Let $N \geq 1$ be an integer, $S$ be a scheme, and $G$ be a finite $S$-group scheme of finite presentation.
\begin{enumerate}[noitemsep,label=(\roman*)]
\item \label{weNd-1} Let $S'$ be a flat surjective $S$-scheme. If $G_{S'}$ is a $N$-torsion group over $S'$, $G$ is a $N$-torsion group over $S$. 
\item \label{weNd-1b} Assume that $G$ is locally free over $S$, commutative and killed by $N$. Let $S'$ be a surjective $S$-scheme. If $G_{S'}$ is a $N$-torsion group over $S'$, then $G$ is a $N$-torsion group over $S$. 
\item \label{weNd-2} Suppose that $S=\Sp{R}$ and that $G$ is a $N$-torsion group. Then there is a $N$-torsion group $G_0$ over a Noetherian subring $R_0$ of $R$ such that $G \simeq (G_0)_R$. 
\end{enumerate}
}

\demo{\ref{weNd-1}: by Corollary \ref{alternate-def-N-torsion} and properties of flatness, we may assume that $S',S$ are affine and $\OO(S) \rightarrow \OO(S')$ is a faithfully flat map of local rings. The conclusion follows by descent. 

\ref{weNd-1b}: For every $d \mid N$ and every $s \in S$, the finite map of finite presentation $G_s \overset{\cdot d}{\rightarrow} G[N/d]_s$ becomes locally free of rank $d^2$ after the (flat) base change $S'_s \rar \Sp{\kappa(s)}$, so it is locally free of rank $d^2$. By the fibre criterion for flatness, $G \overset{\cdot d}{\rightarrow} G[N/d]$ is finite locally free of rank $d^2$. Let $p \mid N$ be prime such that $p^2 \nmid N$, then $G[p]$ is finite locally free over $S$ by prime factorization, and the same argument shows that $\Fr: G[p] \rightarrow \ker{\mrm{Ver}}$ is flat. Hence $G$ is a $N$-torsion group over $S$. 

\ref{weNd-2}: This follows by Noetherian approximation \cite[$\text{IV}_3$ (8.8.2), (8.10.5) (x), (11.2.6.1)]{EGA}, since $R$ is the (directed) direct limit of its finitely generated $\Z$-subalgebras (which are Noetherian). }

\lem[lift-pdiv-serre-tate]{Let $p$ be a prime and $n,r \geq 1$. Let $R$ be a ring of characteristic $p^n$, $I \subset R$ be a nilpotent ideal and $G$ be a $p^r$-torsion group over $R$. Let $E'$ be an elliptic curve over $R' := R/I$ and $\iota': G_{R/I} \rar E'[p^r]$ be an isomorphism of $R'$-group schemes. There exists an elliptic curve $E/R$ and a morphism $\lambda: E'/R' \rar E/R$ in the category $\Ell_R$ such that $\lambda \circ \iota': G_{R/I} \rar E_{R/I}[p^r]$ lifts to an isomorphism $\iota: G \rar E[p^r]$.}

\demo{The map $E'[p^r] \overset{(\iota')^{-1}}{\rar} G_{R'} \rar G_R$ is a deformation of the truncated Barsotti--Tate group $E'[p^r]$ from $R'$ to $R$. 

By \cite[Th\'eor\`eme 4.4 (f)]{Illusie} (in the notation of \emph{loc.cit.}, one has $S_0=\Sp{R'/(p)}$, $S=\Sp{R'}$ and $S'=\Sp{R}$), it comes from a deformation from $R'$ to $R$ of the Barsotti--Tate group\footnote{see \cite[D\'efinition 1.5]{Illusie}.} $E'[p^{\infty}]$. That is, there is a Barsotti--Tate group $D$ over $R$ and isomorphisms $f: E'[p^{\infty}] \rar D_{R'}$ of Barsotti--Tate groups over $R'$, $j: G_R \rar D[p^r]$ of fppf abelian sheaves over $\Sp{R}$ such that the following diagram commutes: 
\vspace{-7pt}\[
\begin{tikzcd}[ampersand replacement=\&]
E'[p^r] \arrow{r}{(\iota')^{-1}}\arrow{d}{\subset} \& G_{R'} \arrow{d}{j_{R'}} \arrow{r}\& G\arrow{d}{j}\\
E'[p^{\infty}] \arrow{r}{f} \& D_{R'} \arrow{r} \& D
\end{tikzcd}
\vspace{-7pt}\]

Every $D[p^r]$ is a finite locally free $R$-group scheme by \cite[Sorites 1.6, Remarque 1.3 (c)]{Illusie}. By the Serre--Tate theorem \cite[Theorem 2.9.1]{KM}, there exists an elliptic curve $E/R$ endowed with isomorphisms $\varphi: E' \rar E_{R'}$ and $\delta: D \rar E[p^{\infty}]$ such that $\varphi[p^{\infty}] = \delta_{R'} \circ f$. Then $\iota = \delta[p^r] \circ j: G \rar E[p^r]$ is an isomorphism of group schemes and one has 
\vspace{-7pt}\[\varphi^{-1}\circ \iota_{R'} = \varphi^{-1}\circ \delta_{R'}[p^r] \circ j_{R'}=\varphi^{-1}\circ \delta_{R'}\circ f[p^r] \circ \iota'=\iota',\vspace{-7pt}\]
so we are done.}

\prop[general-smoothness-PG]{Let $N \geq 1$ be an integer and $G$ be a $N$-torsion group over a ring $R$. Let $\P$ be a relatively representable moduli problem over $R$, \'etale over $\Ell_R$, which is representable if $N \leq 2$. Then $\MPp{\P,\P_G(N)}$ is smooth over $R$ of relative dimension one.}

\demo{We may assume that $N>1$ by Corollary \ref{over-Ell-concrete}. The moduli problem $\P \times \P_G(N)$ is representable by Propositions \ref{product-rel-rep} and \ref{PG-representable-when-good}. The claim is Zariski-local over $R$, so we may assume that we can write $N=p^rN'$, where $p \nmid N'$ is a prime, $r \geq 1$, and $N'\ell$ is invertible in $R$ for some prime $\ell \nmid 2N$. By Proposition \ref{PG-elementary}, we can replace $(\P,\P_G(N))$ with $(\P \times \P_{G[N']}(N'),\P_{G[p^r]}(p^r))$ and thus assume that $N=p^r$.  
Arguing as in the proof of Corollary \ref{flat-deg} (using \cite[$\text{IV}_4$ (17.7.7), (17.10.3)]{EGA}), we may assume that $\P=[\Gamma(\ell)]$. Then $\MPp{\P,\P_G(N)}$ is affine of finite presentation over $\MPp{[\Gamma(\ell)]}$, so it is of finite presentation over $R$. By Corollary \ref{weN-descent} \ref{weNd-2}, we may assume that $R$ is Noetherian.  

Suppose first that $(R,\mfk{m})$ is Artinian with residue characteristic $p$. Let $A$ be a $R$-algebra and $I \subset A$ be an ideal with $I^2=0$. Let $x_0 \in \MPp{\P,\P_G(p^r)}(A/I)$: $x_0$ corresponds to an elliptic curve $F$ over $A' := A/I$, $\alpha' \in \P(F/A')$ and an isomorphism $\iota': G_{A/I} \rightarrow F[p^r]$. By Lemma \ref{lift-pdiv-serre-tate}, we can find an elliptic curve $E/A$, a morphism $\varphi: F/A' \rar E/A$ in $\Ell_R$ and an isomorphism $\iota_p: G_A \rar E[p^r]$ such that the following diagram commutes. 
\vspace{-7pt}\[
\begin{tikzcd}[ampersand replacement=\&]
G_{A/I} \arrow{r}{\iota'} \arrow{d} \& F[p^r] \arrow{d}{\varphi}\\
G_A \arrow{r}{\iota} \& E[p^r]
\end{tikzcd}
\vspace{-7pt}\]

Because $\P$ is \'etale, $(\P\varphi): \P(E/A) \rar \P(F/A')$ is a bijection, so there exists $\alpha \in \P(E/A)$ such that $(\P\varphi)\alpha=\alpha'$. Then $x=(E/A, \varphi, \iota) \in \MPp{\P,\P_G(N)}(A)$ lifts $x_0$. Thus $\MPp{\P,\P_G(N)}$ is locally of finite presentation and formally smooth over $R$, hence smooth (\cite[$\text{IV}_4$ (17.5.2)]{EGA}). 

Let us come back to the general case. Let $n \geq 1$, $\mfk{p} \subset R$ be a prime ideal and $R'=R_{\mfk{p}}/\mfk{p}^n$. If $p$ is not invertible in $R'$, then $R'$ is an Artinian local ring of residue characteristic $p$, so by the previous paragraph $\MPp{\P,\P_G(N)}_{R'}$ is smooth over $R'$. If $p$ is invertible in $R'$, the moduli problem $(\P \times \P_G(N))_{R'}$ is affine \'etale over $\Ell_{R'}$ by Proposition \ref{PG-elementary}, so $\MPp{\P, \P_G(N)}_{R'} \rightarrow \Sp{R'}$ is smooth. Thus the morphism $\MPp{\P , \P_G(N)} \rightarrow \Sp{R}$ has smooth fibres, and it is flat at every point by \cite[Exp. IV, Th. 5.6]{SGA1}, so it is smooth. 

Let us check that $X := \MPp{\P,\P_G(N)}$ is smooth of relative dimension one over $R$. By \cite[$\text{IV}_4$ (17.16.3)]{EGA}, working \'etale-locally over $R$, we may assume that $X$ has a $R$-point, i.e. that $G \simeq E[p^r]$ for some elliptic curve $E/R$ with full $\ell$-torsion. Then $E/R$ is a base change of the elliptic curve representing $[\Gamma(\ell)]_{\Z[1/\ell]}$, so we may assume that $R$ is smooth over $\Z[1/\ell]$ of relative dimension one. Above $R[1/p]$, $[\Gamma(\ell)]_R \times \P_{E[p^r]}(p^r)$ is \'etale, so $X$ is smooth of relative dimension one over $R$ at every point of residue characteristic not $p$. Since $R$ is flat over $\Z$ and $X$ is smooth over $R$, $X_{R[1/p]}$ is dense in $X$, so we are done by \cite[$\text{IV}_4$ (17.10.2)]{EGA}.} 

\cor[general-smoothness-equivalence]{Let $R$ be a ring and $N \geq 1$ be an integer. Let $G$ be a finite locally free commutative $R$-group scheme killed by $N$. Let $\P$ be an affine \'etale relatively representable moduli problem over $R$, which is representable if $N \leq 2$. The following are equivalent:
\begin{enumerate}[noitemsep,label=(\roman*)]
\item \label{repeq-1} There is a flat surjective $R$-scheme $S$ and an elliptic curve $E/S$ such that $G_S$ is isomorphic to $E[N]$ and $\P(E/S) \neq \emptyset$. 
\item \label{repeq-2} There is a faithfully flat \'etale $R$-algebra $R'$ and an elliptic curve $E/R'$ such that $G_{R'}$ is isomorphic to $E[N]$ and $\P_{E/R'}$ is a surjective $R'$-scheme. 
\item \label{repeq-3} $\MPp{\P,\P_G(N)}$ is smooth surjective over $R$ of relative dimension one.
\item \label{repeq-4} $\MPp{\P,\P_G(N)}$ is flat surjective over $R$. 
\item \label{repeq-5} $\MPp{\P,\P_G(N)}$ is surjective over $R$. 
\end{enumerate}
Any of these conditions implies that $G$ is a $N$-torsion group.
}

\demo{By Propositions \ref{product-rel-rep}, \ref{PG-representable-when-good} and Corollary \ref{over-Ell-concrete}, the moduli problem $\P \times \P_G(N)$ is representable and affine of finite presentation, so $\MPp{\P,\P_G(N)}$ is an affine $R$-scheme of finite presentation. The implications \ref{repeq-3} $\Rightarrow$ \ref{repeq-2} $\Rightarrow$ \ref{repeq-1} follow from \cite[$\text{IV}_4$ (17.16.3) (ii)]{EGA}; \ref{repeq-3} clearly implies \ref{repeq-4}; \ref{repeq-4} implies \ref{repeq-1} by considering the universal elliptic curve; and, \ref{repeq-1} implies \ref{repeq-5} since $E/S$ defines a $S$-point of $\MPp{\P,\P_G(N)}$. Finally, if \ref{repeq-5} holds, $G$ is a $N$-torsion group by Corollary \ref{weN-descent} \ref{weNd-1b}, and therefore \ref{repeq-3} holds by Proposition \ref{general-smoothness-PG}.} %

\subsection{Quasi-finiteness of $\P_G(N)$ and ordinary group schemes}
\label{subsect:ordinary}

Let $G$ be a finite locally free commutative group scheme over a ring $R$ killed by some integer $N \geq 1$. Assume for simplicity that there is a flat surjective $R$-scheme $S$ such that $G_S$ is the $N$-torsion subgroup of some elliptic curve $E/S$. The main thread of this section is to find under what conditions on $G$ the moduli problem $\P_G(N)$ is quasi-finite or flat over $\Ell_R$. It turns out (by Proposition \ref{quasifinite-ordinary-flat}) that these conditions are equivalent to the fact that all the specializations of $E$ in characteristic dividing $N$ are ordinary, i.e. that the specializations of $G$ at every geometric point of $\Sp{R}$ are isomorphic to $\mu_N \times \Z/N\Z$. A secondary goal of this section is to see that this condition is already strong enough to imply that $G$ is a $N$-torsion group.    

\rem{The scheme of endomorphisms of the constant $\Z$-group scheme $\Z/N\Z$ (cf. Lemma \ref{hom-finite-free-is-representable}) is clearly the constant $\Z$-ring scheme attached to the ring $\Z/N\Z$. By Cartier duality (e.g. \cite[I.2.9]{Oort}), the obvious map $(\Z/N\Z)_{\Z} \rightarrow \underline{\mrm{End}}_{\Gp}(\mu_N)$ is also an isomorphism. It is also not difficult to see that evaluation at $1$ defines an isomorphism $\underline{\mrm{Hom}}_{\Gp}(\Z/N\Z,\mu_N) \rightarrow \mu_N$ of group schemes.}

\lem[morphism-schemes-muN-ZNZ]{Let $N \geq 1$. The scheme $\mu_N^{\vee} := \underline{\mrm{Hom}}_{\Gp}(\mu_N,\Z/N\Z)$ is an affine \'etale $\Z$-group scheme. It is finite above $\Z[1/N]$, and the unit section $\Sp{\Z} \rightarrow \mu_N^{\vee}$ is an open immersion.}

\demo{The group functor $\mu_N^{\vee}$ is well-defined and representable by an affine $\Z$-group scheme of finite presentation by Lemma \ref{hom-finite-free-is-representable}. Moreover, $\mu_N^{\vee}$ is formally \'etale by \cite[Exp. I, Th. 5.5]{SGA1}, so it is \'etale by \cite[$\text{IV}_4$ (17.6.1)]{EGA}. The unit section is then an open immersion by \cite[$\text{IV}_4$ (17.4.2)]{EGA}. Now, with $\OO_N=\Z[1/N,\zeta_N]$ (which is finite flat over $\Z[1/N]$), one has $(\mu_N)_{\OO_N} \simeq (\Z/N\Z)_{\OO_N}$, thus
\vspace{-7pt}\[(\mu_N)^{\vee}_{\OO_N} \simeq \underline{\mrm{End}}((\Z/N\Z)_{\OO_N},(\Z/N\Z)_{\OO_N}) \simeq (\Z/N\Z)_{\OO_N},\vspace{-7pt}\] so $(\mu_N)^{\vee}$ is finite above $\Z[1/N]$.  
}

\rem{Let $a,b \geq 1$ be coprime integers and $N=ab$. Then one has an isomorphism $\mu_N^{\vee} \simeq \mu_a^{\vee} \times \mu_b^{\vee}$ of group schemes.}

\rem[trivial-muvee]{Let $p$ be prime and $r \geq 1$. Then $(\mu_{p^r}^{\vee})_{\F_p}$ is the trivial $\F_p$-group scheme. Indeed, it is \'etale, and any morphism $f: (\mu_{p^r})_{\overline{\F_p}} \rightarrow (\Z/p^r\Z)_{\overline{\F_p}}$ is trivial, since its image is contained in the unit component of $(\Z/p^r\Z)_{\overline{\F_p}}$. Thus, for any $d \geq 1$, $(\mu_{p^r}^{\vee})_{\Z/p^d\Z}$ is trivial, since the unit section is a surjective open immersion.}

\lem[muN-muNvee-commutative]{Let $N \geq 1$. The two bilinear pairings $\mu_N \times \mu_N^{\vee} \rar \Z/N\Z$ defined below are equal: 
\vspace{-14pt}\begin{align*}
\mu_N \times \mu_N^{\vee} &\simeq \mu_N \times \underline{\mrm{Hom}}_{\Gp}(\mu_N,\Z/N\Z) \rightarrow \Z/N\Z,\\
\mu_N \times \mu_N^{\vee} &\simeq \underline{\mrm{Hom}}_{\Gp}(\Z/N\Z,\mu_N) \times \underline{\mrm{Hom}}_{\Gp}(\mu_N,\Z/N\Z) \rightarrow \underline{\mrm{End}}(\mu_N) \simeq \Z/N\Z
\end{align*}\vspace{-14pt}
}

\demo{Because $\mu_N,\mu_N^{\vee}$ are flat affine $\Z$-schemes of finite presentation with \'etale generic fibre, it is enough to show equality on $\C$-points. This translates to the following easy claim: let $G,H$ be two groups isomorphic to $\Z/N\Z$ and $u: G \rightarrow H$, $v: H \rightarrow G$ be two group homomorphisms. Then there exists an $n \in \Z$ such that $u \circ v$ and $v \circ u$ are exactly the multiplication by $n$.}

It is then formal to deduce:

\cor[endo-muN+ZNZ-ring]{The obvious morphism
\vspace{-7pt}\[\begin{pmatrix}\Z/N\Z & \mu_N \\\mu_N^{\vee} & \Z/N\Z\end{pmatrix} \rightarrow \underline{\mrm{End}}(\mu_N \times \Z/N\Z)\vspace{-7pt}\] of ring schemes (given by the obvious operations and Lemma \ref{muN-muNvee-commutative}) is an isomorphism. Moreover, the rule $\det: \begin{pmatrix} a & b\\c & d\end{pmatrix} \in M(S) \mapsto ad-bc \in (\Z/N\Z)(S)$ is multiplicative, and the isomorphism sends the closed open subscheme $\det^{-1}((\Z/N\Z)^{\times})$ to the subfunctor $\underline{\mrm{Aut}}_{\Gp}(\mu_N \times \Z/N\Z)$. 

In particular, $\underline{\mrm{End}}(\mu_N \times \Z/N\Z)$ and $\underline{\mathrm{Aut}}_{\Gp}(\mu_N \times \Z/N\Z)$ are affine quasi-finite flat $\Z$-schemes.}

\defi{Let $N \geq 1$ be an integer. A group scheme $G$ over a scheme $S$ is an \emph{ordinary $N$-torsion group} if it is finite locally free, commutative, killed by $N$, and such that, for every $s \in S$, $G_s$ becomes isomorphic to $\mu_N \times \Z/N\Z$ over an algebraic closure of $\kappa(s)$.}

\rem{Let $N \geq 1$ and $G$ be a finite locally free group scheme over a scheme $S$.
\begin{itemize}[noitemsep,label=\tiny$\bullet$]
\item Suppose that $G$ is commutative and killed by $N$. Then $G$ is an ordinary $N$-torsion group if and only if for every prime $p \mid N$, $G[p^{v_p(N)}]$ is an ordinary $p^{v_p(N)}$-torsion group.
\item Let $E/S$ be an elliptic curve such that $G \simeq E[N]$. By \cite[(2.8.5), (2.8.7.1), Theorem 2.9.3]{KM}, $G$ is ordinary if and only if every specialization of $E$ in characteristic dividing $N$ is ordinary. 
\item Assume that $G$ is commutative and killed by $N$. By Corollary \ref{iso-over-finite-extension}, if $S'$ is a surjective $S$-scheme and $G_{S'}$ is an ordinary $N$-torsion group, then $G$ is an ordinary $N$-torsion group. 
\end{itemize}}

The following result justifies the choice of name in the definition.  

\lem[ord-N-tors]{Let $G$ be an ordinary $N$-torsion group over a scheme $S$. Then $G$ is a $N$-torsion group.}

\demo{This follows from Corollary \ref{weN-descent} \ref{weNd-1b}.}

\prop[isomorphisms-to-ordinary]{Let $G$ be a finite locally free commutative group scheme killed by $N \geq 1$ over a scheme $S$. The scheme $I := \underline{\mathrm{Isom}}_{\Gp}(G,(\mu_N \times \Z/N\Z)_S)$ is flat over $S$. }

\demo{Working Zariski-locally over $S$, we may assume that $S=\Sp{R}$ for some ring $R$. By prime factorization, we may assume that $N=p^r$ for some prime $p$ and some $r \geq 1$. By approximation, we may assume that $S$ is affine and Noetherian. By descent, it is enough to prove the claim after making all base changes $\Sp{\widehat{\OO_{S,s}}} \rightarrow S$ for $s \in S$ in the image of the structure map $I \rightarrow S$. The image of $I \rightarrow S$ consists of those $s \in S$ such that $I_s(\overline{\kappa(s)})$ is not empty, i.e. such that $G_{\overline{\kappa(s)}} \simeq \mu_N \times \Z/N\Z$; therefore, we may assume that $(R,\mathfrak{m})$ is a complete Noetherian local ring with residue field $k$ and $G_{\overline{k}} \simeq \mu_N \times \Z/N\Z$. If $p \notin \mathfrak{m}$, then $G$ is finite \'etale (cf. e.g. \cite[(3.7), II]{TateGps}), so it is \'etale-locally constant, so $I$ is \'etale-locally constant, thus flat. Thus we may assume that $p \in \mathfrak{m}$. Then we have a connected-\'etale exact sequence (cf. e.g. \cite[(3.7), I]{TateGps})
\vspace{-3pt}\[0 \rightarrow C \rightarrow G \rightarrow D\rightarrow 0\vspace{-7pt}\]
where $C$ a finite flat connected $R$-group scheme and $D$ is a finite \'etale $R$-group scheme. The special fibre $C_k$ is geometrically connected by \cite[$\text{IV}_4$ (18.5.14)]{EGA} and \cite[Exp. $\text{VI}_{\text{A}}$, Prop. 2.4]{SGA3}. By uniqueness of the connected-\'etale exact sequence of $G_{\overline{k}}$, $C_{\overline{k}} \simeq \mu_N$ and $D_{\overline{k}} \simeq \Z/N\Z$. By \cite[Exp. X, Prop. 1.4, Cor. 3.8]{SGA3-2}, there is a finite \'etale faithfully flat $R$-algebra $R'$ such that $C_{R'} \simeq \mu_N$ and $D_{R'} \simeq \Z/N\Z$. Thus we have an exact sequence 
\vspace{-7pt}\[0 \rightarrow (\mu_N)_{R'} \rightarrow G_{R'} \overset{\pi}{\rightarrow} (\Z/N\Z)_{R'}\rightarrow 0\vspace{-7pt}\]
with $\pi$ finite faithfully flat by Lemma \ref{reform-FV}. Thus $\pi^{-1}(1)$ is a finite locally free surjective $R'$-scheme, and it admits a $R''$-point for some faithfully flat $R'$-algebra $R''$ (by \cite[$\text{IV}_4$ (17.16.2)]{EGA}). Hence the above exact sequence splits over $R''$, so $G_{R''} \simeq (\mu_N \times \Z/N\Z)_{R''}$, and $I \rightarrow \Sp{R}$ becomes flat after the faithfully flat base change $R \rightarrow R''$ by Corollary \ref{endo-muN+ZNZ-ring}. By descent, $I$ is flat over $R$ and we are done.  
}

\cor[ordinary-structure]{Let $G$ be an ordinary $N$-torsion group over a scheme $S$. Then $G$ is fppf-locally over $S$ isomorphic to $\mu_N \times \Z/N\Z$. Moreover, if $H$ is a finite locally free commutative $S$-group scheme killed by $N$, $I := \underline{\mrm{Isom}}_{\Gp}(G,H)$ is affine of finite presentation, flat and quasi-finite over $S$. The structure map $I \rightarrow S$ factors as $I \overset{\pi}{\rightarrow} S' \overset{\iota}{\rightarrow} S$, where $\pi$ is faithfully flat, $\iota$ is an open immersion, and $H_{S'}$ is an ordinary $N$-torsion group over $S'$.}

\demo{The first claim follows from Proposition \ref{isomorphisms-to-ordinary} by Corollary \ref{iso-over-finite-extension}. That $I \rightarrow S$ is affine of finite presentation is Lemma \ref{hom-finite-free-is-representable}; to show that it is quasi-finite flat, we may work fppf-locally over $S$ and thus assume that $G \simeq \mu_N \times \Z/N\Z$; the claim then follows from Proposition \ref{isomorphisms-to-ordinary}. 

The (set-theoretic) image $S'$ of $I \rightarrow S$ is open by \cite[$\text{IV}_2$ (2.4.6)]{EGA}. A point $s \in S$ is contained in $S'$ if and only if $I_s$ is not empty, if and only if $I_s(\overline{\kappa(s)})$ is non-empty, if and only if $H_{\overline{\kappa(s)}} \simeq \mu_N \times \Z/N\Z$, if and only if $H_s$ is an ordinary $N$-torsion group over $\kappa(s)$. Thus $H_{S'}$ is an ordinary $N$-torsion group.}

\cor[ordinary-locally-free]{Let $G,H$ be two ordinary $N$-torsion groups over a scheme $S$ and let $I$ be the $S$-scheme $\underline{\mathrm{Isom}}_{\Gp}(G,H)$. 
\begin{itemize}[noitemsep,label=\tiny$\bullet$]
\item If $N$ is invertible in $S$, $I$ is finite \'etale over $S$ of rank $|\GL{\Z/N\Z}|$. 
\item If $N=p^rN'$, where $p$ is prime and locally nilpotent in $S$, $r \geq 1$ and $p \nmid N'$, then $I$ is finite locally free over $S$ of rank $|\GL{\Z/N'\Z}|p^r\varphi(p^r)^2$, where $\varphi$ is Euler's totient function.   
\end{itemize}
}

\demo{By prime factorization, we may assume that $N$ is a power of the prime $p$. When $p$ is invertible in $S$, $G$ and $H$ are \'etale $N$-torsion groups, so $I$ is \'etale-locally isomorphic to $\underline{\mathrm{Aut}}((\Z/N\Z)^{\oplus 2})$, hence finite \'etale of rank $|\GL{\Z/N\Z}|$ by descent. When $p$ is locally nilpotent on $S$, working fppf-locally on $S$, we may assume that $G=H=\mu_N \times \Z/N\Z$. Then $I \simeq \begin{pmatrix} (\Z/N\Z)^{\times} & \mu_N \\ 0 & (\Z/N\Z)^{\times}\end{pmatrix}$ by Corollary \ref{endo-muN+ZNZ-ring} and Remark \ref{trivial-muvee}, and the conclusion follows.}

\cor[ordinary-locally-free-cond]{Let $G$ be an ordinary $N$-torsion group over a scheme $S$ and $H$ be a finite locally free commutative $S$-group scheme killed by $N$. Let $I := \underline{\mrm{Isom}}_{\Gp}(G,H)$. Then $I$ is finite over $S$ if and only the following two conditions hold:
\begin{itemize}[noitemsep,label=\tiny$\bullet$]
\item The image $S'$ of $I \rightarrow S$ (which is open in $S'$) is closed. 
\item The invertible locus of $N$ in $S'$ (which is open) is closed. 
\end{itemize}}

\demo{Recall that $S'$ is actually an open subscheme of $S$. If the invertible locus of $N$ in $S'$ is closed, then $S'$ can be covered with open subschemes $U$ on which either $N$ is invertible, or there is a prime $p \mid N$ nilpotent on $U$, and in particular $I \rightarrow S'$ is finite locally free by Corollary \ref{ordinary-locally-free}. If $S'$ is also closed in $S$, then $I \rightarrow S$ is also finite. 

Conversely, if $I \rightarrow S$ is finite, then, since it is flat of finite presentation, it is locally free, so the rank of $I \rightarrow S$ is a locally constant function on $S$. Since $S'$ is the locus where this rank is non-zero, it is open and closed in $S$. Moreover, the invertible locus of $N$ in $S'$ is the locus where this map has rank $|\GL{\Z/N\Z}|$ by Corollary \ref{ordinary-locally-free}, so it is also open and closed. } 

\cor[ordinary-moduli]{Let $G$ be an ordinary $N$-torsion group over a ring $R$. The moduli problem $\P_G(N)$ is affine of finite presentation and quasi-finite flat over $\Ell_R$. Moreover, for any object $E/S$ of $\Ell_R$, the $S$-scheme $\P_G(N)_{E/S}$ is finite if and only if the following closed subspaces of $S$ are open:
\begin{itemize}[noitemsep,label=\tiny$\bullet$]
\item the set of $s \in S$ such that $N$ vanishes in $\kappa(s)$,
\item the set of $s \in S$ such that $N$ vanishes in $\kappa(s)$ and $E_s$ is supersingular.
\end{itemize}}

\demo{The moduli problem $\P_G(N)$ is affine of finite presentation (by Proposition \ref{PG-representable-when-good}) and quasi-finite flat (by Corollary \ref{ordinary-structure}) over $\Ell_R$. If $E/S$ is an object of $\Ell_R$, the image of $\P_G(N)_{E/S} \rightarrow S$ is precisely the complement of the set of $s \in S$ such that $N$ vanishes in $\kappa(s)$ and $E_s$ is supersingular. The conclusion follows from Corollary \ref{ordinary-locally-free-cond}.  }

\rem{Let $R$ be a regular excellent domain of characteristic zero where $N$ is not invertible. It is tempting to try to use the above results to try to construct a compactification $X_G(N)$ of $Y_G(N)$ when $G$ is an ordinary $N$-torsion group over $R$. However, the moduli problem $\P_G(N)$ is not finite, so the methods of \cite[Chapter 8]{KM} do not obviously apply. This manifests in another way: the $j$-invariant $Y_G(N) \rar \bA^1_R$ is not finite, so the construction of Section \ref{situation-at-infty} does not obviously work. We could of course ignore this issue, since the $j$-invariant is at least quasi-finite, and simply consider the normalization $X_G(N)$ of the $j$-invariant $Y_G(N) \rar \bP^1_R$, but this only highlights the further difficulties that we are facing.

Assume for instance that $R$ is a Jacobson Dedekind domain (say, the $S$-integers in a number field). If $X_G(N)$ is a normal scheme (which seems reasonable, as the normalization of a morphism of smooth $R$-schemes), then it is Cohen--Macaulay and every nonempty open subscheme has dimension two, so that $X_G(N) \rar \mathbb{P}^1_R$ is finite and flat, hence locally free of constant rank. But the rank of the $j$-invariant $Y_G(N)_{R/\mfk{p}} \rar \mathbb{A}^1_{R/\mfk{p}}$, where $\mfk{p} \in \Sp{R}$, is not a constant function of $\mfk{p}$, as can be deduced from Corollaries \ref{ordinary-locally-free} and \ref{ordinary-locally-free-cond}. Therefore, the image of the map $Y_G(N) \rar X_G(N)$ does not remain dense after certain base changes $R \rar R/\mfk{p}$. This makes understanding the geometry of $X_G(N)$ more delicate, as there does not seem to be any clear moduli interpretation for the points of $X_G(N) \backslash Y_G(N)$ with finite ordinary $j$-invariant.

Another difficulty is that $X_G(N)$ contains $j$-invariants corresponding to supersingular elliptic curves in characteristic dividing $N$, which do not appear in $Y_G(N)$; understanding the normalization process at these points also seems nontrivial. }

\prop[ordinary-implies-relative-flat]{Let $G$ be an ordinary $N$-torsion group over a ring $R$. Let $m \geq 1$ divide $N$, so that $G[m]$ is an ordinary $m$-torsion group over $R$. For any object $E/S$ of $\Ell_R$, the forgetful map $\P_G(N)_{E/S} \rar \P_{G[m]}(m)_{E/S}$ sending an isomorphism $\iota: G_T \rar E_T[N]$ to $\iota[m]: G[m]_T \rar E_T[m]$ is affine and flat of finite presentation. The complement of its image consists exactly of the points above the $s \in S$ with residue characteristic $p \mid N$ not dividing $m$ such that $E_{\kappa(s)}$ is not ordinary. 
}

\demo{The first part of the claim follows from Corollary \ref{ordinary-structure}. Let $S$ be a $R$-scheme and $E/S$ be an elliptic curve. By Proposition \ref{PG-representable-when-good} and Corollary \ref{ordinary-moduli}, $\P_G(N)_{E/S}$ and $\P_{G[m]}(m)_{E/S}$ are affine of finite presentation, flat and quasi-finite over $S$. So all that remains is show that the forgetful map is flat and describe its image. By prime factorization, we may assume that $N=p^r$ for some prime $p$ and $m=p^s$ for some $0 \leq s < r$ (the cases $r=s$ and $r=0$ are clear). By the fibre criterion for flatness and descent, we may assume that $R=k$ is an algebraically closed field and $S=\Sp{k}$.

We treat the four cases separately: 
\begin{itemize}[noitemsep,label=$-$]
\item if $p$ is invertible in $k$, both schemes are constant, so the map is flat. That the map is surjective reduces to the well-known fact that $\GL{\Z/N\Z} \rightarrow \GL{\Z/m\Z}$ is surjective.    
\item if $k$ has characteristic $p$ and $E$ is supersingular, $\P_G(N)_{E/k}$ is empty, and $\P_{G[m]}(m)_{E/k}$ is empty if and only if $s=0$, so the map is always flat, and the claims about the image follow). 
\item if $k$ has characteristic $p$ and $E$ is ordinary, $E[N] \simeq \mu_N \times \Z/N\Z$, so, by Corollary \ref{endo-muN+ZNZ-ring} and Remark \ref{trivial-muvee}, this reduces to the fact that the morphism 
\vspace{-5pt}\[\begin{pmatrix}(\Z/N\Z)^{\times} & \!\!\!\!\mu_N \\ 0 & \!\!\!\!(\Z/N\Z)^{\times}\end{pmatrix} \rightarrow \begin{pmatrix} (\Z/m\Z)^{\times} & \!\!\!\!\mu_m \\ 0 & \!\!\!\!(\Z/m\Z)^{\times}\end{pmatrix}, \;\;\; \begin{pmatrix} a & \zeta \\ 0 & b\end{pmatrix} \mapsto \begin{pmatrix} a \!\!\!\!\mod{m} & \zeta^{N/m} \\ 0 & b\!\!\!\!\mod{m}\end{pmatrix}\vspace{-5pt}\] 
 is a flat surjective morphism of $k$-schemes, which is clear.    
\end{itemize}}

\prop[ord-supersing-automorphisms]{Let $E$ be a supersingular elliptic curve over an algebraically closed field $k$ of characteristic $p$. Let $r \geq 1$ be an integer and $F_r$ be the $k$-group scheme $\underline{\mrm{Aut}}_{\Gp}(E[p^r])$, which is affine of finite type. There is a closed immersion $\phi: \mathbb{G}_a \rightarrow F_r$ of group schemes such that the quotient scheme $\mathbb{G}_a \backslash F_r$ is finite over $k$. Moreover, for any $r \geq 1$, the forgetful morphism $F_{r+1} \rar F_r$ has finite image.} 

\demo{

Since $\bA^1_k$ represents the functor $S \mapsto \OO(S)$ (from the category of $k$-schemes to that of rings of characteristic $p$), it is endowed with a structure of ring scheme. The usual group scheme $\alpha_p=\Sp{k[x]/(x^p)}$ thus identifies with the kernel of the morphism of ring schemes $\bA^1_k \rightarrow \bA^1_k, \;z \in \OO(S) \mapsto z^p \in \OO(S)$: it is thus a $\bA^1_k$-module in the obvious sense, and one has a morphism of ring schemes $\bA^1_k \rightarrow \underline{\mathrm{End}}(\alpha_p)$, which one easily checks is a monomorphism\footnote{It is actually an isomorphism, but we do not need this.}. Recall the short exact sequence (cf. e.g. \cite[(15.5)]{Oort})
\vspace{-7pt}\[0 \rightarrow \alpha_p \rar E[p] \rar \alpha_p \rightarrow 0\vspace{-10pt}\] 

Denote the closed immersion $\alpha_p \rar E[p] \rar E[p^r]$ and the fppf morphism $E[p^r] \overset{p^{r-1}}{\rar} E[p] \rar \alpha_p$ by $\iota'$ and $\pi'$ respectively. In particular, one has $\pi' \circ \iota' = 0$. The rule 
\vspace{-7pt}\[u \in \mrm{End}((\alpha_p)_S) \mapsto \mrm{id}+\iota' \circ u \circ \pi' \in \mrm{Aut}(E[p^r]_S)\vspace{-7pt}\]
defines by composition a morphism of $k$-group schemes $\phi: \bA^1_k \rightarrow F_r$. Since $\iota'$ is a closed immersion and $\pi'$ is fppf, $\phi$ is a monomorphism of affine $k$-group schemes of finite type (by Lemma \ref{hom-finite-free-is-representable}), hence a closed immersion by \cite[Exp. $\text{VI}_{\text{B}}$, Cor. 1.4.2]{SGA3}. The quotient $X := \bA^1_k \backslash F_r$ exists and is a $k$-scheme of finite type by \cite[Exp. $\text{VI}_{\text{A}}$, Th. 3.2]{SGA3}. To check that $X$ is finite and thus that $\phi(\bA^1_k)$ is the unit component of $F_r$, it is enough to check that $F_r$ has dimension one by \cite[Exp. $\text{VI}_{\text{A}}$, Prop. 2.4]{SGA3}. 

The $p$-divisible group $D$ attached to $E$ has dimension and codimension $1$, and $E[p]$ contains a single copy of $\alpha_p$ because $\widehat{\OO_{E,0}}$ is a discrete valuation $k$-algebra with residue field $k$, so $E$ has exactly one closed subscheme of degree $p$ concentrated above the unit\footnote{This fact is well-known, but I am grateful to O. Gabber for bringing my attention to this argument.}. By \cite[Theorem 1]{Gabber-Vasiu}, the group scheme $F_r$ is thus one-dimensional\footnote{This also directly implies that $n_D=1$, although there is likely a more direct proof.}. The morphism $F_{r+1} \rar F_r$ has a finite image by \cite[Corollary 3(a)]{Gabber-Vasiu}, and we are done.}

\cor[supersingular-is-bad]{Let $G$ be a finite commutative group scheme over an algebraically closed field $k$ of characteristic $p$, and $\P$ denote the moduli problem $\P_G(p^r)$ for some $r \geq 1$. Let $E/k$ be an elliptic curve. Then $\P_{E/k}(k)$ is either empty, or
\begin{itemize}[noitemsep,label=\tiny$\bullet$]
\item finite if $E$ is ordinary,
\item infinite if $E$ is supersingular.
\end{itemize} }

\demo{By definition, $\P_{E/k}(k)$ is $\mrm{Isom}_{\GpSch_k}(G,E[p^r])$. If this set is not empty, then $G \simeq E[p^r]$, so $\P_{E/k}(k)$ is in bijection with $\mrm{Aut}_{\GpSch_k}(E[p^r])$. If $E$ is ordinary, $E[p^r]$ is ordinary and we apply Corollary \ref{ordinary-moduli}; if not, $E$ is supersingular and we apply Proposition \ref{ord-supersing-automorphisms}.}

\prop[quasifinite-ordinary-flat]{Let $E/R$ be an elliptic curve, and $N \geq 1$ be an integer. Let $\P$ be the moduli problem $\P_{E[N]}(N)$ over $R$. The following are equivalent:
\begin{enumerate}[noitemsep,label=(\roman*)]
\item \label{qfo-1} $\P$ is quasi-finite over $\Ell_R$,
\item \label{qfo-2} $\P$ is flat over $\Ell_R$,
\item \label{qfo-3} every specialization of $E$ in characteristic dividing $N$ is ordinary.
\end{enumerate}  
}

\demo{If \ref{qfo-3} holds, then $E[N]$ is an ordinary $N$-torsion group scheme over $R$, which implies \ref{qfo-1} and \ref{qfo-2} by Corollary \ref{ordinary-moduli}. It is enough to check that either \ref{qfo-1} or \ref{qfo-2} implies \ref{qfo-3} when $R$ is an algebraically closed field $k$ of characteristic $p \mid N$.

If $\P$ is quasi-finite, then $\P(E/k)\simeq \prod_{q \mid N}{\Aut{E[q^{v_q(N)}]}}$ is finite, so $E[p^{v_p(N)}]$ has finitely many automorphisms, so $E$ is ordinary by Proposition \ref{ord-supersing-automorphisms}. If $\P$ is flat, let $\ell \nmid 2N$ be a prime, then $j: \MPp{\P,[\Gamma(\ell)]_k} \rightarrow \bA^1_k$ is flat by Corollary \ref{over-Ell-concrete}, so its image is infinite. Thus there exists an ordinary elliptic curve $F/k$ such that $F[N] \simeq E[N]$. Hence $E[p] \simeq F[p]$ has a non-trivial $k$-point, thus $E$ is ordinary.}

\subsection{The curve $Y_{E[N]}(N)$ for a supersingular elliptic curve $E$}
\label{subsect:supersingular}
In this section, $k$ is a field of prime characteristic $p$, and $\P$ denotes a relatively representable moduli problem over $k$ which is affine \'etale over $\Ell_k$.  

\prop[nothing-on-mp]{Suppose that $\P$ is representable, let $\cE/\M$ be the universal elliptic curve and let $r \geq 1$. Let $G$ be a finite $k$-group scheme such that $G_{\overline{k}} \simeq F[p^r]$ for some supersingular elliptic curve $F/\overline{k}$. Then the map $f: \MPp{\P,\P_G(p^r)} \rar \M$ factors through a closed subscheme $Z \subset \M$ which is finite \'etale over $k$ and supported on the supersingular locus of $\M$.} 

\demo{By Corollary \ref{over-Ell-concrete}, $\M$ is smooth over $k$ with quasi-finite $j$-invariant, so by \cite[Lemma 12.5.4]{KM} its supersingular locus is a closed subscheme $Z$ which is finite \'etale over $k$. By Corollary \ref{general-smoothness-equivalence}, $\M_r := \MPp{\P,\P_G(p^r)}$ is smooth over $k$, hence geometrically reduced, so its scheme-theoretic image $Y$ under $f$ is geometrically reduced. So it is enough to show that $f(\M_r(\overline{k})) \subset Z(\overline{k})$.  
Now, let $P \in \M_r(\overline{k})$; then $P$ corresponds to an isomorphism $\iota: G_{\overline{k}} \rar E[p^r]$ for some elliptic curve $E/\overline{k}$. Then $E[p] \simeq G_{\overline{k}}[p] \simeq F[p]$ has no non-trivial $\overline{k}$-point, so $F$ is supersingular (by \cite[Proposition 12.3.6]{KM}), and we are done. 
}

\cor[structure-relative-mp]{Keep the assumptions and notations of Proposition \ref{nothing-on-mp}. The $\M$-scheme $\M_r := \MPp{\P,\P_G(p^r)}$ is the disjoint union over all the supersingular points $x \in \M$ of the schemes $\underline{\mrm{Isom}}_{\Gp}(G_{\kappa(x)},\cE_{\kappa(x)}[p^r]) \rar \Sp{\kappa(x)} \rar \M$. } 

\medskip

\demo{By Proposition \ref{product-rel-rep}, $\M_r$ identifies to the $\M$-scheme $\P_G(p^r)_{\cE/\M}$. If $Z \subset \M$ is the supersingular locus, one has by Proposition \ref{nothing-on-mp} 
\vspace{-10pt}\[\M_r = \M_r \times_{\M} Z \simeq \P_G(p^r)_{\cE_Z/Z} \simeq \coprod_{x \in Z}{\P_G(p^r)_{\cE_x/\kappa(x)}}\vspace{-11pt}\]
and the conclusion follows.}

\medskip

\cor[structure-forgetful-map-supersing]{Let $G,r$ be as in Proposition \ref{nothing-on-mp} with either $p^r > 2$ or $\P$ representable. The scheme-theoretic image of the forgetful map $\MPp{\P,\P_G(p^r)} \rightarrow \MPp{\P,\P_{G[p^{r-1}]}(p^{r-1})}$ is a finite \'etale $k$-scheme.}

\demo{This is a morphism of smooth affine $k$-schemes by Proposition \ref{general-smoothness-PG}, so its scheme-theoretic image $Z$ is affine and geometrically reduced, and we only need to check that $Z$ is finite. We can replace $\P$ by $\P \times [\Gamma(p+1)]_k$ and thus assume that $\P$ is representable. Then Corollary \ref{structure-relative-mp} and Proposition \ref{ord-supersing-automorphisms} imply that
\vspace{-4pt}\[\MPp{\P,\P_G(p^r)}(\overline{k}) \rightarrow \MPp{\P,\P_{G[p^{r-1}]}(p^{r-1})}(\overline{k})\vspace{-4pt}\]
has finite image.   }

\medskip

\cor[iso-scheme-is-smooth]{Assume that $k$ is algebraically closed. Let $E/k$ be a supersingular elliptic curve and $r \geq 1$ be an integer. Then the group scheme $F_r := \underline{\mrm{Aut}}_{\Gp}(E[p^r])$ is smooth over $k$ of relative dimension one. The closed immersion $\phi: \mathbb{G}_a \rightarrow F_r$ defined in Proposition \ref{ord-supersing-automorphisms} is an open immersion onto the connected component of the identity in $F_r$, so that the quotient $\mathbb{G}_a \backslash F_r$ is a finite constant $k$-scheme.}

\demo{Let $F_r^0 \subset F_r$ be the connected component of the identity. Then Corollaries \ref{general-smoothness-equivalence} and \ref{structure-relative-mp} imply that $\MPp{[\Gamma(p+1)]_k,\P_{E[p^r]}(p^r)}$ is a smooth $k$-scheme of relative dimension one which contains an open subscheme isomorphic to $F_r$. Hence $F_r$ is smooth over $k$ of relative dimension one, so $F_r^0$ is reduced and the surjective closed immersion $\phi: \mathbb{G}_a \rightarrow F_r^0$ is an isomorphism. The conclusion then follows from \cite[Exp. $\text{VI}_{\text{A}}$, 5.5]{SGA3}.   }

\rem{The given proofs of Corollaries \ref{structure-relative-mp} and \ref{iso-scheme-is-smooth} rely on Proposition \ref{general-smoothness-PG}, and thus on the theory of deformations of (truncated) Barsotti--Tate groups. They also admit more direct proofs, respectively by adapting the proof of \cite[Theorem 12.4.3]{KM} and by showing directly\footnote{This was pointed out to me by O. Gabber.} that $\dim{\mrm{Lie}(F_r)}=1$ (then using Proposition \ref{ord-supersing-automorphisms}). This implies that $\MPp{\P,\P_G(p^r)}$ is smooth over $k$ in the situation of Proposition \ref{nothing-on-mp} then, by a more careful study of the ordinary case, Proposition \ref{general-smoothness-PG} over fields. However, it is unclear whether one can handle more general bases with this approach. }

\prop[quotient-of-iso-scheme]{Assume that $k$ is algebraically closed and let $E/k$ be a supersingular elliptic curve. Then $\OO := \mrm{End}(E)$ is a maximal order in the unique quaternion algebra over $\Q$ ramified at $p$ and $\infty$, and $\widehat{\OO} := \mrm{End}(E) \otimes \Z_p$ is the maximal order in the division quaternion algebra over $\Q_p$, with unique two-sided maximal ideal $\mfk{P}$. Moreover, the obvious map $\widehat{\OO} \rar \mrm{End}(E[p^{\infty}])$ is an isomorphism. 

Let $r \geq 1$, $F_r = \underline{\mrm{Aut}}_{\Gp}(E[p^r])$, and $F_r^0$ be its unit connected component. For any $r \geq 1$, the natural map $(\OO/p^r\OO)^{\times} \simeq (\widehat{\OO}/p^r\widehat{\OO})^{\times} \rar F_r^0(k)\backslash F_r(k)$ is surjective with kernel $(1+p^{r-1}\mfk{P})$.}

\demo{The ring $\OO = \mrm{End}(E)$ is a maximal order in a definite quaternion algebra $B$ over $\Q$ (cf. \cite[Proposition 42.1.6, Theorem 42.1.9]{Voight}), and $\widehat{\OO}$ is thus maximal in $B \otimes \Q_p$. It is well-known\footnote{see e.g. \cite[(2.2)]{Tate-pdiv} and \cite[Proposition 1.7]{Drinfeld}} that one has abstractly $\mrm{End}(E[p^{\infty}]) \simeq \widehat{\OO}$, so the obvious injection $j: \widehat{\OO} \rightarrow \mrm{End}(E[p^{\infty}])$ makes $\mrm{End}(E[p^{\infty}])$ into a super-order of $\widehat{\OO}$: by e.g. \cite[Chapter 13]{Voight}, $j$ is an isomorphism. 

For the rest of the argument, since $j(E) \in \F_{p^2}$, we may assume that $k=\overline{\F}_p$, and we use Dieudonn\'e theory as in \cite{Manin}. Let $W$ be the discrete valuation ring of Witt vectors of $k$ with Frobenius automorphism $\varphi$, and let $D(-)$ denote the contravariant Dieudonn\'e module, which is a $W$-module carrying semi-linear actions of the Frobenius $F$ and the Verschiebung $V$. Then $D(E[p^r])$ (resp. $D(\alpha_p)$) is a free $(W/p^r)$-module of rank two (resp. a free $W/p$-module of rank one on which $F$ and $V$ act trivially, cf. \cite[Remark 3 on Theorem 2]{Oort-Tate}). 

Moreover, by \cite[II.4, IV.3 Theorem 4.1, III.5 Example 1]{Manin}, $D(E[p^{\infty}])$ is isomorphic to $W[F,V]/(FV-p,VF-p,F-V)W[F,V]$. Hence we can write $D(E[p^{\infty}]) = W \cdot (Fe_{\infty})\oplus W \cdot e_{\infty}$ with $Ve_{\infty}=Fe_{\infty}$. As a consequence, one has $D(E[p^r]) \simeq D(E[p^{\infty}])/(p^r)$, hence
\vspace{-4pt}\[D(E[p^r]) = (W/p^r) \cdot (Fe_r) \oplus (W/p^r) \cdot e_r,\qquad Fe_r=Ve_r,\vspace{-4pt}\] 
 
where the inclusion $E[p^r] \rightarrow E[p^s]$ for $\infty \geq s \geq r \geq 1$ induces the map $e_{s} \mapsto e_r$ on Dieudonn\'e modules.

For $1 \leq r \leq \infty$ and $u \in \mrm{End}(E[p^r])$, let $M_r(u) \in \mathcal{M}_2(W/(p^r))$ be the matrix of $D(u)$ in the basis $(Fe_r,e_r)$. Then $M_r: \mrm{End}(E[p^r]) \rightarrow \mathcal{M}_2(W/(p^r))$ is an injective anti-homomorphism of rings; for $1 \leq s \leq r$, one has $M_s(u[p^s]) = M_r(u) \mod{p^s}$; and the image of $M_r$ is exactly the collection $E_r$ of matrices in the basis $(Fe_r,e_r)$ of elements in $\mrm{End}_{W/(p^r)}(D(E[p^r]))$ that commute to $F$ and $V$. 
  
The quotient $W^{\varphi^2} \rightarrow (W/(p))^{\varphi^2}$ is surjective (since $(W/(p))^{\varphi^2} \simeq \F_{p^2}$ and $W^{\varphi^2} \supset W(\F_{p^2})$), so, if $a \in W$ is such that $p^{r-1} \mid \varphi^2(a)-a$, one has $a = a_0+p^{r-1}a_1$ with $a_0 \in \Z_{p^2}, a_1 \in W$. It follows by a direct calculation that one has  
\vspace{-8pt}\begin{align*}
E_r &= \left\{\begin{pmatrix}\varphi(b) & a\\p\varphi(a) & b\end{pmatrix} \in \mathcal{M}_2(W/p^r) \mid b \in \Z_{p^2}/(p^r), a \in (\Z_{p^2} + p^{r-1}W)/(p^r)\right\}, && r\geq 1,\\
E_{\infty} &= \left\{\begin{pmatrix}\varphi(b) & a\\p\varphi(a) & b\end{pmatrix} \in \mathcal{M}_2(W) \mid a,b \in \Z_{p^2}\right\}.
\end{align*}\vspace{-12pt}

Let us now use the notations of the proof of Proposition \ref{ord-supersing-automorphisms}. The closed immersion $\iota': \alpha_p \rightarrow E[p^r]$ induces a non-zero map on Dieudonn\'e modules; since $D(\alpha_p)$ is a $k$-line killed by $F$, $D(\iota'): D(E[p^r])/F\cdot D(E[p^r]) \rar D(\alpha_p)$ is an isomorphism. Similarly, $D(\pi'): D(\alpha_p) \rar D(E[p^r])$ is injective onto the $k$-line generated by $p^{r-1}Fe_r$. Recall that $\bA^1_k$ has a structure of ring scheme and acts on $\alpha_p$; for $z \in \bA^1_k(k) \simeq k$, $z$ acts by multiplication by $z$ on $T_0(\alpha_p)$, so $D(z)$ is the multiplication by $z$ by \cite[1.4.3.2 Theorem]{CCO}. Therefore, one has 
\vspace{-4pt}\[M_r(F_r^0(k)) = \left\{\begin{pmatrix} 1 & p^{r-1}x \\0 & 1\end{pmatrix} \mid x \in W/(p)\right\}=M_{\infty}(1+p^{r-1}\mfk{P}),\vspace{-6pt}\]
and the conclusion follows.
 
}

\bigskip

\cor[fibre-struct-supersing+repmoduli]{Assume that $\P$ is representable and that $k$ is algebraically closed. Let $r \geq 1$ and let $G$ be the $p^r$-torsion subgroup scheme of some supersingular elliptic curve over $k$. The connected components of $\M_r := \MPp{\P,\P_G(p^r)}$ are isomorphic to $\bA^1_k$. 

Let $(E_i)_{i \in I}$ be a set of representatives for supersingular elliptic curves over $k$ up to isomorphism. For each $i \in I$, pick an isomorphism $u_i: G \rar E_i[p^r]$, and let $\OO_{i}=\mrm{End}(E_i)$ and $\mfk{P}_i$ be its maximal two-sided ideal of residue characteristic $p$; for $\alpha_i \in \P(E_i/k)$ and $\beta_i \in (\OO_i/(p^r))^{\times}$, let $P_i(\alpha_i,\beta_i) \in \M_r(k)$ denote the point corresponding to the enriched elliptic curve $(E_i, \alpha_i, \beta_i \circ u_i)$. 

Any connected component of $\M_r$ contains some $P_i(\alpha_i,\beta_i)$. Moreover, the points $P_i(\alpha_i,\beta_i)$ and $P_j(\alpha'_j,\beta'_j)$ are in the same connected component (resp. equal) if and only if $i=j$ and, for some $t \in \OO_i^{\times}$, one has $(\P t)(\alpha'_j) =\alpha_j$ and $t \circ \beta_i\equiv \beta'_j \mod{p^{r-1}\mfk{P}_i}$ (resp. $t \circ \beta_i = \beta'_j$).}

\demo{For any supersingular elliptic curve $E/k$, one has $G \simeq E[p^r]$ by \cite[Theorem 2.9.3]{KM}. This implies the claim by Corollaries \ref{iso-scheme-is-smooth} and \ref{structure-relative-mp} and Proposition \ref{quotient-of-iso-scheme}.} 

\cor[fibre-struct-supersing+repmoduli2]{Keep the notations and assumptions of Corollary \ref{fibre-struct-supersing+repmoduli}. Let, for each $i \in I$, $A_i \subset \P(E_i/k)$ (resp. $B_i \subset (\OO_i/(p^r))^{\times}$) be a set of representatives for $\P(E_i/k)/\OO_i^{\times}$ (resp. $\OO_i^{\times} \backslash (\OO_i/(p^{r-1}\mfk{P}_i))^{\times}$). Then the maps 
\vspace{-5pt}\[\coprod_{i \in I}{\P(E_i/k) \times B_i} \rightarrow \pi_0(\M_r),\qquad \coprod_{i \in I}{A_i \times (\OO_i/(p^r))^{\times}} \rightarrow \pi_0(\M_r)\vspace{-5pt}\]
sending $(i,\alpha_i,\beta_i)$ to the connected component of $P_i(\alpha_i,\beta_i)$ are bijective. 
}

\cor[fibre-struct-supersing+relrepmoduli]{Keep the notations of Corollary \ref{fibre-struct-supersing+repmoduli2}; but instead of assuming that $\P$ is representable, assume that either $p^r > 2$ or $\P$ is representable. Then the connected components of $\M_r$ are isomorphic to $\bA^1_k$ and the map  
\vspace{-3pt}\[\coprod_{i \in I}{\P(E_i/k) \times B_i} \rightarrow \pi_0(\M_r),\qquad (i,\alpha_i,\beta_i) \mapsto [P_i(\alpha_i,\beta_i)] \vspace{-5pt}\]
is a bijection. 
}

\demo{We already proved this if $\P$ is representable, so we drop the assumption that $\P$ is representable and assume instead that $p^r > 2$. All objects in Corollary \ref{fibre-struct-supersing+repmoduli} remain well-defined. Note that $\M_r$ (resp. $\pi_0(\M_r)$) then identifies with the quotient of the smooth affine $k$-scheme $\M'_r := \MPp{\P \times [\Gamma(p+1)]_k, \P_G(p^r)}$ (resp. of $\pi_0(\M'_r)$) by $\GL{\Z/(p+1)\Z}$. Corollary \ref{fibre-struct-supersing+repmoduli2} then implies that $\GL{\Z/(p+1)\Z}$ acts freely on $\pi_0(\M'_r)$, so, if $U$ is the union of a set of representatives in $\pi_0(\M'_r)$ for $\pi_0(\M'_r)/\GL{\Z/(p+1)\Z}$, the composition $U \subset \M'_r \rightarrow \M_r$ is an isomorphism. This implies the conclusion. }

\medskip

\cor[forgetful-supersing]{Let $r \geq 2$, and assume that $\P$ is representable if $p^r=4$. Let $G$ be a finite $k$-group scheme such that $G_{\overline{k}} \simeq F[p^r]$ for some supersingular elliptic curve $F/\overline{k}$. Then the forgetful map $\MPp{\P,\P_G(p^r)} \rar \MPp{\P,\P_G(p^{r-1})}$ (which has finite set-theoretic image by Corollary \ref{structure-forgetful-map-supersing}) induces a surjection $\pi_0(\MPp{\P,\P_G(p^r)}_{\overline{k}}) \rar \pi_0(\MPp{\P,\P_G(p^{r-1})}_{\overline{k}})$.}

\bigskip

\rem{Let $W$ be a complete regular local ring of characteristic zero with algebraically closed residue field $k$ of characteristic $p > 2$ and let $\mathcal{E}/W$ be an elliptic curve with supersingular special fibre. It is natural to ask how to define a reasonable compactification $X_{\cE[p]}(p)$ of $Y_{\cE[p]}(p)$ over $W$. The construction of \cite[Chapter 8]{KM}, consisting in taking the normalization $Y^{\nu}$ of $Y_{\cE[p]}(p) \rightarrow \bP^1_W$, is unsatisfactory for the following reason: $Y^{\nu} \rightarrow \bP^1_W$ is integral, so it has discrete fibres. Since the $j$-invariant $Y_{\cE[p]}(p)_k \rar \bA^1_k$ is locally constant, $Y_{\cE[p]}(p)_k \rightarrow Y^{\nu}_k$ is (at least topologically) locally constant, so $Y^{\nu}$ cannot be called a compactification of $Y_{\cE[p]}(p)$.   

Suppose that some suitable compactification $X_{\cE[p]}(p)$ exists (for instance, assume that $X_{\cE[p]}(p)$ contains a dense open subscheme $Y'$ finite flat over $Y_{\cE[p]}(p)$), the $j$-invariant $X_{\cE[p]}(p) \rar \bP^1_W$ is between proper $W$-schemes, so it has closed image; said image contains the generic point of $\bP^1_W$ so $X_{\cE[p]}(p)_k \rar \bP^1_k$ is proper surjective between $k$-schemes of pure dimension one. Now, $X_{\cE[p]}(p)_k$ contains as an open subscheme $Y'_k$, whose (set-theoretic) image under the $j$-invariant is finite. Thus, $X_{\mathcal{E}[p]}(p)_k$ has to contain irreducible components of dimension one that project onto $\mathbb{P}^1_k$, and thus only meet $Y'_k$ at finitely many points. 

There does not seem to be a clear moduli interpretation in the context of isomorphisms of torsion subschemes of elliptic curves for such hypothetical components. However, one possibility suggested by the approach of the book \cite{KM} is to replace the somewhat na\"ive notion of isomorphism of group schemes by something else, possibly closer to a Drinfeld level structure. }

\section{Galois twists of modular curves}
\label{moduli-is-twisting}
In this section, we fix an integer $N \geq 3$ and a field $k$ of characteristic not dividing $N$ with separable closure $k_s$. We also fix a left representation $\rho: \mrm{Gal}(k_s/k) \rar \GL{\Z/N\Z}$.

\subsection{$Y_{\rho}(N)$ and its rational points}

\defi{The group $\GL{\Z/N\Z}$ acts on the following diagram of quasi-projective $k$-schemes (and trivially on the second row)
\vspace{-4pt}\[
\begin{tikzcd}[ampersand replacement=\&]
Y(N)_k \arrow{r}{\mrm{in}}\arrow{d}{j}\& X(N)_k \arrow{d}{\overline{j}}\& C_N \arrow{l}{\mrm{cl}} \arrow{d}\\
\mathbb{A}^1_k \arrow{r}{\mrm{in}} \& \mathbb{P}^1_k \& \Sp{k} \arrow{l}{\infty}
\end{tikzcd}
\vspace{-4pt}\] where $C_N \subset X(N)$ is the cuspidal subscheme. Twisting by $\rho$ (in the sense of Propositions \ref{cocycle-twist} and \ref{twist-equiv}) thus gives a diagram of quasi-projective $k$-schemes:
\vspace{-4pt}\[
\begin{tikzcd}[ampersand replacement=\&]
Y_{\rho}(N)_k \arrow{r}{\mrm{in}}\arrow{d}{j}\& X_{\rho}(N) \arrow{d}{\overline{j}}\& C_{N,\rho} \arrow{l}{\mrm{cl}} \arrow{d}\\
\mathbb{A}^1_k \arrow{r}{\mrm{in}} \& \mathbb{P}^1_k \& \Sp{k} \arrow{l}{\infty}
\end{tikzcd}
\vspace{-4pt}\]
Moreover, the first cell is a Cartesian diagram, $j,\overline{j}$ are finite locally free of constant rank, and $C_{N,\rho}$ is finite \'etale over $k$. Moreover, $X_{\rho}(N)$ is smooth over $\Sp{k}$.}

\demo{This follows from Propositions \ref{cocycle-twist} and \ref{twist-equiv}, the properties of $X(N)$ discussed in the previous section, and descent (as in Proposition \ref{fpqc-descent-prop}).}

\lem[right-representation]{Let $G$ be an (\'etale) $N$-torsion group over $k$ (see Definition \ref{torsion-group}). Let $(A,B)$ be a $\Z/N\Z$-basis of $G(k_s)$. For every $\sigma \in \mrm{Gal}(k_s/k)$, let $D_G(\sigma) \in \GL{\Z/N\Z}$ be such that $\begin{pmatrix}\sigma(A)\\\sigma(B)\end{pmatrix} = D_G(\sigma)\begin{pmatrix}A\\B\end{pmatrix}$. Then, for $\sigma,\sigma' \in \mrm{Gal}(k_s/k)$, one has $D(\sigma\sigma')=D(\sigma')D(\sigma)$.}

\demo{Simple computation.}

\bigskip

From this result, we obtain a description of the rational points of $Y_{\rho}(N)$. 

\prop[rational-Yrho]{For any subextension $L$ of $k_s/k$, $Y_{\rho}(N)(L)$ is in natural bijection with the $L$-isomorphism classes of triples $(E,P,Q)$ satisfying the following properties: 
\begin{itemize}[noitemsep,label=\tiny$\bullet$]
\item $E$ is an elliptic curve defined over $L$,
\item $(P,Q)$ is a $\Z/N\Z$-basis of $E[N](k_s)$,
\item For all $\sigma \in \mrm{Gal}(k_s/L)$, $\begin{pmatrix}\sigma(P)\\\sigma(Q)\end{pmatrix} = \rho^{-1}(\sigma)\begin{pmatrix}P \\ Q\end{pmatrix}$.
\end{itemize}
If $x \in Y_{\rho}(N)(L)$ is given by the triple $(E,P,Q)$, then $j(x) \in \bA^1(L)$ is the $j$-invariant of $E$.}

\demo{Let $k \subset L \subset k_s$ be a subextension. Then $Y_{\rho}(N)(L)$ is in bijection with the triples $(E,P,Q)$, where $E$ is an elliptic curve over $k_s$, and $(P,Q)$ is a basis of $E[N](k_s)$ such that, for every $\sigma \in \mrm{Gal}(k_s/L)$, the triples $(E,P,Q)$ and $(\sigma(E),a\sigma(P)+b\sigma(Q),c\sigma(P)+d\sigma(Q))$, where $\rho(\sigma)=\begin{pmatrix}a & b\\c & d\end{pmatrix}$, are isomorphic over $k_s$ (by the representability of $[\Gamma(N)]_k$). In particular, $j(E) \in k_s$ is fixed by $\mrm{Gal}(k_s/L)$, hence we can assume that $E$ is defined over $L$, and that, for all $\sigma \in \mrm{Gal}(k_s/L)$ with $\rho(\sigma)=\begin{pmatrix}a & b\\c & d\end{pmatrix}$, there is a $k_s$-automorphism $t_{\sigma}$ of $E_{k_s}$ sending $(P,Q)$ to $(a\sigma(P)+b\sigma(Q),c\sigma(P)+d\sigma(Q))$.  

We claim that $\sigma \longmapsto t_{\sigma}^{-1}$ satisfies the cocycle property of Proposition \ref{cocycle-twist}. Since the identity to check is the equality of two $k_s$-automorphisms of $E_{k_s}$, it is enough, by rigidity, to check that it holds on $E[N](k_s)$, and hence on $P$ and $Q$. For $\sigma,\sigma' \in \mrm{Gal}(k_s/L)$, a direct calculation shows that $t_{\sigma\sigma'}\left(\begin{pmatrix}P \\ Q\end{pmatrix}\right) = \sigma t_{\sigma'}\sigma^{-1}t_{\sigma}\left(\begin{pmatrix}P \\ Q\end{pmatrix}\right)$,
thus the cocycle property holds.

Let $E'=E_{t^{-1}}$ be the twisted elliptic curve over $L$:  
we still have $P,Q \in E'[N](k_s)$, $(E',P,Q)$ corresponds to the same point of $Y(N)(k_s)$ as $(E,P,Q)$, and another direct calculation shows that for $\sigma \in \mrm{Gal}(k_s/L)$, one has 
\vspace{-7pt}\[\begin{pmatrix}\sigma_{E'}(P)\\\sigma_{E'}(Q)\end{pmatrix}=\rho(\sigma)^{-1}\begin{pmatrix}P \\ Q\end{pmatrix}.\vspace{-7pt}\]

In other words, we proved that every point in $Y_{\rho}(N)(L)$ came from a triple satisfying the required conditions. Conversely, every such triple clearly produces a $k_s$-rational point on $Y_{\rho}(N)$ stable under $\mrm{Gal}(k_s/L)$, i.e. a $L$-point of $Y_{\rho}(N)$. 

To conclude, all we need to do is show that any two triples $(E,P,Q),(E',P',Q')$ satisfying the requested conditions and representing the same point in $Y(N)(k_s)$ are $L$-isomorphic. Indeed, these triples are isomorphic over $k_s$ (by the representability of $[\Gamma(N)]$), so let $f: E_{k_s} \rar E'_{k_s}$ be a $k_s$-isomorphism mapping $(P,Q)$ to $(P',Q')$. The induced map $f[N]: E[N](k_s) \rar E'[N](k_s)$ commutes to $\mrm{Gal}(k_s/L)$, hence by rigidity the same holds for $f$, so that, by Corollary \ref{inf-galois-desc}, $f$ is the base change of an $L$-isomorphism $E \rightarrow E'$. 

The last statement is clear from the construction and the definition of $Y_{\rho}$. 
}

\prop[quadratic-twist-equiv]{Let $\chi: \mrm{Gal}(k_s/k) \rar \{\pm 1\}$ be a character, and let \[\rho'=\chi \otimes \rho: \mrm{Gal}(k_s/k) \rar \GL{\Z/N\Z}.\] The natural $k_s$-isomorphism between the following two diagrams
\[
\begin{tikzcd}[ampersand replacement=\&]
Y_{\rho}(N)_k \arrow{r}{\mrm{in}}\arrow{d}{j}\& X_{\rho}(N)_k \arrow{d}{\overline{j}}\& (C_N)_{\rho} \arrow{l}{\mrm{cl}} \arrow{d} \& Y_{\rho'}(N)_k \arrow{r}{\mrm{in}}\arrow{d}{j}\& X_{\rho'}(N)_k \arrow{d}{\overline{j}}\& (C_N)_{\rho'} \arrow{l}{\mrm{cl}}\arrow{d} \\ 
\mathbb{A}^1_k \arrow{r}{\mrm{in}} \& \mathbb{P}^1_k \& \Sp{k} \arrow{l}{\infty} \& \mathbb{A}^1_k \arrow{r}{\mrm{in}} \& \mathbb{P}^1_k \& \Sp{k} \arrow{l}{\infty}
\end{tikzcd}
\]
coming from the definition of a Galois twist is defined over $k$.

Moreover, for any finite subextension $L \subset k_s$ of $k$, the induced map $Y_{\rho}(N)(L) \rar Y_{\rho'}(N)(L)$ is given, through the description of Proposition \ref{rational-Yrho}, by the rule $(E,P,Q) \longmapsto (E \otimes \chi, P, Q)$, where $E \otimes \chi$ is the quadratic twist of $E$ by $\chi_{|\mrm{Gal}(k_s/L)}$.}

\demo{By Corollary \ref{inf-galois-desc} (our schemes are of finite type over $\Sp{k}$), the isomorphism is defined over $k$ as long as it commutes to $\mrm{Gal}(k_s/k)$. Since $X_{\rho}(N),X_{\rho'}(N)$ are proper smooth over $k$, it is enough to check that the induced bijection $\tau: Y_{\rho}(N)(k_s) \rar Y_{\rho'}(N)(k_s)$ commutes with the action of $\mrm{Gal}(k_s/k)$. 

If $E/k_s$ is an elliptic curve and $(P,Q)$ is a $\Z/N\Z$-basis of $E[N](k_s)$, let $x \in Y_{\rho}(N)(k_s)$ and $x' =\tau(x)\in Y_{\rho'}(N)(k_s)$ be the geometric points attached to $(E,P,Q)$. A direct calculation shows that the automorphism $\chi(\sigma)$ of $\sigma(E)$ identifies the triple corresponding to $\tau(\sigma(x))$ to $\sigma(x')$.

Let now $k \subset L\subset k_s$ be an intermediate extension and $(E,P,Q)$ be a triple as in Proposition \ref{rational-Yrho} defining some $L$-point of $Y_{\rho}(N)$. Then $(E \otimes \chi,P,Q)$ is a triple as in \emph{loc.cit.} parametrizing a $L$-point of $Y_{\rho'}(N)$. Since the triples $(E \otimes \chi,P,Q)$ and $(E,P,Q)$ are $k_s$-isomorphic, \emph{loc.cit.} implies that $\tau((E,P,Q)) = (E \otimes \chi,P,Q)$ in $Y_{\rho'}(N)(L)$. 
 
}

\subsection{Comparison with $X_G(N)$}

The following result is well-known (e.g. it follows from \cite[Exp. V, Th. 4.1, \S 7, Prop. 8.1]{SGA1}). 

\prop[tannakian-group]{
The functor $G \longmapsto G(k_s)$, from finite \'etale commutative group schemes over $k$ to finite discrete $\mrm{Gal}(k_s/k)$-modules, is an equivalence of categories.}

\nott{From now on, $m \geq 1$ denotes an integer coprime to $N$ such that, if $k$ has prime characteristic $p$, one has $p^2 \nmid m$. }

\prop[xrho-elementary]{Since $\GL{\Z/N\Z}$ acts on the inclusion $Y(N,\Gamma_0(m)) \subset X(N,\Gamma_0(m))$ of $\bP^1_k$-schemes (through the $j$-invariant), we can define the twist $Y_{\rho}(N,\Gamma_0(m)) \subset X_{\rho}(N,\Gamma_0(m))$ of $Y(N,\Gamma_0(m))_k \subset X(N,\Gamma_0(m))_k$ by $\rho$. The $j$-invariant can be twisted to a finite locally free morphism $j_{\rho}: X_{\rho}(N,\Gamma_0(m)) \rightarrow \bP^1_k$ such that $j_{\rho}^{-1}(\bA^1_k) = Y_{\rho}(N,\Gamma_0(m))$. 

The $k$-scheme $X_{\rho}(N,\Gamma_0(m))$ is smooth, except, when $m$ vanishes in $k$, at the points where $j_{\rho}$ is a supersingular value. For every subextension $L/k$ of $k_s/k$, the points $x \in Y_{\rho}(N,\Gamma_0(m))(L)$ identify with the $L$-isomorphism classes of triples $(E,(P,Q),C)$ such that:
\begin{itemize}[noitemsep,label=$-$]
\item $E/L$ is an elliptic curve and $(E,P,Q)$ is a triple over $L$ as in Proposition \ref{rational-Yrho}, 
\item $C$ is a cyclic subgroup of $E$ of degree $m$,
\end{itemize}}

\demo{Since the $X(N,\Gamma_0(m))_k$ are finite over $\bP^1_k$, they are projective. They have an action of $\GL{\Z/N\Z}$ which commutes to the degeneracy maps and Atkin--Lehner involutions. Thus, the existence of the twisted modular curves, their regularity properties, and the fact that the $j$-invariant can be twisted all follow from Propositions \ref{cocycle-twist} and \ref{twist-equiv} and the previous results on the $X(N,\Gamma_0(m))_k$. 

To prove that we identified correctly the rational points of $X_{\rho}(N,\Gamma_0(m))$ (as well as the claims about the behavior of the degeneracy maps and the Atkin--Lehner automorphisms), it is enough to show that every $(E,P,Q,C) \in Y_{\rho}(N,\Gamma_0(m))(k_s)$ fixed by $\mrm{Gal}(k_s/L)$ is $k_s$-isomorphic to a unique quadruple $(E',P',Q',C')$ (up to $L$-isomorphism) as in the statement of the Proposition.

Note that if two such triples $(E,(P,Q),C)$ and $(E',(P',Q'),C)$ are $k_s$-isomorphic, the isomorphism is unique (since its action on the $N$-torsion is prescribed). 
 To this end, we argue as in Proposition \ref{rational-Yrho}. In particular, if $(E,(P,Q),C)$ and $(E',(P',Q'),C)$ are over $L$, the $k_s$-isomorphism is stable under $\mrm{Gal}(k_s/L)$ hence is a $L$-isomorphism. 
 
Conversely, if $(E,(P,Q),C) \in Y_{\rho}(N,\Gamma_0(m))(k_s)^{\mrm{Gal}(k_s/L)}$, then there is a triple $(E_0,P_0,Q_0)$ over $L$ as in Proposition \ref{rational-Yrho} and a $k_s$-isomorphism $\iota: (E_0,P_0,Q_0) \rightarrow (E,P,Q)$. Then $C_0 := \iota^{-1}(C) \in [\Gamma_0(m)]((E_0)_{k_s}/k_s)$ is actually fixed by $\mrm{Gal}(k_s/L)$, so $C_0 \in [\Gamma_0(m)](E_0/L)$ and $(E,(P,Q),C)$ is represented by the triples $(E_0,(P_0,Q_0),C_0)$, which satisfies all the required properties.}

\cor[xrho-degn-al]{In the situation of Proposition \ref{xrho-elementary}, the degeneracy maps and Atkin--Lehner automorphisms between the various $X(N,\Gamma_0(m))$ (defined in Proposition \ref{ext-degn}) can be twisted to finite locally free morphisms between the twists $X_{\rho}(N,\Gamma_0(m))$. The relations proved throughout Section \ref{gamma0m-structure} hold\footnote{Specifically, this refers to the third item of Lemma \ref{elem-degeneracies}, the first item of Lemma \ref{AL-elem}, and Proposition \ref{AL-auto}, after, in all instances, replacing the moduli scheme for $\P=\P(N)$ by the twist of its compactification.}.}

\demo{The Atkin--Lehner automorphisms and the degeneracy maps between the $X(N,\Gamma_0(m))$ are finite locally free maps and commute to $\GL{\Z/N\Z}$, so they can be twisted to finite locally free morphisms between the $X_{\rho}(N,\Gamma_0(m))$ by Propositions \ref{cocycle-twist} and \ref{twist-equiv}. These results also imply (by Proposition \ref{ext-degn}) that the aforementioned relations still hold for the twists of the compactified moduli schemes.}

\rem{Let $L/k$ be a subextension of $k_s/k$ and $d,t \geq 1$ be such that $dt\mid m$. The twist $(D_{d,t})_{\rho}: X_{\rho}(N,\Gamma_0(m)) \rightarrow X_{\rho}(N,\Gamma_0(t))$ sends $(E,(P,Q),C) \in X_{\rho}(N,\Gamma_0(m))(L)$ (in the description of Proposition \ref{xrho-elementary}) to $(E/C\{d\}, (P,Q) \mod{C\{d\}}, (C/C\{d\})\{t\})$. 

Similarly, the construction of Proposition \ref{AL-auto} also applies to describe the twists of the Atkin--Lehner automorphism.}

\prop[group-is-twist-unpolarized]{Let $G$ be a $N$-torsion group over $k$. Choose a basis $(A,B)$ of $G(k_s)$ and define $D_G$ as in Lemma \ref{right-representation}, so that $\sigma \in \mrm{Gal}(k_s/k) \longmapsto D_G(\sigma)^{-1} \in \GL{\Z/N\Z}$ is a group homomorphism. There is an isomorphism $\iota_m: X_G(N,\Gamma_0(m)) \rar X_{D_G^{-1}}(N,\Gamma_0(m))$ of projective $k$-schemes characterized by the following property: let $L/k$ be a subextension of $k_s/k$ and $x = (E/L,\varphi,C) \in Y_G(N,\Gamma_0(m))(L)$, then $\iota_m(x) \in X_{D_G^{-1}}(N,\Gamma_0(m))(L)$ is given in the description of Proposition \ref{xrho-elementary} by the triple $(E,(\varphi(A),\varphi(B)),C)$. 
In particular, $\iota_m$ preserves the $j$-invariant. More generally, the family $(\iota_m)_m$ commutes with the Atkin--Lehner automorphisms and all the degeneracy maps. 
}

\demo{Since $(A,B)$ is a $\Z/N\Z$-basis of $G(k_s)$, the rule  
\vspace{-8pt}\[\varphi \in \P_G(N)_{k_s}(F/S) \mapsto (\varphi(A),\varphi(B)) \in [\Gamma(N)]_{k_s}(F/S)\vspace{-8pt}\]
is an isomorphism of moduli problems of level $N$ over $k_s$. Hence it induces isomorphisms $\tau_m: Y_G(N,\Gamma_0(m))_{k_s} \rar Y(N,\Gamma_0(m))_{k_s}$, which extend to isomorphisms $\overline{\tau}_m$ between the compactified moduli schemes by Proposition \ref{compactification-functor}. The $\overline{\tau}_m$ are compatible with the Atkin--Lehner automorphisms and the degeneracy maps by Remark \ref{degn-AL-infty-natural}.  
 
By definition of the Galois twist, it is enough to show that, for any $\sigma \in \mrm{Gal}(k_s/k)$, the following diagram commutes:
\vspace{-8pt}\[
\begin{tikzcd}[ampersand replacement=\&]
X_G(N,\Gamma_0(m)) \times_k \Sp{k_s} \arrow{r}{\overline{\tau}_m}\arrow{d}{(\mrm{id},\underline{\sigma}^{-1})} \& X(N,\Gamma_0(m)) \times_k \Sp{k_s} \arrow{d}{D_G^{-1}(\sigma)\circ (\mrm{id},\underline{\sigma}^{-1})}\\
X_G(N,\Gamma_0(m)) \times_k \Sp{k_s} \arrow{r}{\overline{\tau}_m} \&  X(N,\Gamma_0(m)) \times_k \Sp{k_s}
\end{tikzcd}
\vspace{-10pt}\]

Since $Y_G(N,\Gamma_0(m))_{k_s} \rightarrow X_G(N,\Gamma_0(m))_{k_s}$ is scheme-theoretically dense and $X(N,\Gamma_0(m))_{k_s}$ is separated, it is enough to check that the diagram commutes when restricted to $Y_G(N,\Gamma_0(m))$. 

In other words, given $k$-morphisms $f: S \rar Y_G(N,\Gamma_0(m))$, $g: S \rar \Sp{k_s}$, we need to prove that $D_G^{-1}(\sigma) \circ \mrm{pr}_1 \circ \tau_m(f,g) = \tau_m(f,\underline{\sigma}^{-1} \circ g)$. This is a direct calculation using the descriptions of $\tau_m$ and $D_G$. }

\rem{In the situation of Lemma \ref{right-representation}, consider an elliptic curve $E/k$. Then, by Proposition \ref{tannakian-group}, for any subextension $k \subset L \subset k_s$,
\vspace{-10pt}\[\varphi \in \mrm{Isom}_{\GpSch_L}(G_L,E[N]_L) \mapsto (\varphi(A),\varphi(B)) \in E[N](k_s)^2\vspace{-9pt}\]
is injective onto the collection of bases $(P,Q)$ of $E[N](k_s)$ such that for every $\sigma \in \mrm{Gal}(k_s/L)$, one has $D_G(\sigma)\begin{pmatrix}P\\Q\end{pmatrix}=\begin{pmatrix}\sigma(P)\\\sigma(Q)\end{pmatrix}$. Furthermore, these maps respect the natural actions of $\mrm{Gal}(k_s/k)$ and $\mrm{Aut}(E_{k_s})$. After quotienting by $\mrm{Aut}(E_{k_s})$, we obtain the isomorphism induced by $\iota_m$ on $L$-rational points with $j$-invariant $j(E)$. }

\prop[group-is-twist-withpol]{Keep the notations of Proposition \ref{group-is-twist-unpolarized}. The twist $(\mu_N^{\times})_{(\det{D_G})^{-1}}$ of $(\mu_N^{\times})_k$ by the group homomorphism $(\det{D_G})^{-1}: \mrm{Gal}(k_s/k) \rar (\Z/N\Z)^{\times}$ is an \'etale $(\Z/N\Z)^{\times}$-torsor over $k$. Moreover, there is a unique isomorphism $\eta: \mrm{Pol}(G) \rar (\mu_N^{\times})_{(\det{D_G})^{-1}}$ of $(\Z/N\Z)^{\times}$-torsors over $k$ such that, for any $\alpha \in \mrm{Pol}(G)(k_s)$, one has 
\vspace{-6pt}\[\eta(\alpha)=\alpha(A,B) \in \mu_N^{\times}(k_s) \simeq (\mu_N^{\times})_{(\det{D_G})^{-1}}(k_s).\vspace{-6pt}\]
 Furthermore,
\begin{enumerate}[noitemsep,label=(\roman*)]
\item\label{gtw-1} the following diagram commutes:
\vspace{-8pt}\[
\begin{tikzcd}[ampersand replacement=\&]
X_G(N,\Gamma_0(m)) \arrow{r}{\overline{\tau}_m}\arrow{d}{\det} \& (X(N,\Gamma_0(m))_k)_{D_G^{-1}} \arrow{d}{(\mrm{We})_{D_G^{-1}}}\\
\mrm{Pol}(G) \arrow{r}{\eta} \&  (\mu_N^{\times})_{(\det{D_G})^{-1}}
\end{tikzcd}
\vspace{-10pt}\]
\item\label{gtw-2} let $d,t\geq 1$ be divisors of $m$ and $D_{d,t}$ be the degeneracy map from Proposition \ref{ext-degn}. The following diagram commutes:
\vspace{-8pt}\[
\begin{tikzcd}[ampersand replacement=\&]
X_G(N,\Gamma_0(t)) \arrow{ddd}{\det} \arrow{rrr}{\overline{\tau}_t} \& \& \& (X(N,\Gamma_0(t))_k)_{D_G^{-1}} \arrow{ddd}{\mrm{We}_{D_G^{-1}}}\\
\& X_G(N,\Gamma_0(m)) \arrow{ul}{D_{d,t}} \arrow{r}{\overline{\tau}_m}\arrow{d}{\det} \& (X(N,\Gamma_0(m))_k)_{D_G^{-1}} \arrow{d}{\mrm{We}_{D_G^{-1}}}\arrow{ur}{(D_{d,t})_{D_G^{-1}}}\&\\
\& \mrm{Pol}(G) \arrow{ld}{\cdot d} \arrow{r}{\eta} \&  (\mu_N^{\times})_{(\det{D_G})^{-1}}\arrow{rd}{\circ \underline{d}}\& \\
\mrm{Pol}(G) \arrow{rrr}{\eta} \& \&\& (\mu_N^{\times})_{(\det{D_G})^{-1}}
\end{tikzcd}
\vspace{-10pt}\]

\end{enumerate}
}

\demo{Because $(\mu_N^{\times})_k$ is a finite \'etale $(\Z/N\Z)^{\times}$-torsor over $k$ and $(\Z/N\Z)^{\times}$ acts by isomorphisms of torsors, the twist $(\mu_N^{\times})_{(\det{D_G})^{-1}}$ is well-defined and a finite \'etale $(\Z/N\Z)^{\times}$-torsor. The rule $\lambda \in \mrm{Pol}(G)(k_s) \mapsto \lambda(A,B) \in (\mu_N^{\times})_{(\det{D_G})^{-1}}(k_s)$ is a morphism of $(\Z/N\Z)^{\times} \times \mrm{Gal}(k_s/k)$-sets, so $\eta$ is a well-defined isomorphism of $(\Z/N\Z)^{\times}$-torsors by Proposition \ref{tannakian-group}.  

To check that the diagram \ref{gtw-1} commutes, since the smooth non-cuspidal locus of $X_G(N,\Gamma_0(m))_k$ is scheme-theoretically dense (and $(\mu_N^{\times})_{(\det{D_G})^{-1}}$ is separated), it is enough to check on smooth non-cuspidal $k_s$-points (e.g. \cite[$\text{IV}_4$ (17.15.10) (iii)]{EGA}), in which case the claim is clear by the explicit descriptions of all four maps.  

For \ref{gtw-2}: the bottom cell commutes since $\eta$ is an isomorphism of torsors, the top cell commutes by Proposition \ref{group-is-twist-unpolarized}, the outer and inner cells commute by \ref{gtw-1}, the leftmost cell commutes by Proposition \ref{degen-AL-Weil} and the rightmost cell commutes by Propositions \ref{degen-AL-Weil}, \ref{degen-AL-GL2} and \ref{twist-equiv}. }

\lem[group-is-twist-polarized]{Keep the set-up of Proposition \ref{group-is-twist-withpol}. Let $\lambda \in \mathrm{Pol}(G)(k)$ and let $\zeta_k: \OO_N \rightarrow k_s$ (recall that $\OO_N=\Z[1/N,\zeta_N]$) be the morphism defined by $\lambda(A,B) \in \mu_N^{\times}(k_s)$. Then
\begin{enumerate}[noitemsep,label=(\roman*)]
\item\label{gtp-1} $\det{D_G}$ is the cyclotomic character mod $N$,
\item\label{gtp-2} the $k$-scheme $\left[\Sp{\OO_N \otimes k}\right]_{(\det{D_G})^{-1}}$ is the disjoint union of copies of $\Sp{k}$ indexed by $\mrm{Hom}(\OO_N,k_s)$ as follows: a morphism $\OO_N \rightarrow k_s$ is a $x \in \mu_N^{\times}(k_s) \simeq \left[\Sp{\OO_N \otimes k}\right]_{(\det{D_G})^{-1}}(k_s)$, where the latter identification comes from the Galois twist. 	
\item\label{gtp-3} the following diagram commutes:
\vspace{-8pt}\[
\begin{tikzcd}[ampersand replacement=\&]
(\Z/N\Z)^{\times} \arrow{rr}{a \longmapsto \zeta_k\circ \underline{a}}\arrow{d}{a \mapsto a \cdot \lambda} \& \& \mrm{Hom}(\OO_N,k_s) \arrow{d}{\text{\ref{gtp-2}}}\\
\mrm{Pol}(G)(k) \arrow{rr}{\lambda' \longmapsto \lambda'(A,B)} \& \& \left[\Sp{\OO_N \otimes k}\right]_{(\det{D_G})^{-1}}(k)
\end{tikzcd}
\vspace{-12pt}\]
\end{enumerate}}

\demo{Since $\lambda$ is surjective, $\lambda(A,B)$ is a primitive $N$-th root of unity and $\zeta_k$ is well-defined. For $\sigma \in \mrm{Gal}(k_s/k)$, one has $\sigma(\lambda(A,B))=\det{D_G(\sigma)} \cdot \lambda(A,B)$, which then implies \ref{gtp-1}. One then deduces \ref{gtp-2} and \ref{gtp-3} by a direct calculation on the Galois sets (using Proposition \ref{tannakian-group}). } 

\cor{In the situation of Lemma \ref{group-is-twist-polarized}, let $n \in (\Z/n\Z)^{\times}$ and let $\underline{n} \in \mrm{Aut}(\OO_N)$ such that $\underline{n}(\zeta)=\zeta^n$ for all $N$-th roots of unity $\zeta \in \OO_N$. Then $\overline{\tau}_m$ induces an isomorphism
\vspace{-7pt}\[X_G^{n \cdot \lambda}(N,\gamma_0(m)) \rightarrow (X(N,\Gamma_0(m))_k)_{(\det{D_G})^{-1}}\times_{\left[\Sp{\OO_N \otimes k}\right]_{(\det{D_G})^{-1}}} \Sp{k},\vspace{-7pt}\]
where the second projection is given through Lemma \ref{group-is-twist-polarized} \ref{gtp-2} by the morphism $\zeta_k \circ \underline{n}: \OO_N \rightarrow k_s$.   
}

\appendix

\section{Descent and Galois twists}
\subsection{Sheaves in the fpqc topology}
Recall the following definition \cite[Definition 022B]{Stacks}: let $T$ be a scheme, and $(\pi_i)_{i \in I}$ be a family of flat scheme morphisms to $T$, with $\pi_i: T_i \rar T$. This family is a \emph{fpqc covering} if for every affine open subscheme $U \subset T$, there is a finite set $F$, a map $j: F \rar I$, and affine open subsets $W_f \subset T_{j(f)}$ for each $f \in F$, such that $U=\bigcup_{f \in F}{\pi_{j(f)}(W_f)}$. 

The following descent result is classical (e.g. \cite[\S 8.1 Proposition 1]{BLR}). 

\prop[morsheaf]{Let $T \rightarrow S$ be a morphism of schemes. The functor $V \in \Sch_S \longmapsto \mrm{Mor}_S(V,T)$ satisfies the sheaf property for fpqc covers.}

\cor[bi-morsheaf]{Let $S$ be a scheme, and $U,V$ be two $S$-schemes. The functors 
\vspace{-6pt}\[\underline{\mrm{Mor}}_S(U,V): T \in \Sch_S \longmapsto \mrm{Mor}_T(U_T,V_T),\qquad  \underline{\mrm{Isom}}_S(U,V): T \in \Sch_S \longmapsto \mrm{Isom}_T(U_T,V_T)\vspace{-6pt}\] satisfy the sheaf property for fpqc covers. }

\demo{The functor $\underline{\mrm{Mor}}_S(U,V)$ identifies with $T \in \Sch_S \mapsto \mrm{Mor}_S(U_T,V)$. Since the rule $T \in \Sch_S \mapsto U_T \in \Sch_S$ sends a fpqc cover to a fpqc cover, we are done. Since $\underline{\mrm{Isom}}_S(U,V))$ is a subfunctor of $\underline{\mrm{Mor}}_S(U,V)$, it is separated for fpqc covers. Now, let $(\pi_i)_{i \in I}$ be a fpqc cover of $T$, with $\pi_i: T_i \rightarrow T$, and let $\alpha \in \underline{\mrm{Mor}}_S(U,V)(T)$ be such that for each $i \in I$, $\pi^{\ast}\alpha \in \underline{\mrm{Isom}}_S(U,V)(T_i)$. Then, for each $i$, $\alpha \times_{V_T} V_{T_i}: U \times_{V_T} V_{T_i} \rightarrow V_{T_i}$ is an isomorphism. Since the $V_{T_i} \rar V_T$ form a fpqc cover, $\alpha$ is an isomorphism by e.g. \cite[Lemma 02L4]{Stacks} and we are done.}

\prop[group-morsheaf]{Let $S$ be a scheme, and $U,V$ be two $S$-group schemes. The functors 
\vspace{-6pt}
\begin{align*}
\underline{\mrm{Mor}}_{S-\Gp}(U,V)&: T \in \Sch_S \longmapsto \mrm{Mor}_{\GpSch_T}(U_T,V_T),\\
\underline{\mrm{Isom}}_{S-\Gp}(U,V)&: T \in \Sch_S \longmapsto \mrm{Isom}_{\GpSch_T}(U_T,V_T)
\end{align*}
\vspace{-20pt}

\noindent satisfy the sheaf property for fpqc covers.}

\demo{Since $\underline{\mrm{Isom}}_{S-\Gp}(U,V)$ is the intersection in the sheaf $\underline{\mrm{Mor}}_S(U,V)$ of the subsheaf $\underline{\mrm{Isom}}_S(U,V)$ and the sub-presheaf $\underline{\mrm{Mor}}_{S-\Gp}(U,V)$, it is enough to show that $\underline{\mrm{Mor}}_{S-\Gp}(U,V)$ is a subsheaf of $\underline{\mrm{Mor}}_S(U,V)$. Indeed, $\underline{\mrm{Mor}}_{S-\Gp}(U,V)$ is the equalizer of the two morphisms 
\vspace{-7pt}\[F_1, F_2:\underline{\mrm{Mor}}_S(U,V) \Rightarrow \underline{\mrm{Mor}}_S(U \times_S U,V),\qquad F_1(\varphi) = m_V \circ (\varphi \times \varphi); \quad F_2(\varphi) = \varphi \circ m_U, \vspace{-7pt}\]
 where $m_U, m_V$ are the group laws in $U,V$ respectively.}

\prop[desc-affine-schemes]{Let $S$ be a scheme, and $T \rar S$ be a fpqc cover. Let $F$ be a fpqc sheaf on $\Sch_S$ such that its restriction to the category $\Sch_T$ is representable by an affine $T$-scheme\footnote{That is, the map $X \rar T$ is affine.} $X$. Then $F$ is representable by an affine $S$-scheme $X'$, and $X \simeq X' \times_S T$.}

\demo{ 

\emph{Step 1: $S=\Sp{A}$} \\ After refining $T$, we may assume that $T=\Sp{B}$ and $X$ is an affine scheme, where $B$ is a faithfully flat $A$-algebra. By Yoneda's lemma, there exists $\alpha \in F(X)$ satisfying the following property: for any $T$-scheme $Z$, the map $f \in \mrm{Mor}_T(Z,X) \mapsto F(f)\alpha \in F(Z)$ is bijective. Let $\sigma:X \rightarrow T$ denote the structure map.

Let $Y$ be a $T \times_S T$-scheme and $p_1, p_2: Y \rightarrow T$ be the two projections of the structure map. One has a series of natural isomorphisms
\vspace{-7pt}\[\mrm{Mor}_{T \times_S T}(Y, X \times_S T) \simeq \mrm{Mor}_{p_1}(Y,X) \simeq F(Y) \simeq \mrm{Mor}_{p_2}(Y,X) \simeq \mrm{Mor}_{T \times_S T}(Y, T \times_S X), \vspace{-7pt}\]  
which defines by Yoneda's lemma an isomorphism $c: T \times_S X \rightarrow X \times_S T$ of $T \times_S T$-schemes. Hence the second component of $c$ is the composition of the second projection and $\sigma$, and let $f$ denote the first component of $c$. Let us show that the two morphisms $f \circ (\mrm{id}_T \times f), f \circ \mrm{pr}_{13}: T \times_S T \times_S X \rightarrow X$ agree. One has by construction $F(f)\alpha = F(\mrm{pr}_2: T \times_S X \rightarrow X)\alpha$, hence
\vspace{-10pt}
\begin{align*}%
F(f \circ \mrm{pr}_{13})\alpha&=F(\mrm{pr}_2 \circ \mrm{pr}_{13})F(f)\alpha=F(\mrm{pr}_3: T \times_S T \times_S X \rar X)\alpha,\\
F(f \circ (\mrm{id}_T \times f))\alpha&=F(\mrm{pr}_2 \circ (\mrm{id}_T \times f))\alpha=F(f\circ \mrm{pr}_{23})\alpha=F(\mrm{pr}_2 \circ \mrm{pr}_{23})\alpha\\
&=F(\mrm{pr}_3)\alpha.
\end{align*}\vspace{-20pt}

Viewing $T \times_S T \times_S X$ as a $T$-scheme by the first projection, $f \circ (\mrm{id}_T \times f), f \circ \mrm{pr}_{13}$ are morphisms of $T$-schemes, hence they are equal. This implies in particular that the following two morphisms of $T \times_S T \times_S T$-schemes agree:  
\vspace{-7pt}\[T \times_S T \times_S X \overset{\mrm{id} \times_S c}{\longrightarrow} T \times_S X \times_S T \overset{c \times \mrm{id}}{\longrightarrow} X \times_S T \times_S T,\vspace{-7pt}\]
\vspace{-14pt}\[\pi'_{13} = (f \circ \mrm{pr}_{13},\mrm{id},\sigma \circ \mrm{pr}_{3}): T \times_S T \times_S X \rar X \times_S T \times_S T.\]

By \cite[Proposition 023N]{Stacks}\footnote{This proposition is formulated for modules over rings, but the explicit construction shows that it works for algebras as long as the descent datum is also a ring homomorphism, hence it works for affine schemes as well.}, there is an affine $S$-scheme $Y$ and a $T$-isomorphism $(\chi,\sigma): X \rar Y \times_S T$ such that the following diagram of $T \times_S T$-schemes commutes:
\vspace{-7pt}\[
\begin{tikzcd}[ampersand replacement=\&]
T \times_S X \arrow{r}{c}\arrow{d}{\mrm{id}\times (\chi,\sigma)} \& X \times_S T \arrow{d}{(\sigma,\chi) \times \mrm{id}}\\
T \times_S (Y \times_S T) \arrow{r}{\mrm{can}} \&  (T \times_S Y) \times_S T
\end{tikzcd}
\vspace{-7pt}\]
Let $\alpha'=F((\chi,\sigma)^{-1})\alpha \in F(Y \times_S T)$; it follows from the above diagram that the pull-backs of $\alpha'$ under both projections $Y \times_S T \times_S T \rar Y \times_S T$ agree, hence $\alpha'$ is the pull-back of some $\alpha_0 \in F(Y)$. Thus $f \in \mrm{Mor}_S(Z,Y) \mapsto F(f)\alpha_0 \in F(Z)$ defines a morphism $u$ of fpqc sheaves on $\Sch_S$, with $u(Z)$ being a bijection whenever $Z$ has a map to $T$. Hence $u$ is locally an isomorphism, hence $u$ is an isomorphism.

\emph{Step 2: Glueing affine data} 

\noindent
By Step 1, for every affine open subscheme $U \subset S$, there exists an affine morphism $\pi_U: X'_U \rightarrow U$ and $\alpha_U \in F(X'_U)$ such that for every $U$-scheme $W$, $f \in \mrm{Mor}_U(W,X'_U) \longmapsto F(f)\alpha_U \in F(W)$ is a bijection. Moreover, for every inclusion $V \subset U$ of affine open subschemes of $S$, there is a unique map $j_{VU}: X'_V \rar X'_U$ of $U$-schemes such that $F(j_{VU})\alpha_U=\alpha_V$. Now, since $X'_U \times_U V$ (along with the pull-back of $\alpha_U$) represents the restriction of $F$ to $\Sch_V$, $j_{VU}$ is an open immersion with image $\pi_U^{-1}(V)$.

Finally, given affine open subsets $W \subset V \subset U$ of $S$, one has $j_{VU} \circ j_{WV} = j_{WU}$ (they are $U$-morphisms to $X'_U$ under which $\alpha_U$ has the same pull-back). Thus, for affine open subschemes $U,V$ of $S$, we may define an isomorphism $t_{UV}: \pi_U^{-1}(U \cap V) \rightarrow \pi_V^{-1}(V \cap U)$ of $U \cap V$-schemes by the following rule: for an affine $W \subset U \cap V$, $(t_{UV})_{|W}$ is the composition $\pi_U^{-1}(W) \underset{j_{WU}}{\simeq} X'_W \underset{j_{WV}}{\simeq} \pi_V^{-1}(W)$. Then $((X'_U)_U,(\pi_U^{-1}(U \cap V))_{U,V},(t_{U,V})_{U,V})$ is a glueing data in the sense of \cite[Section 01JA]{Stacks}, so the $X'_U$ glue to an $S$-scheme $X'$ by \emph{loc. cit.}, and the $\alpha_U$ glue to some $\alpha \in F(X')$. 

The rule $\phi: \underline{\mrm{Mor}}_{S}(-,X') \rar F, f \mapsto F(f)\alpha$ is a morphism of fpqc sheaves on $\Sch_S$, and $\phi(Z)$ is an isomorphism for every $S$-scheme admitting a map to an affine open $U \subset S$. Thus $\phi$ is a local isomorphism, hence an isomorphism and $(X',\alpha)$ represents $F$. In particular, the restriction of $F$ to $\Sch_T$ is represented by $X'_T$, so that $X'_T \simeq X$. The structure map $X' \rightarrow S$ is Zariski-locally affine over the base, hence affine. 
}

\bigskip

We conclude this section by recalling that many properties of morphisms are fpqc-local on the base. 

\prop[fpqc-descent-prop]{Let $\mathcal{P}$ be one of the following properties of morphisms of schemes: isomorphism, surjective, affine, quasi-compact, closed immersion, open immersion, locally of finite type, locally of finite presentation, locally quasi-finite, finite, flat, \'etale, smooth of relative dimension $d$, Cohen--Macaulay of relative dimension $d$, proper. 

Let $f: X \rar Y$ be a morphism of schemes and assume that there is a fpqc cover $Y_i \rar Y$ such that the base change of $f$ by any $Y_i \rar Y$ satisfies property $\mathcal{P}$. Then $f$ satisfies property $\mathcal{P}$.}

\demo{See \cite[Section 02YJ]{Stacks}.}

\subsection{Galois twists of schemes}

In this section, $k$ denotes a field with separable closure $k_s$. 

First, recall the following well-known results of Galois descent.

\lem[galois-semi-lin]{Let $U$ be a $k$-scheme and $L/k$ be finite Galois. Then the scheme $X_{L \otimes_k L}$ is the disjoint union of copies of $X_L$ indexed by $\mrm{Gal}(L/k)$, where the $\sigma$-th copy of $X_L$ is sent into $X_{L \otimes_k L}$ by the closed open immersion $(\mrm{pr}_1,\mrm{pr}_2,\underline{\sigma}\circ \mrm{pr}_2)$, where $\underline{\sigma}$ is the automorphism of $\Sp{L}$ induced by $\sigma$.}

\demo{It is enough to show that the morphism 
\vspace{-12pt}\[\phi: a \otimes b \in L \otimes_k L \mapsto (a\sigma(b))_{\sigma} \in \prod_{\sigma \in \mrm{Gal}(L/k)}{L}\vspace{-8pt}\] 
of $L$-algebras is an isomorphism. It is enough to show that $\phi$ is injective, since both sides have dimension $[L:k]$ as $L$-vector spaces. This follows in turn from the well-known fact (e.g. \cite[Chap. V, \S 6, no. 1]{Bourbaki}) that the functions $L\rightarrow L$ induced by the $\sigma \in \mrm{Gal}(L/k)$ are $L$-linearly independent. }

\lem[galois-descent-cor]{Let $L/k$ be a finite Galois extension, and $U,V$ be two $k$-schemes. The group $\mrm{Gal}(L/k)$ acts naturally on $\mrm{Mor}_L(U_L, V_L)$. The base change map 
\vspace{-7pt}\[f \in \mrm{Mor}_k(U,V) \longmapsto f \times \mrm{id}_L \in \mrm{Mor}_L(U_L,V_L)^{\mrm{Gal}(L/k)}\vspace{-7pt}\] is bijective. }

\demo{The map is obviously well-defined. It is injective since $k \rightarrow L$ is faithfully flat (by Proposition \ref{bi-morsheaf}). By \emph{loc.cit.}, it suffices to show that, if $f \in \mrm{Mor}_L(U_L,V_L)^{\mrm{Gal}(L/k)}$ and $\pi_1,\pi_2: \Sp{L} \times_k \Sp{L} \rar \Sp{L}$ are the projections, then $\pi_1^{\ast}(f)=\pi_2^{\ast}(f)$. This follows from the decomposition of Lemma \ref{galois-semi-lin}. }

\cor[inf-galois-desc]{Let $L/k$ an algebraic Galois extension. Let $U$ be a quasi-compact quasi-separated $k$-scheme and $V$ be a $k$-scheme locally of finite type. The group $\mrm{Gal}(L/k)$ acts naturally on $\mrm{Mor}_L(U_L, V_L)$. The base change map 
\vspace{-6pt}\[f \in \mrm{Mor}_k(U,V) \longmapsto f \times \mrm{id}_L \in \mrm{Mor}_L(U_L,V_L)^{\mrm{Gal}(L/k)}\vspace{-12pt}\]
 is an isomorphism.}

\demo{By \cite[$\text{IV}_3$ (8.8.2) (i)]{EGA}, the base change map $\underset{\longrightarrow}{\lim}\;\mrm{Mor}_{L'}(U_{L'},V_{L'}) \rightarrow \mrm{Mor}_L(U_L,V_L) $ is bijective (all the partial maps being injective by Proposition \ref{bi-morsheaf}), where $L'$ runs over finite Galois subextensions of $L/k$. Thus we are reduced to the case where $L/k$ is finite Galois, which is dealt with in Lemma \ref{galois-descent-cor}.}

\prop[cocycle-twist]{Let $k$ be a field with separable closure $k_s$. For $\sigma \in \mrm{Gal}(k_s/k)$, denote by $\underline{\sigma}$ the associated automorphism of $\Sp{k_s}$, so that $\mrm{Gal}(k_s/k)$ acts on the right on $\Sp{k_s}$. Let $X$ be a quasi-projective $k$-scheme. 
Let $\rho: \mrm{Gal}(k_s/k) \rar \mrm{Aut}_{k_s}(X_{k_s})$ be a continuous map satisfying the following cocycle condition: for all $\sigma,\sigma' \in \mrm{Gal}(k_s/k)$, 
\vspace{-8pt}\[\rho(\sigma\sigma')=\rho(\sigma)\circ (\mrm{id},\underline{\sigma^{-1}})\circ \rho(\sigma') \circ (\mrm{id},\underline{\sigma}).\vspace{-9pt}\] 
Then, there is a couple $(X_{\rho},j_X)$, unique up to unique isomorphism, satisfying the following properties:
\begin{itemize}[noitemsep,label=$-$]
\item $X_{\rho}$ is a quasi-projective $k$-scheme,
\item $j_X: X_{\rho} \times_k \Sp{k_s} \rar X \times_k \Sp{k_s}$ is an isomorphism of $k_s$-schemes,
\item for any $\sigma \in \mrm{Gal}(k_s/k)$, the following diagram commutes:  
\vspace{-7pt}\[
\begin{tikzcd}[ampersand replacement=\&]
X_{\rho} \times_k \Sp{k_s} \arrow{r}{j_X}\arrow{d}{(\mrm{id},\underline{\sigma}^{-1})} \& X \times_k \Sp{k_s} \arrow{d}{\rho(\sigma)\circ (\mrm{id},\underline{\sigma}^{-1})}\\
X_{\rho} \times_k \Sp{k_s} \arrow{r}{j_X} \&  X \times_k \Sp{k_s}
\end{tikzcd}
\vspace{-7pt}
\]
\end{itemize}

In particular, the following diagram commutes, for any $\sigma \in \mrm{Gal}(k_s/k)$:
\vspace{-7pt}\[
\begin{tikzcd}[ampersand replacement=\&]
X_{\rho}(k_s) \arrow{r}{j_X}\arrow{d}{\sigma} \& X(k_s) \arrow{d}{\rho(\sigma) \circ \sigma}\\
X_{\rho}(k_s) \arrow{r}{j_X} \& X(k_s)
\end{tikzcd}
\vspace{-7pt}
\]
Let $\mathscr{C}$ be a property of morphisms of $k$-schemes which is \'etale-local on the base and stable by pre- and post-composition by isomorphisms (for instance, any of the properties listed in Proposition \ref{fpqc-descent-prop}). If $X \rar \Sp{k}$ satisfies $\mathscr{C}$, then so does $X_{\rho} \rar \Sp{k}$.
If $X \rar \Sp{k}$ is projective, so is $X_{\rho} \rar \Sp{k}$.  
}

\rem{When $G$ acts in fact on $X$ by $k$-automorphisms, the cocycle condition simply means that $\rho: \mrm{Gal}(k_s/k) \rar G$ is a group homomorphism.}

\demo{(cf. \cite[\S 2.6]{Borel-Serre}) If the pair $(X_{\rho},j_X)$ exists (and $X_{\rho}$ is just some $k$-scheme), then $X_{\rho}$ is automatically a quasi-projective $k$-scheme by \cite[Lemma 0BDE]{Stacks}. The last diagram is a consequence of the preceding claims. The second-to-last assertion is formal, and the last claim is \cite[Lemma 0BDG]{Stacks}. 

We now focus on uniqueness up to unique isomorphism. Assume that $(U,j_U)$ and $(V,j_V)$ are two pairs satisfying the requested conditions. Then the map $j=j_V^{-1} \circ j_U: U_{k_s} \rar V_{k_s}$ is an isomorphism of $k_s$-schemes which is fixed by $\mrm{Gal}(k_s/k)$. Since $U,V$ are quasi-compact, $j=j_1 \times \mrm{id}_{k_s}$ for some $j_1: U \rightarrow V$ by Corollary \ref{inf-galois-desc}. Since $j_1 \times \mrm{id}_{k_s}$ is an isomorphism, $j_1$ is an isomorphism. Moreover, if $u$ is an automorphism of a pair $(X_{\rho},j_X)$ satisfying the given conditions, then $u \times \mrm{id}_{k_s} = \mrm{id}_{X_{k_s}}$, so $u=\mrm{id}_X$ by Proposition \ref{bi-morsheaf}. 

Finally, we discuss existence. Because $\rho$ is continuous, $\mrm{Aut}_{k_s}(X_{k_s})$ is discrete, and $X$ is of finite type over $k$, there is a finite Galois extension $L/k$ such that $\rho$ factors as a cocycle $\mrm{Gal}(L/k) \rightarrow \mrm{Aut}_L(X_L)$. The rule $\sigma \in \mrm{Gal}(L/k) \mapsto f(\sigma)=\rho(\sigma)\circ (\mrm{id},\underline{\sigma^{-1}}) \in \mrm{Aut}_k(X_L)$ is a group homomorphism. 

For every $\sigma \in \mrm{Gal}(L/k) \backslash \{\mrm{id}\}$, $f(\sigma)$ has no fixed point (by considering its second projection). Therefore, the quotient $q: X_L \rightarrow X'$ by $\mrm{im}(f)$ exists by \cite[II (4.5.4), (4.6.6)]{EGA} and \cite[Exp. V, Prop. 1.8]{SGA1}, and $q$ is finite \'etale by e.g. \cite[Theorem A7.1.1]{KM}.

By Lemma \ref{galois-semi-lin}, there is a unique morphism $(\varphi,\pi): X_{L \otimes_k L} \rightarrow X_L \times \underline{\mrm{Gal}(L/k)}$ such that $(\varphi,\pi) \circ (\mrm{pr}_1,\mrm{pr}_2,\underline{\sigma}\circ \mrm{pr}_2)=(f(\sigma)^{-1},\sigma)$ for all $\sigma \in \mrm{Gal}(L/k)$. Then $\varphi$ is actually a morphism of $L$-schemes, where the $L$-structure on the source is the last projection. Lemma \ref{galois-semi-lin} implies that $(\varphi,\pi)$ is an isomorphism. One has 
\vspace{-8pt}\[(f(\sigma') \times\mrm{id}_L) \circ (\mrm{pr}_1,\mrm{pr}_2,\underline{\sigma}\circ \mrm{pr}_2) = (\mrm{pr}_1,\mrm{pr}_2,\underline{\sigma'\sigma}\circ \mrm{pr}_2) \circ f(\sigma')\vspace{-9pt}\]
 as morphisms $X_L \rightarrow X_{L \otimes_k L}$, hence $(\varphi,\pi) \circ (f(\sigma') \times \mrm{id}_L) = (\mrm{id}_{X_L} \times \sigma') \circ (\varphi,\pi)$. Therefore $\varphi$ is the quotient of $X_{L \otimes L}$ by $\mrm{im}(f) \times \{\mrm{id}\}$. By \cite[Exp. V, Prop. 1.9]{SGA1}, there is a $L$-isomorphism $j_{X,L}: X'_L \rightarrow X_L$ such that $j_{X,L} \circ (q \times \mrm{id}_L) = \varphi$.  
 Then the pair $(X',j_{X,L} \times_L \Sp{k_s})$ satisfies the requirements. }

\prop[twist-equiv]{Let $f: X \rightarrow Y$ be a morphism of quasi-projective $k$-schemes. Let 
\vspace{-9pt}\[\rho_X: \mrm{Gal}(k_s/k) \rightarrow \mrm{Aut}_{k_s}(X_{k_s}), \qquad \rho_Y: \mrm{Gal}(k_s/k) \rightarrow \mrm{Aut}_{k_s}(Y_{k_s})\vspace{-9pt}\]
 be cocycles in the sense of Proposition \ref{cocycle-twist} such that, for every $\sigma \in \mrm{Gal}(k_s/k)$, one has $f_{k_s} \circ \rho_X(\sigma)=\rho_Y(\sigma) \circ f_{k_s}$. There is a unique morphism $f_{\rho}: X_{\rho} \rar Y_{\rho}$ of $k$-schemes such that the following diagram commutes:
\vspace{-8pt}\[
\begin{tikzcd}[ampersand replacement=\&]
X_{\rho} \times_k \Sp{k_s} \arrow{r}{j_X}\arrow{d}{f_{\rho} \times_k \mrm{id}_{k_s}} \& X \times_k \Sp{k_s} \arrow{d}{f \times_k \mrm{id}_{k_s}}\\
Y_{\rho} \times_k \Sp{k_s} \arrow{r}{j_Y} \&  Y \times_k \Sp{k_s}
\end{tikzcd}
\vspace{-12pt}\]

Moreover, $f_{\rho}$ acquires all the properties of $f$ that are \'etale-local on the base. The construction $(X,\rho) \mapsto (X_{\rho},j_X), f \mapsto f_{\rho}$ is functorial in the obvious sense, and it preserves finite limits.}

\demo{The morphism $f_{\rho}^0 = j_Y^{-1} \circ f_{k_s} \circ j_X: (X_{\rho})_{k_s} \rightarrow (Y_{\rho})_{k_s}$ is $\mrm{Gal}(k_s/k)$-invariant, hence is the base change of a (uniquely defined) morphism $f_{\rho}: X_{\rho} \rightarrow Y_{\rho}$ by Corollary \ref{inf-galois-desc}. The proof of functoriality is similar.

Let $\mathscr{D}$ be a finite diagram of ``quasi-projective $k$-schemes with cocycles'' and $X$ be its limit in the category $\Sch_k$, which also has an associated cocycle. Then $X$ is quasi-projective. Now, 
\vspace{-7pt}\[F_{\mathscr{D}_{\rho}}: T \in \Sch_k \longmapsto \lim\limits_{\substack{\longleftarrow\\D \in \mathscr{D}}}{\mrm{Mor}_k(T,D_{\rho})},\qquad F_{X_{\rho}}: T \in \Sch_k \longmapsto \mrm{Mor}_k(T,X_{\rho})\vspace{-7pt}\]
 are fpqc sheaves by Proposition \ref{morsheaf}, and fonctoriality, yields a map $\lambda: F_{X_{\rho}} \rar F_{\mathscr{D}_{\rho}}$. 

The base change of $\lambda$ to $k_s$ identifies with the morphism $\underline{\mrm{Mor}}_{k_s}(-,X_{k_s}) \rightarrow \lim\limits_{\substack{\longleftarrow\\D \in \mathscr{D}}}{\underline{\mrm{Mor}}_{k_s}(-,D_{k_s})}$, so it is an isomorphism. Thus $\lambda$ is a local isomorphism of sheaves, hence an isomorphism. }

\cor[twist-ring-functions]{Let $X$ be a quasi-projective $k$-scheme such that $\OO(X)$ is a finitely generated $k$-algebra. Let $\rho: \mrm{Gal}(k_s/k) \rar \mrm{Aut}_{k_s}(X_{k_s})$ satisfy the cocycle condition of Proposition \ref{cocycle-twist}. Then the map $\rho': \mrm{Gal}(k_s/k) \rightarrow \mrm{Aut}_{k_s}(\Sp{\OO(X)} \times_k \Sp{k_s})$ induced by $\rho$ satisfies the cocycle condition, and $\Sp{\OO(X_{\rho})}$ identifies with the twist of $\Sp{\OO(X)}$ by $\rho'$. }

\demo{$k_s$ is flat over $k$, so $\OO(X) \otimes_k k_s \rightarrow \OO(X_{k_s})$ is an isomorphism, and one has for all $\sigma\in \mrm{Gal}(k_s/k)$ a commutative diagram
\vspace{-10pt}\[
\begin{tikzcd}[ampersand replacement=\&]
X_{\rho} \times_k \Sp{k_s} \arrow{rrr}{j_X}\arrow{ddd}{(\mrm{id},\underline{\sigma}^{-1})} \arrow{dr} \& \& \& X \times_k \Sp{k_s} \arrow{dl}\arrow[swap]{ddd}{\rho(\sigma)\circ (\mrm{id},\underline{\sigma}^{-1})}\\
\& \Sp{\OO(X_{\rho})} \times_k \Sp{k_s} \arrow{r}{j_X}\arrow[swap]{d}{(\mrm{id},\underline{\sigma}^{-1})} \& \Sp{\OO(X)} \times_k \Sp{k_s} \arrow[swap]{d}{\rho(\sigma)\circ (\mrm{id},\underline{\sigma}^{-1})}\&\\
\& \Sp{\OO(X_{\rho})} \times_k \Sp{k_s} \arrow{r}{j_X} \&  \Sp{\OO(X)} \times_k \Sp{k_s} \&\\
X_{\rho} \times_k \Sp{k_s} \arrow{rrr}{j_X}\arrow{ru} \&\&\&  X \times_k \Sp{k_s} \arrow{lu}
\end{tikzcd}
\vspace{-12pt}\]
 the unlabeled arrows being canonical morphisms $Y \rar \Sp{\OO(Y)}$, whence the conclusion.
}

\bibliographystyle{amsplain}
{\scriptsize \bibliography{rev-biblio-moduli-smf}}

\end{document}